\documentclass[11pt,oneside]{amsart}

\usepackage{paralist}
\usepackage{graphics} 
\usepackage{epsfig} 
\usepackage{xspace}

\usepackage{array}

\usepackage{appendix}
\usepackage[french,english]{babel} 
\usepackage[T1]{fontenc} 
\usepackage[utf8]{inputenc} 
\usepackage{geometry}
\usepackage{graphicx}
\usepackage{tabto}
\usepackage{amsmath}
\usepackage{calligra}
\usepackage{amsthm}
\usepackage{bbm}
\usepackage{cancel}
\usepackage{mathrsfs}
\usepackage{mathtools}
\usepackage[algo2e]{algorithm2e}
\usepackage{version}
\usepackage{algorithm}
\usepackage{arydshln}
\usepackage[math]{cellspace}

\usepackage{multicol}
\usepackage{listings}
\usepackage{tkz-euclide}

\usepackage{matlab-prettifier}

\usepackage{bm}
\usepackage{subfig}
\usepackage{amssymb}

\usepackage{amssymb}
\newcommand{\numberset}{\mathbb}

\newcommand{\R}{\numberset{R}}
\newcommand{\F}{\numberset{F}}
\newcommand{\I}{\numberset{I}}
\newcommand{\PP}{\numberset{P}}

\DeclareFontFamily{U}{matha}{\hyphenchar\font45}
\DeclareFontShape{U}{matha}{m}{n}{
	<-6> matha5 <6-7> matha6 <7-8> matha7
	<8-9> matha8 <9-10> matha9
	<10-12> matha10 <12-> matha12
}{}
\DeclareSymbolFont{matha}{U}{matha}{m}{n}
\DeclareFontFamily{U}{mathx}{\hyphenchar\font45}
\DeclareFontShape{U}{mathx}{m}{n}{
	<-6> mathx5 <6-7> mathx6 <7-8> mathx7
	<8-9> mathx8 <9-10> mathx9
	<10-12> mathx10 <12-> mathx12
}{}
\DeclareSymbolFont{mathx}{U}{mathx}{m}{n}
\DeclareMathDelimiter{\vvvert} {0}{matha}{"7E}{mathx}{"17}%
\newcommand{\norm}[1]{\left\lVert#1\right\rVert}
\newcommand{\abs}[1]{\left\lvert#1\right\rvert}

\newtheorem{defi}{Definition}
\newtheorem{pro}{Problem}

\let\div\undefined\DeclareMathOperator{\div}{div} 
\let\curl\undefined\DeclareMathOperator{\curl}{curl} 
\DeclareMathOperator*{\Grad}{\boldsymbol\nabla}
\DeclareMathOperator*{\Grads}{\underline{\boldsymbol{\varepsilon}}}
\DeclareMathOperator{\grad}{\nabla} 
\DeclareMathOperator{\grads}{\boldsymbol\nabla_s}
\newcommand{\derivative}[2]{\frac{\partial #1}{\partial #2}}
\newcommand\Huo{(H^1_0(\Omega))^d} 
\newcommand\Ldo{L^2_0(\Omega)} 
\newcommand\Hub{(H^1(\B))^d} 
\newcommand\LdB{L^2(\B)} 
\newcommand\LdO{L^2(\Omega)} 
\newcommand\Oft{\Omega^f_t} 
\newcommand\Ost{\Omega^s_t} 
\newcommand\Of{\Omega^f} 
\newcommand\Os{\Omega^s} 
\newcommand\B{\mathcal B} 
\renewcommand\u{\mathbf{u}} 
\renewcommand\a{\mathbf{a}} 
\renewcommand\v{\mathbf{v}} 
\renewcommand\d{\mathbf{d}} 
\renewcommand\c{\mathbf{c}} 
\newcommand\f{\mathbf{f}} 
\newcommand\g{\mathbf{g}} 
\newcommand\n{\mathbf{n}} 
\newcommand\w{\mathbf{w}} 
\newcommand\x{\mathbf{x}} 
\newcommand\s{\mathbf{s}} 
\renewcommand\S{\mathbf{S}} 
\newcommand\V{\mathbf{V}} 
\newcommand\X{\mathbf{X}} 
\newcommand\Xbar{\overline{\X}} 
\newcommand\Y{\mathbf{Y}} 
\newcommand\Z{\mathbf{Z}} 
\newcommand\LL{\boldsymbol{\Lambda}} 
\newcommand\ssigma{\boldsymbol\sigma} 
\newcommand\llambda{\boldsymbol\lambda} 
\newcommand\mmu{\boldsymbol\mu} 
\newcommand\ds{\mathrm{d}\s} 
\newcommand\dx{\mathrm{d}\x} 
\newcommand\dA{\mathrm{dA}} 
\newcommand\dt{\Delta t} 
\newcommand\dr{\delta\rho} 
\newcommand\T{\mathcal{T}} 
\newcommand\calS{\mathcal{S}} 
\newcommand\quadnode{\mathbf{p}}
\newcommand\quadweigth{\omega}

\newcommand\Pcal{\mathcal{P}}
\newcommand\hatT{\hat{T}}

\usepackage{listings}
\renewcommand\lg{\begin{color}{black}}
	\newcommand\gl{\end{color}}
\newcommand\blue{\begin{color}{black}}
	\newcommand\noblue{\end{color}}
\newcommand\red{\begin{color}{black}}
	\newcommand\nored{\end{color}}
\newcommand\green{\begin{color}{black}}
	\newcommand\nogreen{\end{color}}

\renewcommand{\arraystretch}{1.5}

\begin{document}
	
		
		\title[]{On the interface matrix for fluid-structure interaction problems with fictitious domain approach}
		
		%
		\author[Daniele Boffi]{Daniele Boffi}
			\address{Computer, electrical and mathematical sciences and engineering division, King Abdullah University of Science and Technology, Thuwal 23955, Saudi Arabia and Dipartimento di Matematica \textquoteleft F. Casorati\textquoteright, University of Pavia, Via Ferrata 1, 27100, Pavia, Italy}
			\email{daniele.boffi@kaust.edu.sa}
			\urladdr{kaust.edu.sa/en/study/faculty/daniele-boffi}
		
		\author{Fabio Credali}
		\address{Computer, electrical and mathematical sciences and engineering division, King Abdullah University of Science and Technology, Thuwal 23955, Saudi Arabia and Dipartimento di Matematica \textquoteleft F. Casorati\textquoteright, University of Pavia, Via Ferrata 1, 27100, Pavia, Italy}
		\email{fabio.credali@kaust.edu.sa}
		
		\author{Lucia Gastaldi}
		\address{Dipartimento di Ingegneria Civile, Architettura, Territorio, Ambiente e di Matematica, Universit\`a degli Studi di Brescia, via Branze 43, 25123, Brescia, Italy}
		\email{lucia.gastaldi@unibs.it}
		\urladdr{lucia-gastaldi.unibs.it}
		
		\subjclass{65N30, 65N12, 74F10}
		
		\begin{abstract}
			We study a recent formulation for fluid-structure interaction problems based on the use of a distributed Lagrange multiplier in the spirit of the fictitious domain approach. In this paper, we focus our attention on a crucial computational aspect regarding the interface matrix for the finite element discretization: it involves integration of functions supported on two different meshes. Several numerical tests show that accurate computation of the interface matrix has to be performed in order to ensure the optimal convergence of the method.\medskip
			
			\noindent {\bf Keywords}: fluid-structure interaction, fictitious domain, non-matching grids, mesh intersection.
		\end{abstract}
		
		\maketitle
	
	\section{Introduction}
	The study of fluid-structure interaction problems, developed mainly from the second half of the last century, embraces various scientific fields, such as Mathematics, Physics and Engineering.
	
	Methods for FSI problems can be divided into two large families. 
	
	Boundary-fitted approaches use a mesh for the fluid domain that deforms around a Lagrangian mesh for the structure; these two meshes match on the shared interface. A popular scheme of this category makes use of an arbitrary Lagrangian-Eulerian (ALE) coordinate system \blue\cite{hirt1974arbitrary,donea1977lagrangian,donea2004arbitrary,hughes1981lagrangian}\noblue. The main advantage is that the kinematic constraints are satisfied by construction, but, on the other hand, strong distortions of the mesh could be obtained.
	
	For this reason, the so-called non-boundary-fitted methods were subsequently introduced. In these cases a discretization of the solid structure is immersed in the fluid mesh. The disadvantage that arises is a reduction in accuracy near the interface between fluid and structure. \blue This family includes several methods, developed to solve problems of various kinds since there is no method that is optimal for all situations. Among them, we have the level set formulation \cite{chang1996level} for incompressible, immiscible Navier-Stokes equations separated by a free surface and the Nitsche-XFEM \cite{alauzet2016nitsche,burman2014unfitted}, introduced to study thin-walled elastic structures immersed in an incompressible fluid.
	
	The immersed boundary method \cite{peskin2002immersed,peskin1972flow} is both a mathematical formulation and a numerical scheme. Mathematically, Eulerian and Lagrangian variables are linked making use of interaction equations involving Dirac delta functions. \green The numerical scheme consists in defining the Eulerian variables on a fixed Cartesian mesh and the Lagrangian variables on the reference domain. \nogreen Finally, the fictitious domain method \cite{glowinski2001fictitious,glowinski1997lagrange} was originally introduced to study the particulate flow: indeed, the fluid domain is artificially extended to include also particles, which are considered as solid bodies. With this methodology it is not necessary to generate a new mesh at each time step: the solid body is represented making use of a reference configuration, which is mapped at each time step to the actual position.\noblue
	
	\blue To develop our model\noblue, we started from a formulation \blue based on the \noblue immersed boundary method \red \cite{pesk} \nored and we moved towards the fictitious domain approach.
	In the case of the original finite element immersed boundary method, the evolution of the structure is governed by an ordinary differential equation that involves the position of the solid and the velocity of the fluid. Starting from here, we then introduced a Lagrange multiplier with the effect that the movement of the structure is managed through a bilinear form that contains the multiplier. As part of our research, several theoretical aspects have been addressed over the years: in \cite{stat}, we showed the well-posedness of both continuous and discrete problems, in \cite{2020} the existence and the uniqueness of the solution were proved in the case of the linearized problem. These results are summarized in a unified setting in \cite{2021}. On the other hand, numerical aspects have been studied only from the point of view of time discretization in \cite{boffi2020higher}. 
	We are dealing with a monolithic scheme, based on solving the coupled system without separating the equations for fluid and structure and characterized \green by the fact that the solid and the fluid meshes are totally independent of each other. \nogreen From the computational point of view, the assembling of the interface matrix is based on integrations over the solid mesh involving also the fluid basis functions. This operation can be carried out in two different ways: the first one is based on the computation of the intersection between the two meshes, \blue with the aim of implementing a composite quadrature rule, on the other hand, the second one works directly on the solid elements without involving preliminary geometric computations\noblue.
	
	The computation of intersection between immersed meshes has been addressed in several papers only from the computational point of view as a support for the simulation of fluid-structure interactions or contact mechanics without emphasizing why such a procedure is necessary. For instance, in \cite{krause2016parallel}, the goal is to develop an efficient parallel approach for the variational transfer of information between two meshes, using the computation of the intersection as starting point; in \cite{farrell2011conservative}, the Sutherland-Hodgman algorithm for clipping polygons is used in support of Galerkin projection for conservative interpolation of discrete fields; finally, in \cite{bvrezina2017fast}, fast algorithms for the intersection of simplicial elements of different dimensions in the three-dimensional space are presented as a support for XFEM and Mortar methods.
	
	Our aim is to show how this expensive geometric procedure is necessary in order to reach the optimal convergence rate of the method. Indeed, we are going to show that skipping the computation of the mesh intersection, we cannot reach the optimal convergence even if we increase the precision of the used quadrature rule. 
	
	The structure of this article is the following one. In Section 2, we summarize the theoretical background of our fluid-structure interaction problems starting from the Navier-Stokes equation and emphasizing the fictitious domain framework. \blue In particular\noblue, we focus our attention on the case of a thick structure. Next, we briefly recall the time semi-discretization and its connection with stationary saddle-point problems that can be studied numerically making use of the finite element method. In Section 5, we describe how the coupling matrix can be assembled with the two methods mentioned above, paying also attention on algorithmic aspects. Finally, in Section 6, we report several numerical tests showing how the choice of the procedure for assembling the interface matrix affects the convergence of our method. 
	
	\section{Fictitious domain approach: problem setting}
	Let us consider the case of a fluid-structure interaction where a solid elastic body is immersed in a two or three dimensional fluid. We denote by $\Oft$ and $\Ost$ the time dependent regions in $\R^d$ (with $d=2,3$) occupied respectively by the fluid and the solid at the time instant $t$. \blue In particular\noblue, both cases of a solid of codimension zero and one are allowed in our model, even if in this paper we are going to consider only the first case. We denote by $\Omega$ the union between $\Oft$ and $\Ost$ such that it is not dependent on time, assuming also that the solid interface $\partial\Ost$ cannot touch the boundary of $\partial\Omega$, i.e. $\partial\Omega\cap\partial\Ost=\emptyset$. Moreover, the solid is represented making use of a fixed reference domain $\B$ which is mapped at each time instant to the actual domain $\Ost$ through a mapping $\X:\B\rightarrow\Ost$.
	
	The mapping $\X$ represents the motion of the solid, while the fluid configuration is represented by the velocity $\u^f$ and by the pressure $p^f$. The solid unknown is time dependent and is defined making use of a Lagrangian variable $\s\in\B$, while the fluid unknowns are functions of time and of the Eulerian variable $\x$. In this setting, each point $\x\in\Ost$ can be represented in terms of the Lagrangian variable
	\begin{equation*}
		\x = \X(\s,t)
	\end{equation*}
	and the motion of the solid body is expressed by the kinematic equation
	\begin{equation}\label{eq:kinematik_eq_solid}
		\u^s(\x,t) = \derivative{\X}{t}(\s,t) \quad\text{ for } \x = \X(\s,t).
	\end{equation}
	On the other hand, using the derivative with respect to $\s$, we can introduce the deformation gradient $\F=\grads\X$ and its determinant $\abs{\F}$; since we are going to consider the case of incompressible solid materials, we have that $\abs{\F}$ is constant in time and when the reference domain $\B$ corresponds to the initial position $\Os_0$ of the structure, we get $\abs{\F}=1$.
	
	For the fluid, we denote by $\rho_f$ the density and by $\ssigma^f$ the Cauchy stress tensor defined, using the viscosity $\nu_f>0$ and the symmetric gradient $\Grads$, as
	\begin{equation*}
		\ssigma^f = -p^f \I + \nu_f\Grads(\u^f);
	\end{equation*}
	these objects are involved in the incompressible Navier-Stokes equations which describe the fluid dynamic in $\Oft$ as follows
	\begin{equation}
		\begin{aligned}
			&\rho_f\bigg(\derivative{\u^f}{t}+\u^f\cdot\nabla\u^f\bigg)=\div\ssigma^f\\
			&\div\u^f = 0.
			\label{eq:strong_fluid}
		\end{aligned}
	\end{equation}
	
	For the solid, we assume that the material is incompressible but also viscoelastic, hence we describe it with a Cauchy stress tensor made of two contributions $\ssigma^s = \ssigma_f^s + \ssigma_s^s$: in details, $\ssigma_f^s$ is a term analogous to the fluid stress $\ssigma^f$ but defined using the pressure $p^s$ which is \blue the Lagrange multiplier associated with the incompressibility condition\noblue, while $\ssigma_s^s$ depends on the Piola-Kirchhoff elasticity stress tensor $\PP$ via the Piola transformation; we have respectively
	\begin{equation*}
		\begin{aligned}
			\ssigma_f^s &= -p^s \I + \nu_s\Grads(\u^s)\\
			\ssigma_s^s &= \abs{\F}^{-1}\PP\F^T
		\end{aligned}
	\end{equation*}
	where we denote by $\nu_s>0$ the body viscosity. Furthermore, a potential energy density $W(\F,\s,t)$ can be used for modeling the elastic part of the stress; we have that
	\begin{equation*}
		\PP(\F,\s,t) = \derivative{W}{\F}(\F,\s,t).
	\end{equation*}
	
	At this point, we have all the ingredients we need to express the equations representing the solid behavior, indeed denoting by $\rho_s$ the solid density, we can write
	\begin{equation}
		\begin{aligned}
			&\rho_s\frac{\partial^2\X}{\partial^2 t} = \div_\s\big(\abs{\F}\ssigma_f^s\F^{-T}+\PP(\F)\big) &&\quad\text{ in }\B\\
			&\div\u^s = 0 &&\quad\text{ in }\Ost.
		\end{aligned}
	\end{equation}
	
	To complete the model we need to be careful with respect to the behavior of our system on the interface $\Gamma_t$ and prescribe initial and boundary conditions. In particular, we need to enforce continuity of velocity and Cauchy stress, hence we impose the following transmission conditions
	\begin{equation}
		\begin{aligned}
			&\u^f = \u^s && \text{ on } \partial\Ost\\
			&\ssigma_f\n_f = - \big(\ssigma_f^s+\abs{\F}^{-1}\PP\F^T\big)\n_s && \text{ on } \partial\Ost
		\end{aligned}
	\end{equation}
	where $\n_s$ and $\n_f$ denote the outer unit normal vectors to $\Ost$ and $\Oft$, respectively. Finally, we choose the following conditions
	\begin{equation}
		\begin{aligned}
			&\u^f(0) = \u_0^f &&\text{ in }\Of_0\\
			&\u^s(0) = \u_0^s &&\text{ in }\Os_0\\
			&\X(0) = \X_0 &&\text{ in } \B\\
			&\u^f =0 &&\text{ on }\partial\Omega.
		\end{aligned}
	\end{equation}

	\blue We now recall some standard notations in functional analysis \cite{lions2012non} we are going to use in the following. Given a domain $D$, we denote by $L^2(D)$ the space of square integrable functions on $D$, while for the standard Sobolev spaces we use $W^{s,q}(D)$, where $s\in\R$ indicates the differentiability and $q\in[1,\infty]$ the integrability exponent. When $q=2$, the usual notation $H^s(D)$ is used. Furthermore, $L^2_0(D)$ is the subspace of zero mean valued functions, while $H^1_0(D)$ is the subset of $H^1(D)$ with zero trace on the boundary of $D$. The norm in $H^s(D)$ is denoted by $\norm{\cdot}_{s,D}$ and $(\cdot,\cdot)_D$ stands for the scalar product of $L^2(D)$. The indication of the domain will be omitted when no confusion arises; in particular we will omit $\Omega$, while we will indicate explicitly when quantities are defined on $\B$. \noblue
	
	Once we have defined the model for both fluid and structure, we can introduce the fictitious domain approach basically extending the fluid variables also in the region occupied by the solid. Hence, we call $\u$ and $p$ velocity and pressure defined on the whole $\Omega$ with
	\begin{equation}
		\u = \begin{cases}
			\u^f \quad\text{ in }\Oft\\
			\u^s \quad\text{ in }\Ost
		\end{cases}
		\text{, }\quad
		p = \begin{cases}
			p^f \quad\text{ in }\Oft\\
			p^s \quad\text{ in }\Ost
		\end{cases}\text{.}
	\end{equation}
	Moreover, we take into account that the material velocity of the solid is equal to the velocity of the fictitious fluid, that is
	\begin{equation}\label{eq:kinematic_fictitious}
		\derivative{\X}{t}(\s,t) = \u(\X(\s,t),t) \quad \text{ for } \s\in\B;
	\end{equation}
	\blue this condition, analogous to \eqref{eq:kinematik_eq_solid}, is now a constraint for the velocity $\u$, describing its behavior in $\Ost$. \noblue
	In variational terms, condition \eqref{eq:kinematic_fictitious} is going to be enforced introducing a distributed Lagrange multiplier. For this purpose, we consider a Hilbert space $\LL$ and we define a continuous bilinear form  $\c:\LL\times\Hub\rightarrow\R$ satisfying the property
	\begin{equation}
		\c(\mmu,\Z) = 0 \quad\forall\mmu\in\LL \quad\Rightarrow\quad \Z=0.
	\end{equation}
	There are several choices for the form $\c$; in particular, we can choose $\c$ to be the duality pairing between $\Hub$ and its dual $\LL=(\Hub)^\prime$
	\begin{equation}\label{eq:c_dual}
		\c(\mmu,\Y) = \langle \mmu,\Y \rangle \quad \forall\mmu\in(\Hub)^\prime,\;\forall\Y\in\Hub
	\end{equation}
	or, alternatively, we can set $\LL=\Hub$ and define $\c$ as the scalar product
	\begin{equation}\label{eq:c_h1}
		\c(\mmu,\Y) = (\mmu,\Y)_\B + (\grad_s\mmu,\grad_s\Y)_\B \quad \forall\mmu,\Y\in\Hub.
	\end{equation}
	
	At the end, we can write the variational formulation of our problem; before doing that, we introduce the following notations
	\begin{equation*}
		\begin{aligned}
			&a(\u,\v)=\int_\Omega \nu\Grads(\u):\Grads(\v)\,\dx\quad
			&&\nu= \begin{cases}
				\nu_f & \mbox{ in } \Oft \\
				\nu_s & \mbox{ in } \Ost 
			\end{cases}\\
			&b(\u,\v,\w)=\frac{\rho_f}{2}((\u\cdot\nabla\v,\w)-(\u\cdot\nabla\w,\v))\quad
			&&\dr = \rho_s-\rho_f.
		\end{aligned}
	\end{equation*}
	
	\begin{pro}
		\label{pro:problem_fictitious}
		Given $\u_0\in\Huo$ and $\X_0:\B\longrightarrow\Omega$, $\forall t\in (0,T)$, find $(\u(t),p(t))\in\Huo \times \Ldo$, $\X(t)\in (W^{1,\infty}(\B))^d$ and $\llambda(t)\in \LL$, such that
		\begin{subequations}
			\begin{align}
					&\notag \rho_f\frac{\partial}{\partial t}(\u(t),\v)+b(\u(t),\u(t),\v)+a(\u(t),\v)\\
					&\qquad\qquad\qquad\quad-(\div\v,p(t))+\c(\llambda(t),\v(\X(t)))=0 
				&&\forall\v\in\Huo\\
				&(\div\u(t),q)=0  &&\forall q\in \Ldo\\
				& \label{eq:fict_solid} \dr \bigg(\frac{\partial^2\X}{\partial t^2},\Y \bigg)_\B+(\PP(\F(\s,t)),\grads\Y )_\B-\c(\llambda(t),\Y)=0&&\forall\Y\in \Hub \\
				& \c\bigg(\mmu,\u(\X(\cdot,t),t)-\frac{\partial\X}{\partial t}(t) \bigg)=0 &&\forall\mmu\in\LL \\
				&\u(\x,0)=\u_0(\x) &&\forall\x\in\Omega\\
				&\X(\s,0)=\X_0(\s) &&\forall\s\in\B
			\end{align}
		\end{subequations}
	\end{pro}
	
	In \cite{2020}, it was shown that a linearized version of Problem~\ref{pro:problem_fictitious} admits a unique solution. See also \cite{2021} for a review on the state of the art about this problem.
	
	\section{Time advancing scheme and associated stationary problem}
	In this section we recall the formulation of our problem when it is discretized in time via backward finite differences; in particular, we are going to see that at each time step, the solution of the problem is reduced to a stationary saddle point problem. 
	
	Let us consider a positive integer $N$ and partition the time interval $[0,T]$ into $N$ equal parts so that each node corresponds to \lg $t_n=n\dt$, for $n=0,\dots,N$, where $\dt = T/N$ \gl is the time step size. Hence, the approximation of the derivatives reads as follows
	\begin{equation}
		\begin{aligned}
			& \partial_t\v(t_{n+1})\approx\frac{\v^{n+1}-\v^{n}}{\dt}\mbox{,}&&&&
			\partial_{tt}\v(t_{n+1})\approx\frac{\v^{n+1}-2\v^{n}+\v^{n-1}}{\dt^2}\text{.}
		\end{aligned}
	\end{equation}
	
	As a consequence, the time discretized version of Problem~\ref{pro:problem_fictitious} can be written in the following manner.
	
	\begin{pro}
		\label{pro:problem_time_dis}
		Given $\u_0\in \Huo$ and $\X_0\in (W^{1,\infty}(\B))^d$, for $n=1,\dots,N$ find $(\u^n,p^n)\in\Huo\times \Ldo$, $\X^n\in \Hub$, and $\llambda^n\in\LL$, such that
		\begin{subequations}
			\begin{align}
					&\notag\rho_f\bigg(\frac{\u^{n+1}-\u^{n}}{\dt},\v\bigg)+b(\u^{n},\u^{n+1},\v)+a(\u^{n+1},\v)\\
					&\hspace*{2.7cm}-(\div\v,p^{n+1})+\c(\llambda^{n+1},\v(\X^{n}))=0&&\forall\v\in\Huo\\
						&(\div\u^{n+1},q)=0  &&\forall q\in \Ldo\\
					&\notag \dr \bigg(\frac{\X^{n+1}-2\X^n+\X^{n-1}}{\dt^2},\Y \bigg)_\B+(\PP(\F^{n+1}),\grads\Y )_\B\\
						&\hspace*{6cm}-\c(\llambda^{n+1},\Y)=0 &&\forall\Y\in \Hub\\
				& \label{eq:initiTime}\c\bigg(\mmu,\u^{n+1}(\X^{n})-\frac{\X^{n+1}-\X^{n}}{\dt} \bigg)=0  &&\forall\mmu\in\LL
			\end{align}
		\end{subequations}
	\end{pro}
	where we can define $\X^{-1}$ via an equation of the following type
	\begin{equation*}
		\frac{\X^0-\X^{-1}}{\Delta t} = \u_0^s \quad \text{ in }\B.
	\end{equation*}
	
	As we were saying before, at each time step of Problem~\ref{pro:problem_time_dis}, we solve a saddle point problem which is independent of time: we refer to \cite{stat} for its analysis, in terms of inf-sup conditions and convergence.
	
	In order to present the structure of the stationary problem \blue we are going to solve\noblue, we first reduce our discussion to the case of a linear model of the Piola-Kirchhoff stress tensor, that we define as ${\PP(\F)=\kappa\F=\kappa\grad_s\X}$. Therefore, we get the following saddle point problem.
	
	\begin{pro}
		\label{pro:stationary_general}
		\lg Let $\Xbar\in(W^{1,\infty}(\B))^d$ be invertible with Lipschitz inverse and \gl\blue $\overline{\u}\in\Huo$ such that $\div\overline{\u}=0$\noblue. Given $\f\in (\LdO)^d$, $\g\in (\LdB)^d$ and $\d\in (\LdB)^d$, find $(\u,p)\in\Huo \times \Ldo$, $\X\in \Hub$ and $\llambda\in \LL$, such that
		\begin{subequations}
			\begin{align}
				&\a_f(\u,\v)-(\div\v,p)+\c(\llambda,\v(\Xbar))=(\f,\v) && 
				\forall\v\in\Huo\\
				&(\div\u,q)=0&&\forall q\in \Ldo\\
				& \a_s(\X,\Y)-\c(\llambda,\Y)=(\g,\Y)&&\forall\Y\in \Hub \\
				& \c(\mmu,\u(\Xbar)-\X)=\c(\mmu,\d)&&\forall\mmu\in \LL
			\end{align}
		\end{subequations}
	\end{pro}
	
	Particularly, the bilinear forms $\a_f$ and $\a_s$ are defined as
	\begin{equation}
		\begin{split}
			&\a_f(\u,\v)=\alpha (\u,\v)+a(\u,\v)+b(\overline{\u},\u,\v)\quad \forall\u,\v\in(H^1_0(\Omega))^d\\
			&\a_s(\X,\Y) = \beta (\X,\Y)_\B+\gamma (\grads\X,\grads\Y)_\B \quad \forall\X,\Y\in\Hub,
		\end{split}
	\end{equation}
	where $\Xbar = \X^n$, $\overline{\u}=\u^n$ in the nonlinear terms. \blue Problem~\ref{pro:stationary_general} is obtained from Problem~\ref{pro:problem_time_dis} considering a single step and taking \noblue
	\begin{equation*}
		\begin{aligned}
			&\u = \u^{n+1},\,p = p^{n+1},\,\X = \frac{\X^{n+1}}{\lg\dt\gl},\,\llambda = \llambda^{n+1}\\
			&\f = \frac{\rho_f}{\lg\dt\gl}\u^n\\
			&\g = \frac{\dr}{\Delta t^2}(2\X^n-\X^{n-1})\\
			&\d = -\frac{1}{\Delta t}\X^n\\
			&\alpha = \frac{\rho_f}{\Delta t},\,\beta = \frac{\dr}{\Delta t},\,\gamma = \kappa\Delta t.
		\end{aligned}
	\end{equation*}
	
	\section{Finite element discretization of the stationary problem}
	
	In order to discretize Problem~\ref{pro:stationary_general} using finite elements, we need to consider four finite dimensional subspaces; in detail $\V_h\subset\Huo$, $Q_h\subset\Ldo$, $\S_h\subset\Hub$ and $\LL_h\subset\LL$, where the spaces for velocity and pressure, $\V_h$ and $Q_h$ have to satisfy the inf-sup conditions for the Stokes problem. In addition, we reduce the discussion to the case when the two spaces on the solid, $\S_h$ and $\LL_h$, coincide. 
	
	We observe that if $\c$ is defined as the duality pairing between $\Hub$ and its dual, we can identify it with the scalar product in $(\LdB)^d$ providing that $\LL_h$ is included in $(\LdB)^d$:
	\begin{equation}\label{eq:c_l2_inner}
		\c(\mmu,\Y) = (\mmu,\Y)_\B \quad \forall\mmu\in\LL_h,\forall\Y\in\S_h.
	\end{equation}
	On the other hand, if $\c$ is defined as the inner product of $\Hub$, no modification is necessary.
	
	These finite dimensional spaces are defined making use of two meshes for the spatial discretization of the domains $\Omega$ and $\B$: we denote by $\mathcal{T}_h$ the spatial discretization of $\Omega$ characterized by the size $h_{\mathcal{T}}$, while $\mathcal{S}_h$ is the mesh we use for $\B$, with size $h_{\mathcal{S}}$.
	
	In this setting, we are now able to state the finite element version of Problem~\ref{pro:stationary_general} whose well-posedness has been discussed in \cite{stat} and \cite{2021}.
	
	\begin{pro}
		\label{pro:discrte_stationary_general}
		Let $\Xbar\in(W^{1,\infty}(\B))^d$ be invertible with Lipschitz inverse and \red $\overline{\u}\in\Huo$ such that $\div\overline{\u}=0$\nored. Given $\f\in (\LdO)^d$, $\g\in (\LdB)^d$ and $\d\in (\LdB)^d$, find $(\u_h,p_h)\in\V_h \times Q_h$, $\X_h\in \S_h$ and $\llambda_h\in \LL_h$, such that
		\begin{subequations}
			\begin{align}
				&\a_f(\u_h,\v_h)-(\div\v_h,p_h)+\c(\llambda_h,\v_h(\Xbar))=(\f,\v_h) && 
				\forall\v_h\in\V_h\\
				&(\div\u_h,q_h)=0&&\forall q_h\in Q_h\\
				& \a_s(\X_h,\Y_h)-\c(\llambda_h,\Y_h)=(\g,\Y_h)&&\forall\Y_h\in \S_h \\
				& \c(\mmu_h,\u_h(\Xbar)-\X_h)=\c(\mmu_h,\d)&&\forall\mmu_h\in \LL_h.
			\end{align}
		\end{subequations}
	\end{pro}
	
	Moreover, this can be rewritten in matrix form as
	\begin{equation}
		\label{eq:matrix}
		\left[\begin{array}{@{}cc|c|c@{}}
			A_f & B^\top & 0 & C_f^\top \\ 
			B & 0 & 0 & 0 \\
			\hline
			0 & 0 & A_s & -C_s^\top \\
			\hline
			C_f & 0 & -C_s & 0\\
		\end{array}\right]
		\left[ \begin{array}{c}
			\u\\
			p\\
			\hline
			\X\\
			\hline
			\llambda
		\end{array}\right]=
		\left[ \begin{array}{c}
			\f\\
			\mathbf{0}\\
			\hline
			\g\\
			\hline
			\d
		\end{array}\right],
	\end{equation} \renewcommand{\arraystretch}{1.15}
	\blue where the blocks are defined \noblue using the basis functions of $\V_h$, $Q_h$, $\S_h$ and $\LL_h$: indeed, denoting the basis functions respectively with $\bm{\phi}$, $\psi$, $\bm{\chi}$ and $\bm{\zeta}$, we get
	\begin{equation*}
		\begin{array}{l}
			(A_f)_{ij}=\a_f(\bm{\phi}_j,\bm{\phi}_i)\\[2ex]
			B_{ki}=-(\div\bm{\phi}_i,\psi_k)\\[2ex]
			(A_s)_{ij}=\a_s(\bm{\chi}_j,\bm{\chi}_i)\\[2ex]
			(C_f)_{\ell j}=\c(\bm{\zeta}_\ell,\bm{\phi}_j(\Xbar))\\[2ex]
			(C_s)_{\ell j}=\c(\bm{\zeta}_\ell,\bm{\chi}_j)\mbox{.}
		\end{array}
	\end{equation*}

	\section{Assembling the interface matrix $C_f$}
	
	The $C_f$ block of the stiffness matrix of the finite element discretization \eqref{eq:matrix} plays a crucial role in our method because it represents the interaction between the fluid and the structure. In particular, since we use a fictitious domain approach, the fluid domain is extended in order to include also the region occupied by the solid, hence \blue there is an overlap of the fluid mesh with the solid mesh mapped into fluid domain\noblue.
	
	The matrix $C_f$ \green is associated with \nogreen the bilinear form $\c(\mmu_h,\v_h(\Xbar))$, which can be the scalar product of $(\LdB)^d$ or $\Hub$. Therefore, we have to compute an integral over the solid reference domain $\B$ that involves $\mmu_h\in\LL_h$, defined on $\B$, and $\v_h\in\V_h$ defined on the whole domain $\Omega$.
	
	\begin{figure}[!ht]
		\centering
		\includegraphics[width=0.55\textwidth]{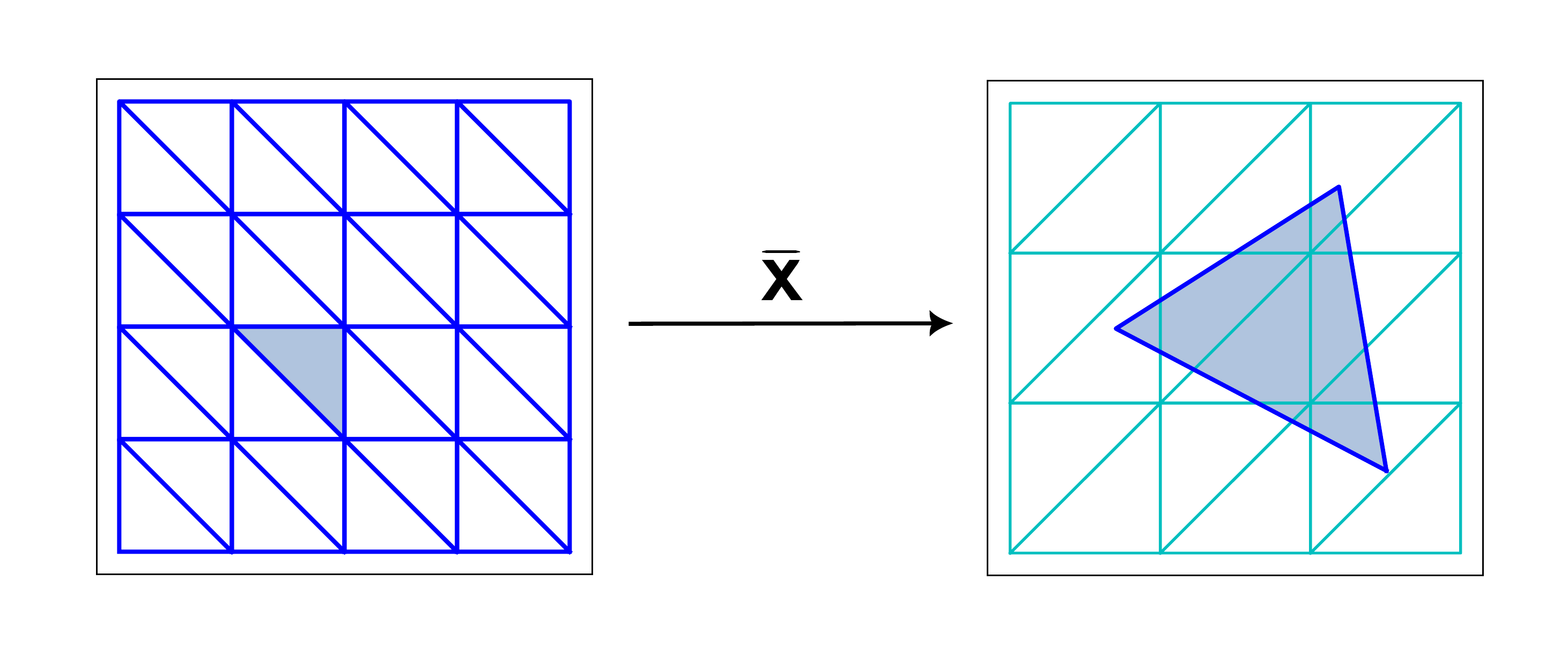}
		\caption{Mapping of a solid element into the fluid mesh.}
		\label{fig:mapping_tri}
	\end{figure}

\begin{figure}
		\centering
		\subfloat[]{ \label{fig:basis_supp}\includegraphics[width=0.31\textwidth]{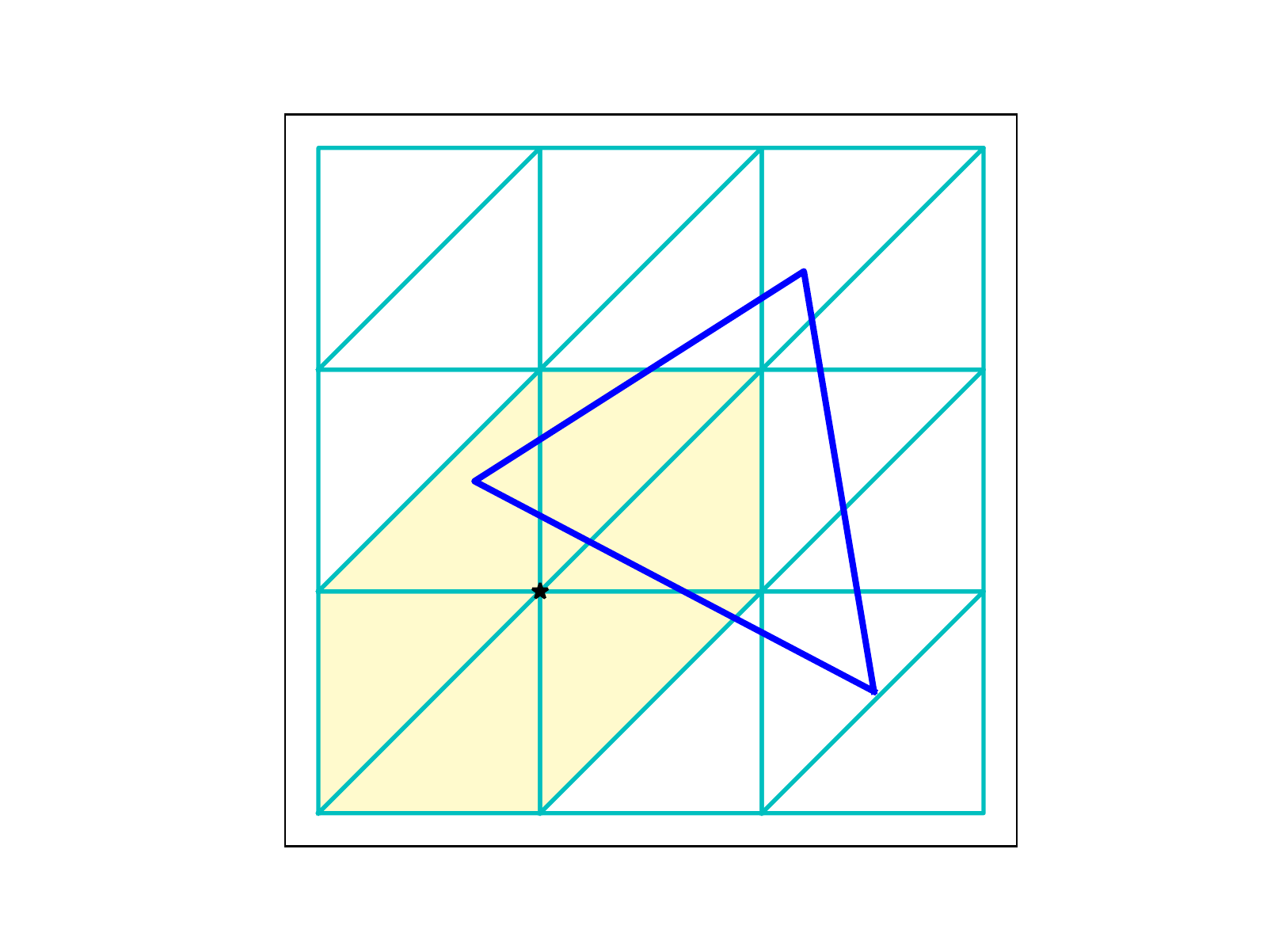}}\quad
		\subfloat[]{ \label{fig:int}\includegraphics[width=0.31\textwidth]{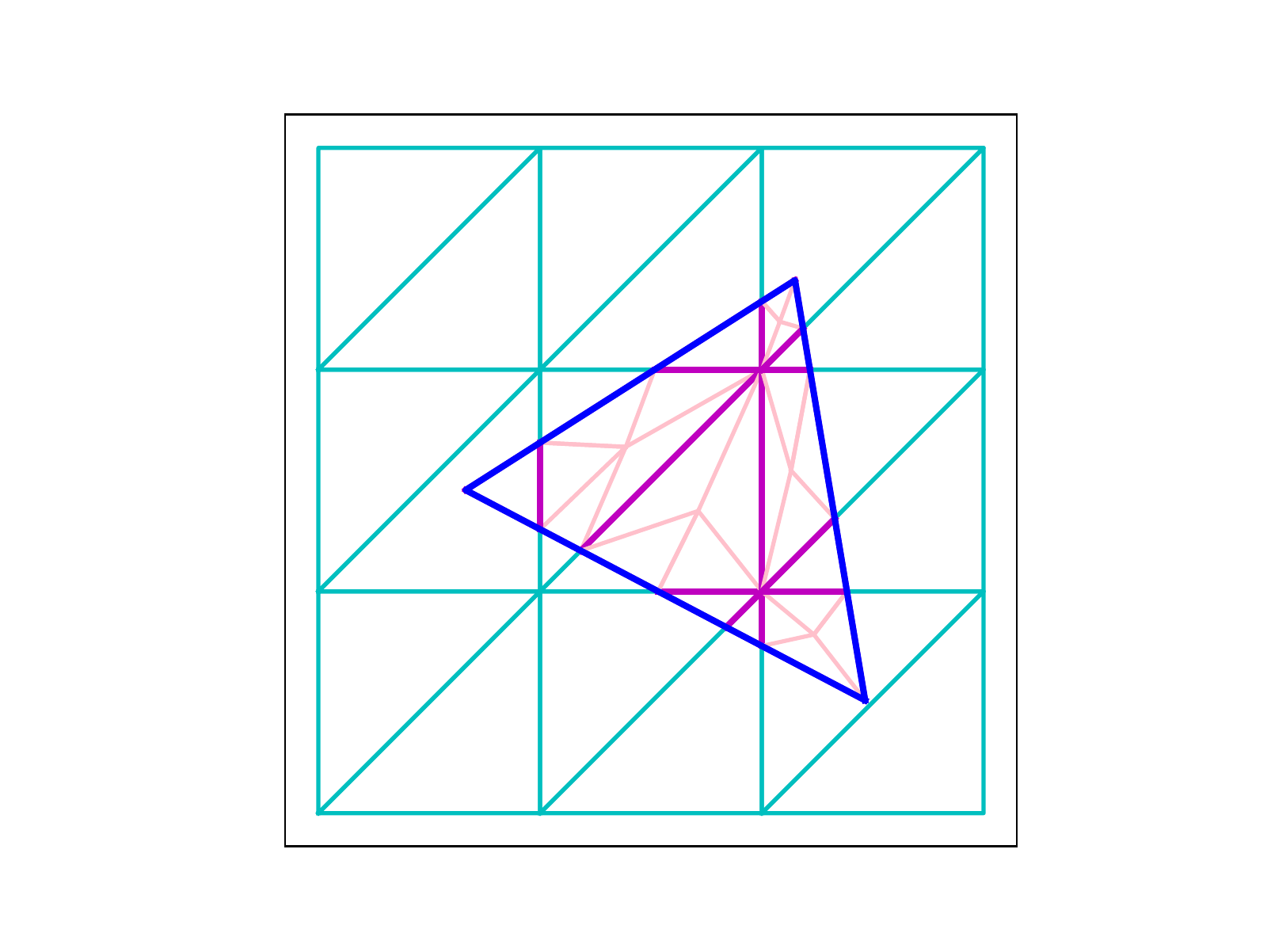}}\quad
		\subfloat[]{\label{fig:noint}\includegraphics[width=0.31\textwidth]{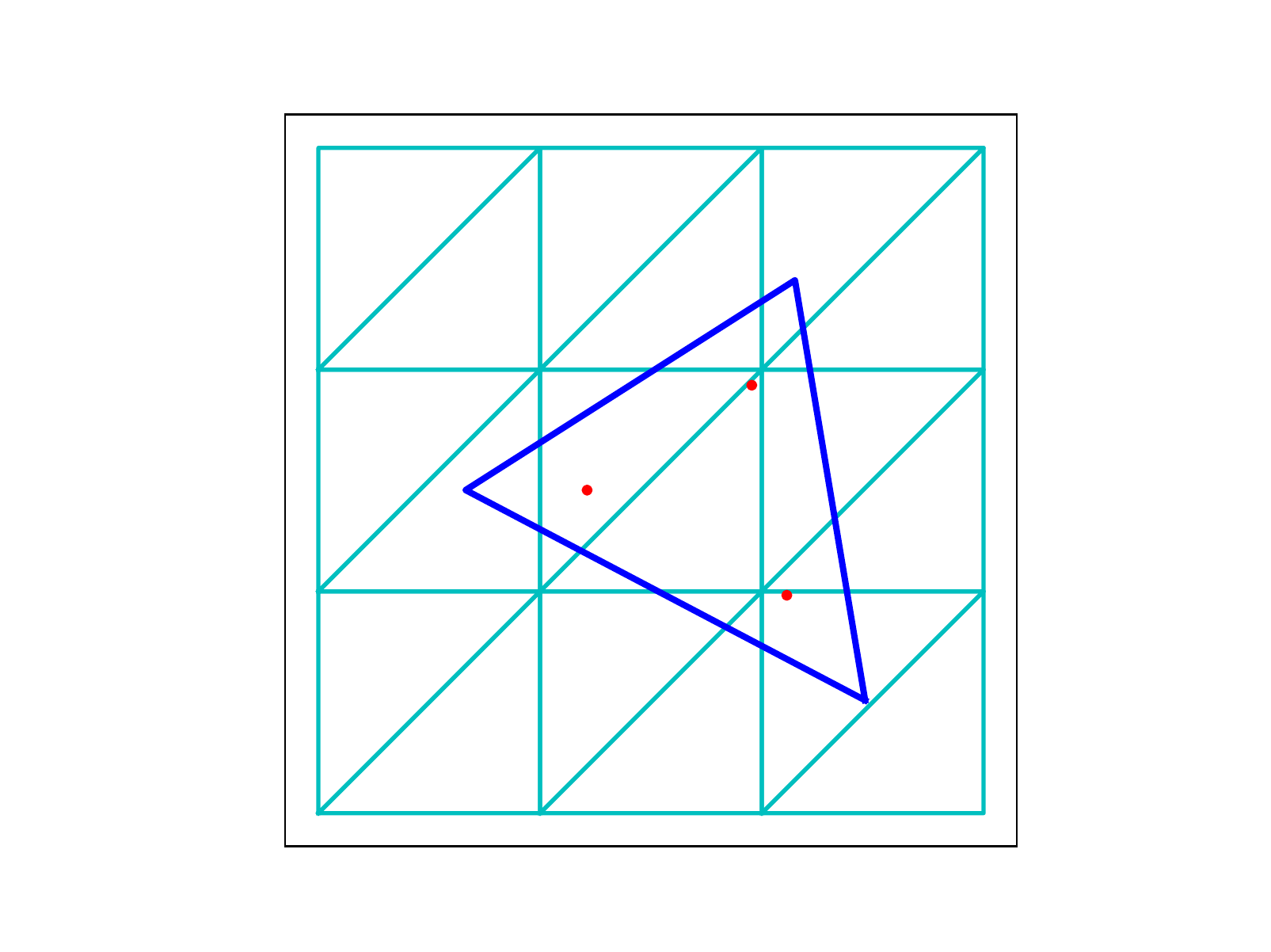}}
		\caption{A graphical example for the two approaches proposed for the assembling of the interface matrix. The cyan mesh is related to the velocity while the big triangle is a mapped element in its actual position. (a) The support of the fluid basis function associated with the black-starred point is represented in yellow. In particular, it is clear that only a portion of the support is included in the solid element. (b) The mapped solid triangle is partitioned into sub-polygons: each of them is related to a single velocity element. If a polygon is not already a triangle, it is triangulated (pink lines). (c) Mapped nodes for a Gauss quadrature rule with order two.}
	\end{figure}

	\green In Figure \ref{fig:mapping_tri}, we show a portion of the meshes $\calS_h$ of $\B$ and $\T_h$ of $\Omega$. In the picture we highlight in gray a particular triangle $T\in\calS_h$ and the corresponding element $\Xbar(T)$ immersed in $\Omega$. Basically, to map the mesh $\calS_h$ of $\B$ into $\Omega$, we simply apply the transformation $\Xbar$ to the nodes, keeping linear the edges connecting them. As a consequence, how to manage the geometrical aspects becomes really important because we have to integrate over $\calS_h$ the velocity shape functions defined on $\T_h$, combined with $\Xbar$. In doing that, we have to take into account the position occupied by $\Xbar(\calS_h)$, see Figure \ref{fig:basis_supp} for an example of mismatch between the supports of fluid and mapped solid basis functions: indeed, the support of the fluid basis function under consideration, indicated in yellow, only partially matches the blue mapped solid triangle. We also notice that this situation does not concern the assembling of $C_s$, also defined through $\c$, but involving the spaces $\S_h$ and $\LL_h$ both defined on $\B$. \nogreen
	
	In this section, we restrict our discussion to the case $d=2$ with triangular meshes in order to present two different approaches for the assembling procedure of $C_f$: the first approach depends on the computation of the intersection between the velocity mesh and the mesh of $\B$ mapped into $\Os$, while the second one skips \lg these \gl geometric computations. 
	
	At this point, we introduce the discrete version of the two choices for the bilinear form $\c$ with respect to the mesh \lg$\calS_h$\gl; \lg if $\c$ is defined like in \eqref{eq:c_l2_inner}, we have 
	\begin{equation} \label{eq:ch_l2}
			\c_h(\mmu_h,\v_h(\Xbar)) = \int_{\lg\calS_h\gl} \mmu_h\cdot\v_h(\Xbar)\,\ds 
			= \sum_{T_s\in\calS_h}\int_{T_s} \mmu_h\cdot \v_h(\Xbar)\,\ds
	\end{equation}
	otherwise, if the definition is the one in \eqref{eq:c_h1}, we can write
	\begin{equation}
		\begin{aligned}\label{eq:ch_h1}
			\c_h(\mmu_h,\v_h(\Xbar)) &= \int_{\lg\calS_h\gl} \mmu_h\cdot \v_h(\Xbar)+ \grads\mmu_h:\grads\v_h(\Xbar)\,\ds \\
			&= \sum_{T_s\in\calS_h}\int_{T_s} \mmu_h\cdot \v_h(\Xbar) + \grads\mmu_h:\grads\v_h(\Xbar)\,\ds.
		\end{aligned}
	\end{equation}

	In the next subsections, we expand these two definitions, denoting by $\{(\quadnode_k,\quadweigth_k)\}_{k=1}^K$ nodes and weights respectively of the quadrature rule under consideration; $\abs{T}$ denotes the area of a triangle $T$. \gl\blue We remark that $\grads\v_h(\Xbar)$ is the gradient of a composite function, indeed we have to compute the gradient with respect to $\s$ of the velocity $\v_h$ applied to the map $\Xbar$, defined on $\B$\noblue.
	
	\subsection{Assembling with mesh intersection}
	\label{sec:Assembling with mesh intersection}
	
	\lg If we want to perform an exact computation of the integrals in \eqref{eq:ch_l2} and \eqref{eq:ch_h1}, we need to use a composite rule on the intersection of the fluid and solid meshes. \gl
	
	The aim of computing the intersection is to obtain a new triangulation for the structure, finer than the original one, such that each new element is immersed in (i.e. interacts with) a single fluid element: in this way, all the basis functions involved for the integration (both solid and fluid) are supported in the structure element under consideration. With this procedure, we can obtain an accurate computation of the interaction because, thanks to the composite quadrature rule defined, we can take fully into account the contribution given by the objects involved.
	
	In particular, the computation of the intersection has to be made testing each actual solid element \lg $\Xbar(T_s)$ \gl with each fluid element \lg $T_f\in\T_h$ \gl; the intersection could be still a triangle or a general polygon: in the first case, we have already a new element of the finer triangulation while, in the latter one, we compute a sub-triangulation connecting the barycenter with the vertices of the polygon. \lg In particular, we denote by $\{P_j\}_{j=1}^J$ the set of the resulting polygons for each $T_s$ and by $\{T_i\}_{i=1}^{I_j}$ the elements of the new triangulation of $P_j$. Therefore, the numerical version of the integrals in \eqref{eq:ch_l2} is given by
	\begin{equation}
		\begin{aligned}
			\int_{T_s} \mmu_h\cdot \v_h(\Xbar)\,\ds
			&= \sum_{j=1}^{J} \int_{P_j} \mmu_h\cdot \v_h(\Xbar)\,\ds
			= \sum_{j=1}^{J} \sum_{i=1}^{I_j} \int_{T_i} \mmu_h\cdot \v_h(\Xbar)\,\ds\\
			&= \sum_{j=1}^{J} \sum_{i=1}^{I_j}\abs{T_i}\sum_{k=1}^{K} \quadweigth_k \mmu_h(\quadnode_k)\cdot \v_h(\Xbar(\quadnode_k)),
		\end{aligned}
	\end{equation}
	while for \eqref{eq:ch_h1}, we have
	\begin{equation}
		\begin{aligned}
			&\int_{T_s} \mmu_h\cdot \v_h(\Xbar)+ \grads\mmu_h:\grads\v_h(\Xbar)\,\ds\\
			&= \sum_{j=1}^{J} \int_{P_j} \mmu_h\cdot \v_h(\Xbar)+ \grads\mmu_h:\grads\v_h(\Xbar)\,\ds\\
			&= \sum_{j=1}^{J} \sum_{i=1}^{I_j} \int_{T_i} \mmu_h\cdot \v_h(\Xbar)+ \grads\mmu_h:\grads\v_h(\Xbar)\,\ds\\
			&= \sum_{j=1}^{J} \sum_{i=1}^{I_j} \abs{T_i} \sum_{k=1}^{K} \quadweigth_k \big[ \mmu_h(\quadnode_k)\cdot \v_h(\Xbar(\quadnode_k)) + \grads\mmu_h(\quadnode_k):\grads\v_h(\Xbar(\quadnode_k)) \big].
		\end{aligned}
	\end{equation}\gl
	
	Algorithm \ref{alg:intersection} summarizes the steps of the procedure that is also graphically represented in Figure \ref{fig:int} where purple lines denote the polygons resulting from the intersection between the blue solid element and the cyan velocity triangles; pink lines denote the triangulation constructed using the barycenter.
	
	
	\begin{algorithm}[h!]
		\caption{Assembling $C_f$ with mesh intersection}
		\begin{flushleft}
			\label{alg:intersection}
			$\lbrace T_f\rbrace_{f=1,\dots,N_F}^{} =\mbox{elements of the fluid mesh}$\\
			$\lbrace T_s\rbrace _{s=1,\dots,N_S}^{} =\mbox{elements of the reference solid mesh}$\\
			\medskip
			Compute mesh intersection $\lbrace T_f\rbrace_f\cap\lbrace \Xbar(T_s)\rbrace _s$:\\
			\quad each $\Xbar(T_s)$ is partitioned into \lg$J\ge 1$ \gl polygons $\lbrace P_{j}\rbrace _{j=1,\dots,J}$\\
			\quad associated with $J$ fluid elements $\lbrace T_{f,j}\rbrace _{j=1,\dots,J}$, i.e. $P_{j}\leftrightarrow T_{f,j}$\\
			\medskip
			\For{$\lbrace T_s\rbrace _{s=1,\dots,N_S}$}{
				
				\For{$\lbrace P_{j}\rbrace _{j=1,\dots,J}$}{
					
					\If{$P_{j}$ is not a triangle}{Triangulate $P_{j}$ into $\lbrace T_i\rbrace_{i=1,\dots,\lg I_j\gl}$}
					\Else{$I_j=1$ and $T_1 = P_{j}$}
					\For{$\lbrace T_i\rbrace_{i=1,\dots,\lg I_j\gl}$}{
						Integrate using the basis functions related to $T_s$ and $T_{f,j}$\\
						Load contribution to the global matrix
					}
					
				}
				
			}
		\end{flushleft}
	\end{algorithm}

	\subsection{Assembling without mesh intersection}
	
	A cheaper approach consists in directly integrating on each element of the structure discretization, therefore the first step consists in computing the quadrature nodes of the element under consideration accordingly to the chosen quadrature rule. At this point, the evaluation of the solid basis functions is a trivial operation.
	
	On the other hand, for the evaluation of the shape functions related to the velocity, we must pay attention because the quadrature nodes could belong to different fluid elements: for this reason, it becomes necessary to understand whose fluid triangle each of the quadrature points belongs to. In the case of triangular meshes, this operation can be performed computing the barycentric coordinates of the nodes with respect to all fluid elements. Once this is done, at each quadrature point we can evaluate the fluid functions related to the element that contains it and then integrate. We emphasize that this operation introduces an additional source of error because we are evaluating the fluid basis functions at a quadrature node assuming that at the other ones they are not supported.
	
	\lg In this case, the quadrature over $T_s$ for \eqref{eq:ch_l2} is given by
	\begin{equation}
		\int_{T_s} \mmu_h\cdot \v_h(\Xbar)\,\ds = \abs{T_s} \sum_{k=1}^{K} \quadweigth_k \mmu_h(\quadnode_k)\cdot \v_h(\Xbar(\quadnode_k));
	\end{equation}
	conversely, for \eqref{eq:ch_h1}, we have
	\begin{equation}
		\begin{aligned}
			&\int_{T_s} \mmu_h\cdot \v_h(\Xbar) + \grads\mmu_h:\grads\v_h(\Xbar)\,\ds\\
			&\quad= \abs{T_s} \sum_{k=1}^{K} \quadweigth_k \big[ \mmu_h(\quadnode_k)\cdot \v_h(\Xbar(\quadnode_k)) + \grads\mmu_h(\quadnode_k):\grads\v_h(\Xbar(\quadnode_k)) \big].
		\end{aligned}
	\end{equation}\gl
	
	A clarifying example is in Figure \ref{fig:noint}, where to fix ideas, we represent the coupling without intersection is the case of a Gaussian quadrature rule with three nodes. Here, if we look at the three points, the basis functions of the related fluid elements are supported also in fluid triangles that we are not considering for this computation. All the steps are reported in Algorithm \ref{alg:no_intersection}.
	\begin{algorithm}[h!]
		\caption{Assembling the $C_f$ matrix block without mesh intersection}
		\begin{flushleft}
			\label{alg:no_intersection}
			$\lbrace T_f\rbrace_{f=1,\dots,N_F}^{} =\mbox{elements of the fluid mesh}$\\
			$\lbrace T_s\rbrace _{s=1,\dots,N_S}^{} =\mbox{elements of the reference solid mesh}$\\
			\medskip
			\For{$\lbrace T_s\rbrace _{s=1,\dots,N_S}$}{
				\red Compute $K$ quadrature nodes $\quadnode_1,\dots,\quadnode_K$ in $T_s$ accordingly to the rule chosen\\\nored
				Evaluate solid basis functions\\
				Find the fluid element containing each mapped quadrature point: \red$\Xbar(\quadnode_k)\leftrightarrow T_{f,k}$ \nored with $k=1,\dots,K$\\
				\For{$\lbrace T_{f,k}\rbrace _{k=1,\dots,K}$}{
					Evaluate the fluid basis functions in \red$\Xbar(\quadnode_k)$ \nored and integrate over $T_s$\\
					Load contribution to the global matrix
				}
			}
		\end{flushleft}
	\end{algorithm}

	\subsection{Other geometries and three-dimensional case}
	It is important to notice that the procedures just described for two dimensional simplices may be also valid for other geometries, such as quadrilateral meshes. Moreover, we can also extend the same ideas to the three dimensional case. For instance, if we consider simplicial meshes in the three-dimensional space, we have to deal with tetrahedra. From the theoretical point of view, our method still works and the two algorithms regarding the assembling techniques for the interface matrix do not change; \lg on the other hand, the computation of the intersection of three dimensional polyhedra is a non-trivial operation and needs the support of specific geometric libraries. \gl
	For the non intersection case, the use of barycentric coordinates for detecting where the quadrature nodes are placed in the fluid mesh is still a valid procedure.
	
	\section{Numerical tests}
	
	\subsection{Model problem}
	To study how the two approaches presented above affect the convergence of our method, we perform several numerical tests on a two dimensional stationary problem, that turns out to be a simplified version of the problem presented in Problem~\ref{pro:stationary_general}.
	\begin{pro}
		\label{pro:stationary_tests}
		Let $\Xbar\in (W^{1,\infty}(\B))^2$ be invertible with Lipschitz inverse. Given $\f\in (\LdO)^2$, $\g\in (\LdB)^2$ and $\d\in (\LdB)^2$, find $(\u,p)\in(H_0^1(\Omega))^2 \times \Ldo$, $\X\in (H^1(\B))^2$ and $\llambda\in (H^1(\B))^2$, such that
		\begin{subequations}
			\begin{align}
				&\a_f(\u,\v)-(\div\v,p)+\c(\llambda,\v(\Xbar))=(\f,\v) && 
				\forall\v\in(H_0^1(\Omega))^2\\
				&(\div\u,q)=0&&\forall q\in \Ldo\\
				& \a_s(\X,\Y)-\c(\llambda,\Y)=(\g,\Y)&&\forall\Y\in (H^1(\B))^2 \\
				& \label{eq:motion_test} \c(\mmu,\u(\Xbar)-\X)=\c(\mmu,\d)&&\forall\mmu\in (H^1(\B))^2
			\end{align}
		\end{subequations}
	\end{pro}
	\noindent where the bilinear forms are defined as
	\begin{equation*}
		\begin{split}
			&\a_f(\u,\v)=(\Grad\u,\Grad\v)\\
			&\a_s(\X,\Y)=(\grads\X,\grads\Y)_\B\\
			&\c(\mmu,\Y)= (\grads\mmu,\grads\Y)_\B + (\mmu,\Y)_\B\\
		\end{split}
	\end{equation*}
	For each of our eight tests, we first select analytical solutions not belonging to the finite element spaces chosen for the approximation, and then we compute the right hand sides as follows
	\begin{equation*}
		\begin{aligned}
			\f&=-\bm{\Delta}\u+\grad p+\c(\llambda,\v)\\
			\g&=-\bm{\Delta}\X-\llambda+\g_{\partial\B}\\
			\d&=\u(\Xbar) - \X
		\end{aligned}
	\end{equation*}
	where $\g_{\partial\B}$ is the boundary term we obtain integrating by parts the equation of the solid without considering the multiplier; for this reason it is related with the normal derivative of $\X$ on $\partial\B$
	\blue\begin{equation*}
		\a_s(\X,\Y)=(\grads\X,\grads\Y)_\B = (-\bm{\Delta}\X,\Y)_\B + (\grads\X\cdot\n_s,\Y)_{\partial\B} \mbox{.}
	\end{equation*}\noblue
	For the definition of the datum $\f$, we emphasize again that the contribution given by the multiplier is weakly imposed only on the solid domain. Consequently, exploiting once again the immersion of the solid in the fluid, in order to assemble this term, we have \red still \nored to use the mesh intersection to ensure an optimal computation of the datum.
	
	The first six cases are characterized by the fact that $\Xbar$ is set to be the identity, hence the reference and the actual solid domain coincide; moreover in Tests 5 and 6 we consider the case of discontinuous pressure with the discontinuity matching both the fluid triangulation and the boundary of the structure and with non-matching interface, respectively. On the other hand, for the last two tests, we consider the case of a non trivial mapping $\Xbar$ with two different bodies. Each test is performed choosing the Stokes pair $\Pcal_1-iso-\Pcal_2/\Pcal_1$ introduced by Bercovier and Pironneau in \cite{bercovier} and its enhanced version $\Pcal_1-iso-\Pcal_2/\Pcal_1+\Pcal_0$ discussed in \cite{mass}, \blue while for the approximation of both the solid variables, we choose piecewise linear elements $\Pcal_1$\noblue. 
	
	\subsection{Bercovier-Pironneau element}\label{sec:BP}
	The Bercovier-Pironneau element is a popular choice in this field of study because it is considered as a cheaper version of the Hood-Taylor element: they share the same set of degrees of freedom, but the Bercovier-Pironneau involves linear shape functions also for the velocity discretization, hence, for this reason it is a low-order choice.
	
	The peculiarity of a choice of this type concerns with the discretization of the domain and the need to work with two different meshes for velocity and pressure: indeed, each pressure element is a macro element containing the four velocity triangles obtained by connecting the middle points. Therefore, we can write $\V_h$ and $Q_h$ as
	\begin{equation}
		\begin{aligned}
			\V_h &= \{ \v\in(H^1_0(\Omega))^2: \v_{\mid T}\in(\Pcal_1(T))^2 \quad \forall T\in\mathcal{T}_{h/2} \} \\
			Q_h &= \{q\in \Ldo: q_{\mid T}\in \Pcal_1(T)\quad\forall T\in\mathcal{T}_{h}\}.
		\end{aligned}
	\end{equation}
	
	The enhanced version with discontinuous pressure has been introduced because it guarantees a local mass conservation with the average of the divergence equal to zero in each element; in this case, the pressure space is defined as 
	\begin{equation}
		Q_h = \{ q\in \Ldo: q=q_1+q_0,\,q_1\in \blue H^1(\Omega)\noblue,\,q_{1\mid T}\in \Pcal_1(T),\,q_{0\mid T}\in \Pcal_0(T)\quad\forall T\in\mathcal{T}_{h} \}.
	\end{equation}

	The use of two different meshes for velocity and pressure has consequences also in the assembling of the interface matrix $C_f$: since it is defined from $\c(\llambda_h,\u_h(\Xbar))$, we have to couple the solid mesh $\calS_h$ with the one related to the velocity $\mathcal{T}_{h/2}$.	
	
	\subsection{Mesh generation}
	For the fluid domain discretization, we use uniform meshes. \blue In particular\noblue, in the case of $\Pcal_1-iso-\Pcal_2/\Pcal_1+\Pcal_0$ elements, we need to pay attention to corner triangles\lg, \gl i.e. triangles with two edges placed on the boundary\lg, \gl because they produce wrong approximations due to the sensitivity to the boundary conditions \lg(the phenomenon is analyzed in \cite{mass})\gl: this situation is corrected via diagonal exchange. On the other hand, for the structure, we work with both uniform and unstructured grids. 
	
	All the meshes, satisfying the shape regularity property, are generated on the unit square $[0,1]^2$ and then mapped \red into \nored $\Omega$ and $\B$. To obtain unstructured meshes, we used the Gmsh finite element mesh generator \cite{geuzaine2009gmsh}.
	
\begin{figure}[!h]
		\centering
		\subfloat[Right.1]{\label{fig:right1} \includegraphics[trim=50 16 40 50,width=0.25\linewidth]{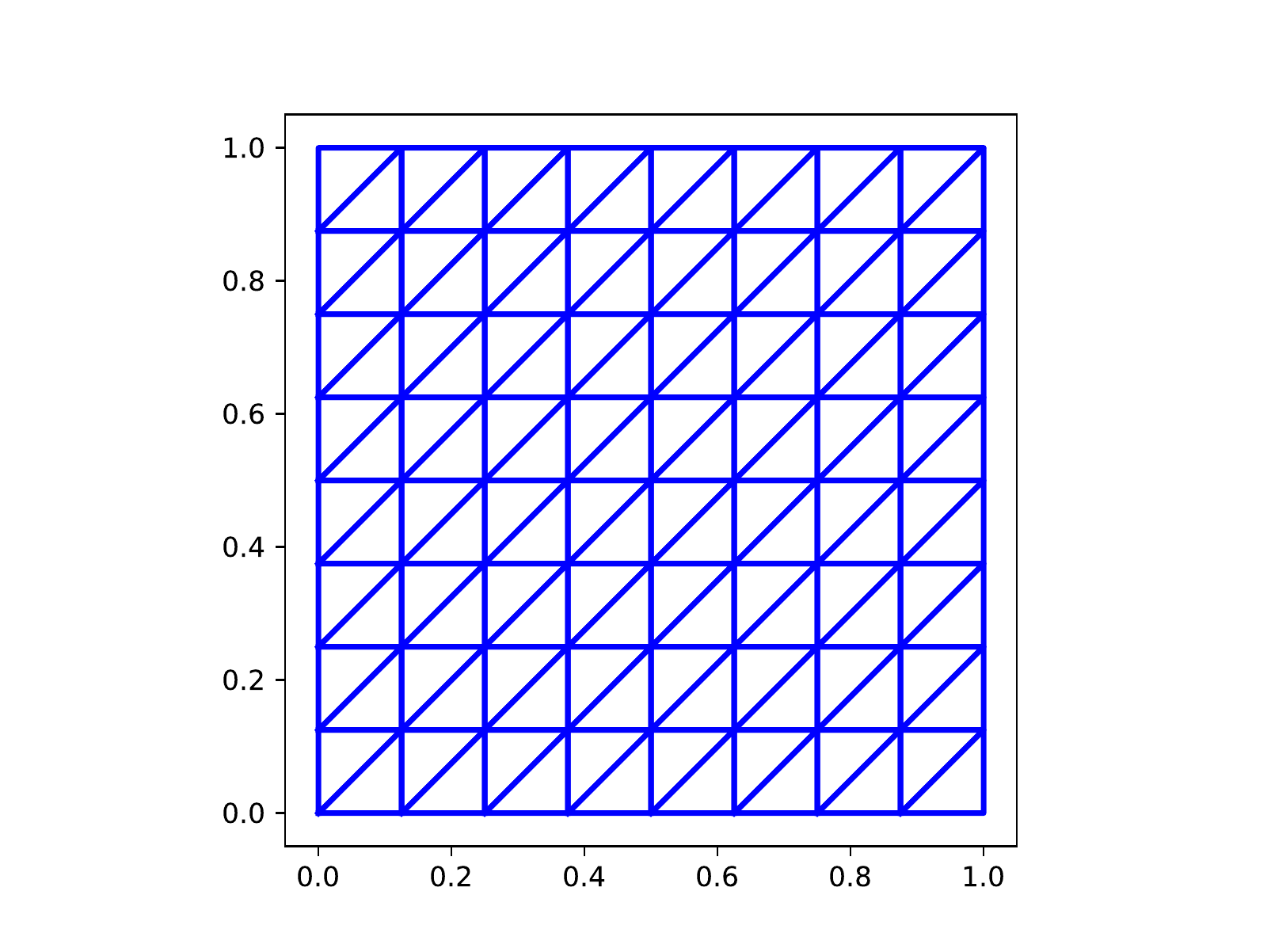}}
		\subfloat[Right.2]{\label{fig:right2} \includegraphics[trim=50 16 40 50,width=0.25\linewidth]{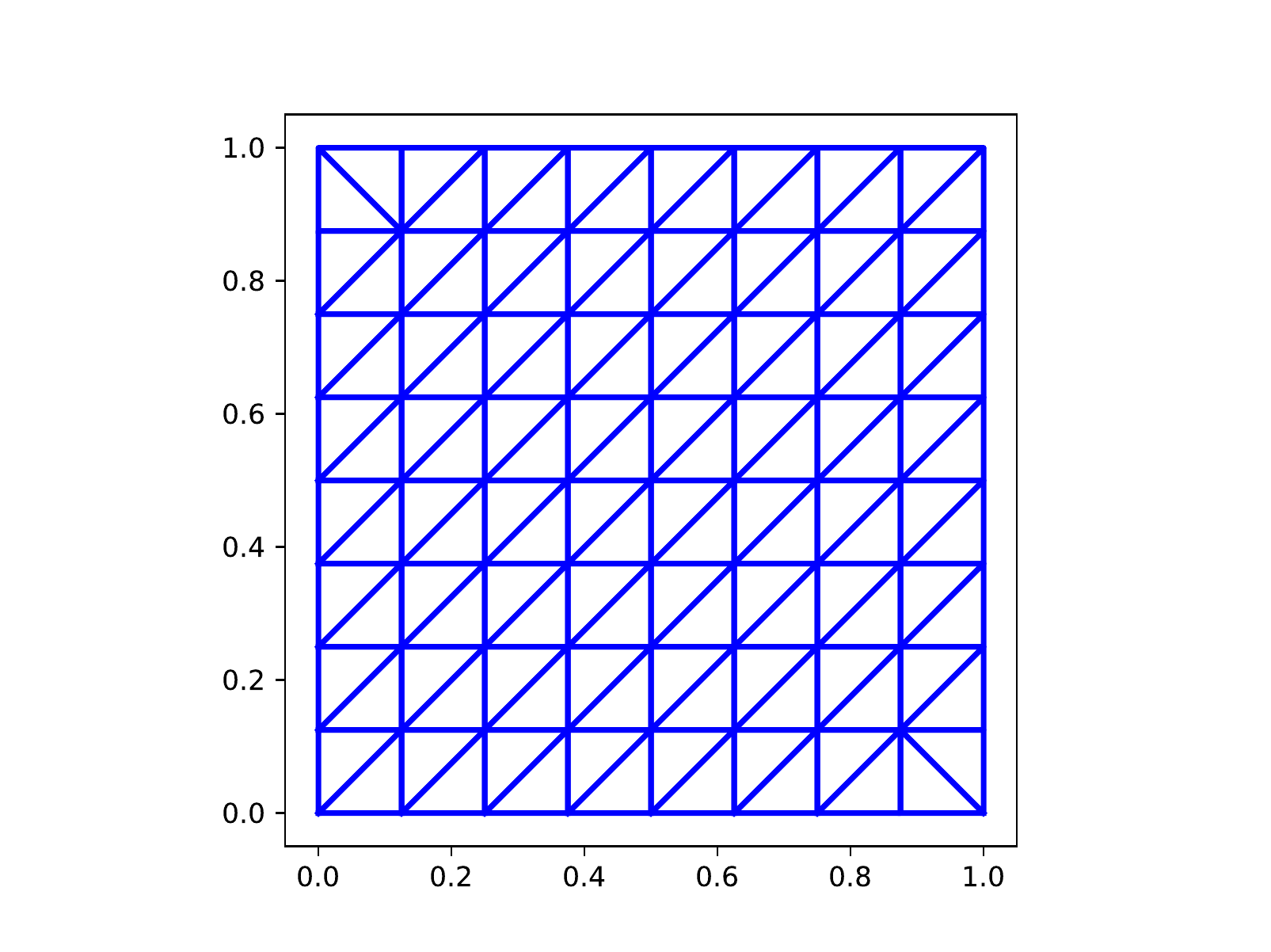}}
		\subfloat[Left]{\label{fig:left}\includegraphics[trim=50 16 40 50,width=0.25\linewidth]{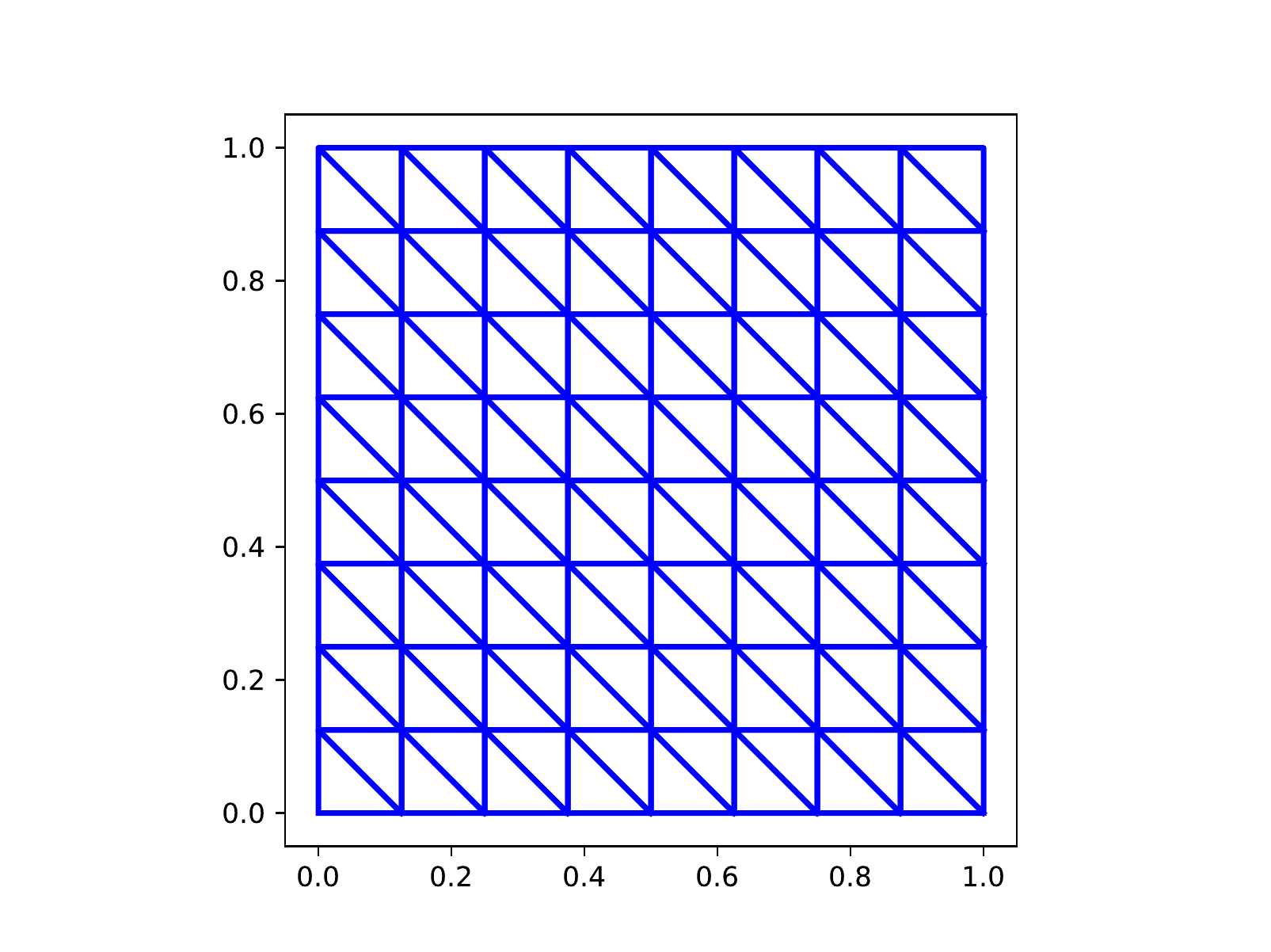}}
		\subfloat[Unstructured]{\label{fig:unstr}\includegraphics[trim=50 16 40 50,width=0.25\linewidth]{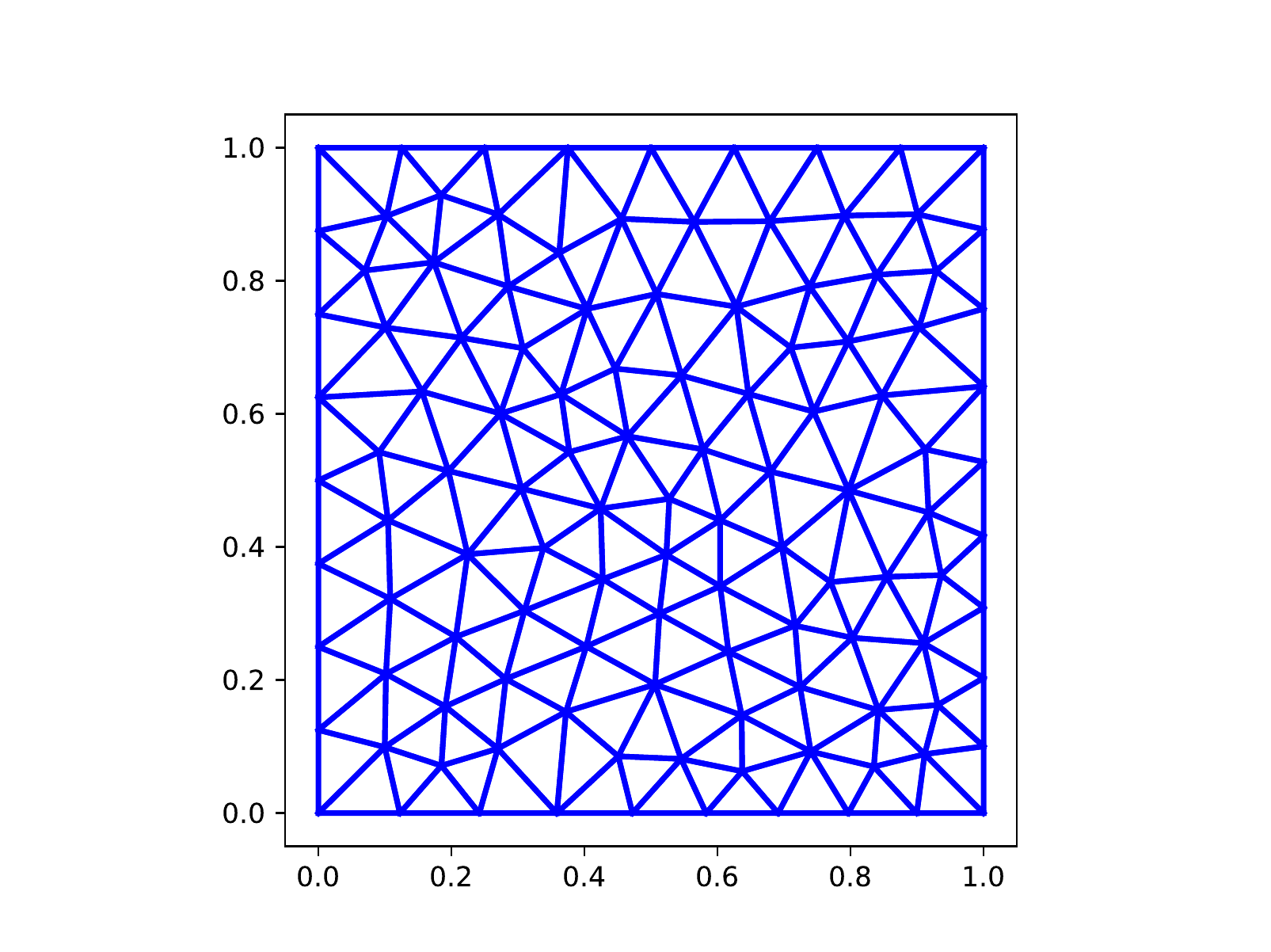}}
		\caption{Some examples of meshes. (a) and (b) are used for the fluid discretization and, in particular, in (b) we can see the correction of the corner triangles \blue we need when we work with the $\Pcal_1+\Pcal_0$ approximation of the pressure\noblue. In (c) and (d), we have two meshes used for the structure. When we work with perfectly matching meshes, we use the type in (a) for both the domains.}
		\label{fig:meshes}
	\end{figure}

	\lg Some examples are presented in Figure \ref{fig:meshes}, while in Table \ref{table:mesh_dofs}, we report the number of degrees of freedom of each unknown in dependence of the chosen discrete spaces and meshes. \gl

	\begin{table}[h!]
		\begin{center}
			\begin{tabular}{ccccc}
				\hline
				DOFs $\u_h$ & DOFs $p_h$ & DOFs $p^0_h$ & DOFs $\X_h$, $\llambda_h$ & DOFs $\X_h^{uns}$, $\llambda_h^{uns}$\\
				\hline
				2,178  		& 289  		& 801  		& 162  		& 232 \\
				8,450  		& 1,089 	& 3,137  	& 578  		& 742 \\
				33,282 		& 4,225		& 12,417  	& 2,178  	& 2,788 \\
				132,028 	& 16,641  	& 49,409  	& 8,450  	& 11,018 \\
				526,338  	&  66,049 	& 197,121  	& 33,282  	& 43,734 \\
				2,101,250 	& 263,169 	& 787,457  	& 132,028 	& 174,316\\
				\hline
			\end{tabular}
			\caption{Mesh parameters and degrees of freedom related to the chosen discrete spaces. The last column regards the solid discretization with unstructured mesh, while the others are related to uniform meshes. With $p_h^0$, we denote the pressure approximation with $\Pcal_1+\Pcal_0$ elements.}
			\label{table:mesh_dofs}
		\end{center}
	\end{table}
	
	\subsection{Mesh intersection in Python}
	All the tests has been performed with Python 3.8 scripts written specifically for this purpose. We focus our attention on the libraries used for the geometric computations.
	
	In Algorithms \ref{alg:intersection} and \ref{alg:no_intersection}, we denoted by $N_F$ and $N_S$ the number of elements of the two meshes (fluid and solid respectively), hence the number of the intersections to be tested is $N_S\times N_F$, which may become large when working with very fine meshes. In order to avoid useless computations, we can check in advance if two elements have non-empty intersection making use of the bounding box technique: if the bounding boxes are disjoint, the intersection of the two elements under consideration is empty and therefore we can move on and examine another pair of elements. These preliminary tests can be done using a tree-search algorithm: in our Python code, this procedure is based on the \verb!Rtree! library.
	
	To actually compute the intersection between two triangles, the \verb!Shapely.geometry! library is very useful because it can manage objects of type \verb!polygon! that allow us to compute and represent the intersection in an easy way: it is enough to construct two \verb!polygon! objects to represent the two triangles under consideration and then use the \verb!intersection()! method to obtain the resulting polygon. 
	
	On the other hand, if we assemble $C_f$ without intersection, we do not need particular tools: if we consider a solid triangle with its quadrature nodes, we can simply convert them with respect to the barycentric coordinates system related to each fluid element and understand in which one they are placed.
	
\begin{figure}
		\begin{center}
			\subfloat[Test 1]{\includegraphics[width=5cm]{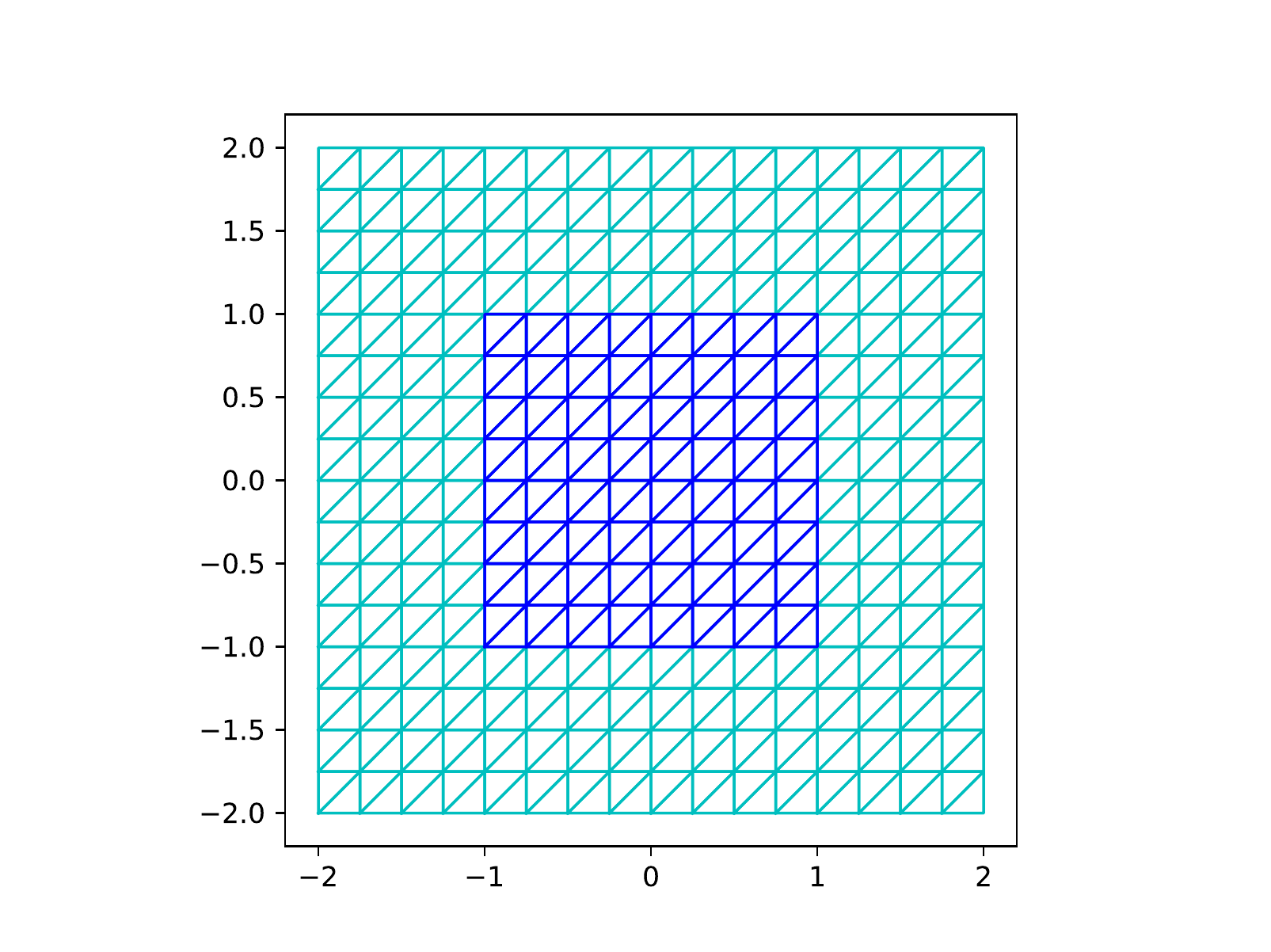}}
			\subfloat[Test 2]{\includegraphics[width=5cm]{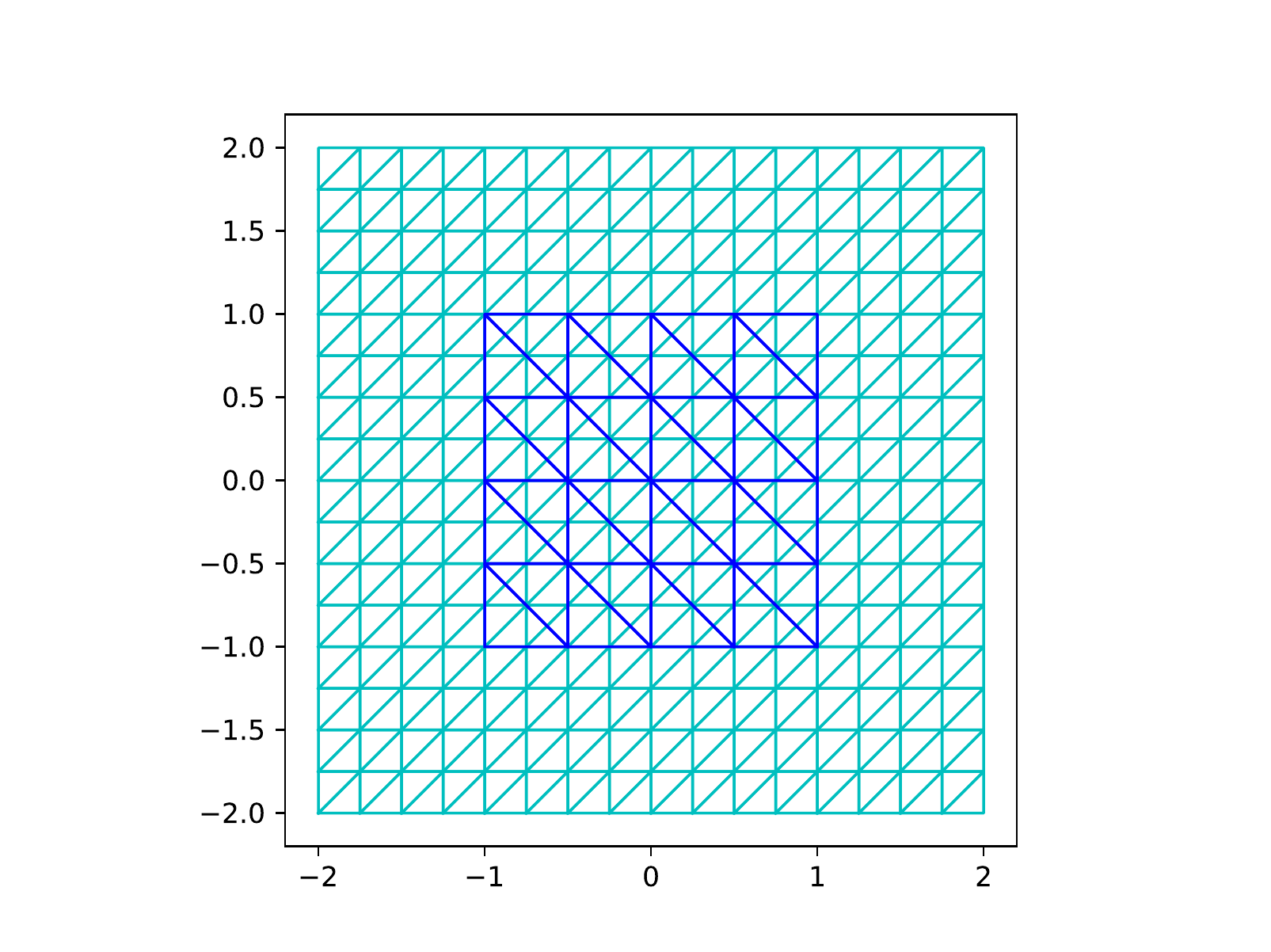}}
			\subfloat[Test 3]{\includegraphics[width=5cm]{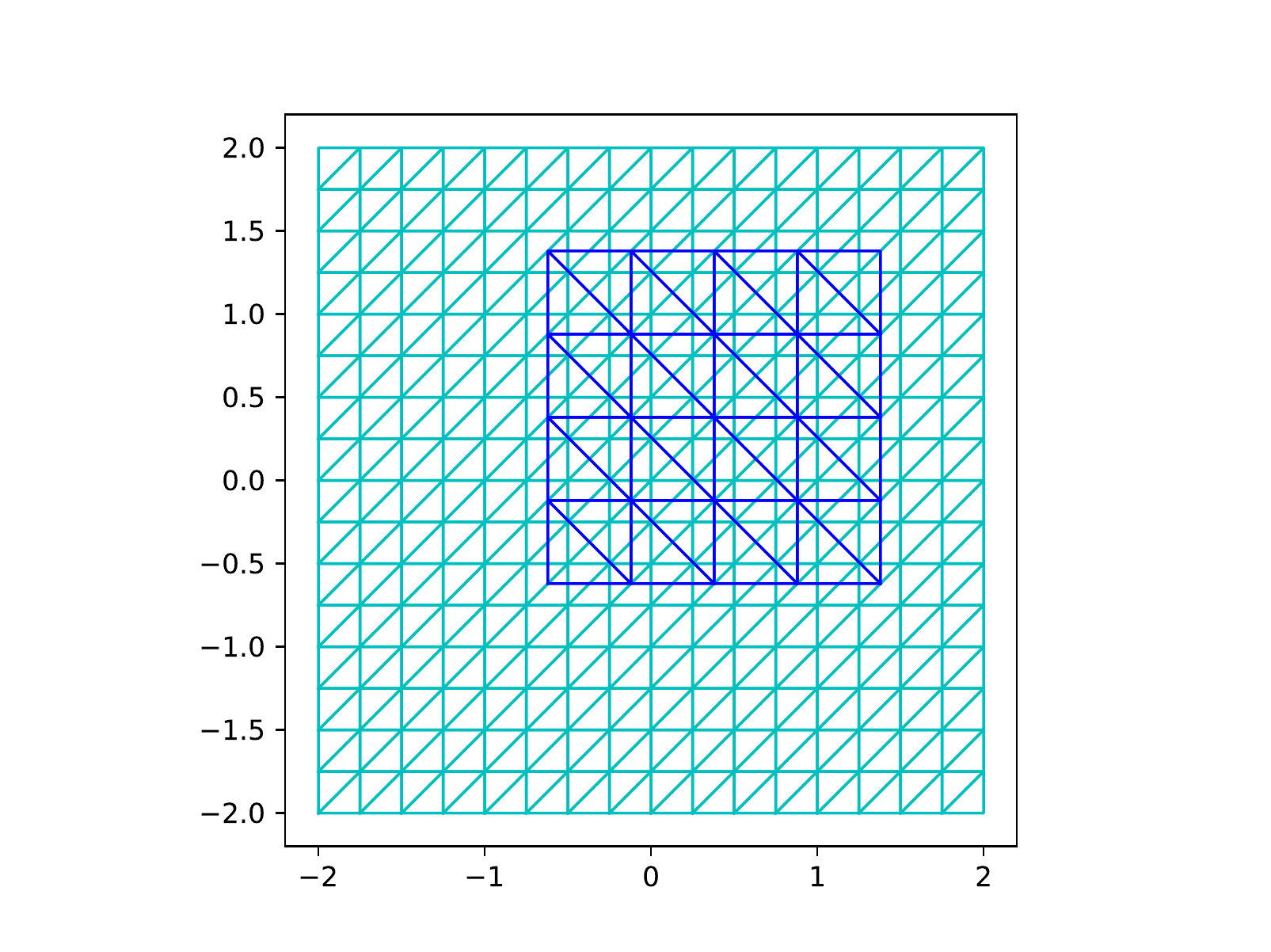}}\\
			\subfloat[Test 4]{\includegraphics[width=5cm]{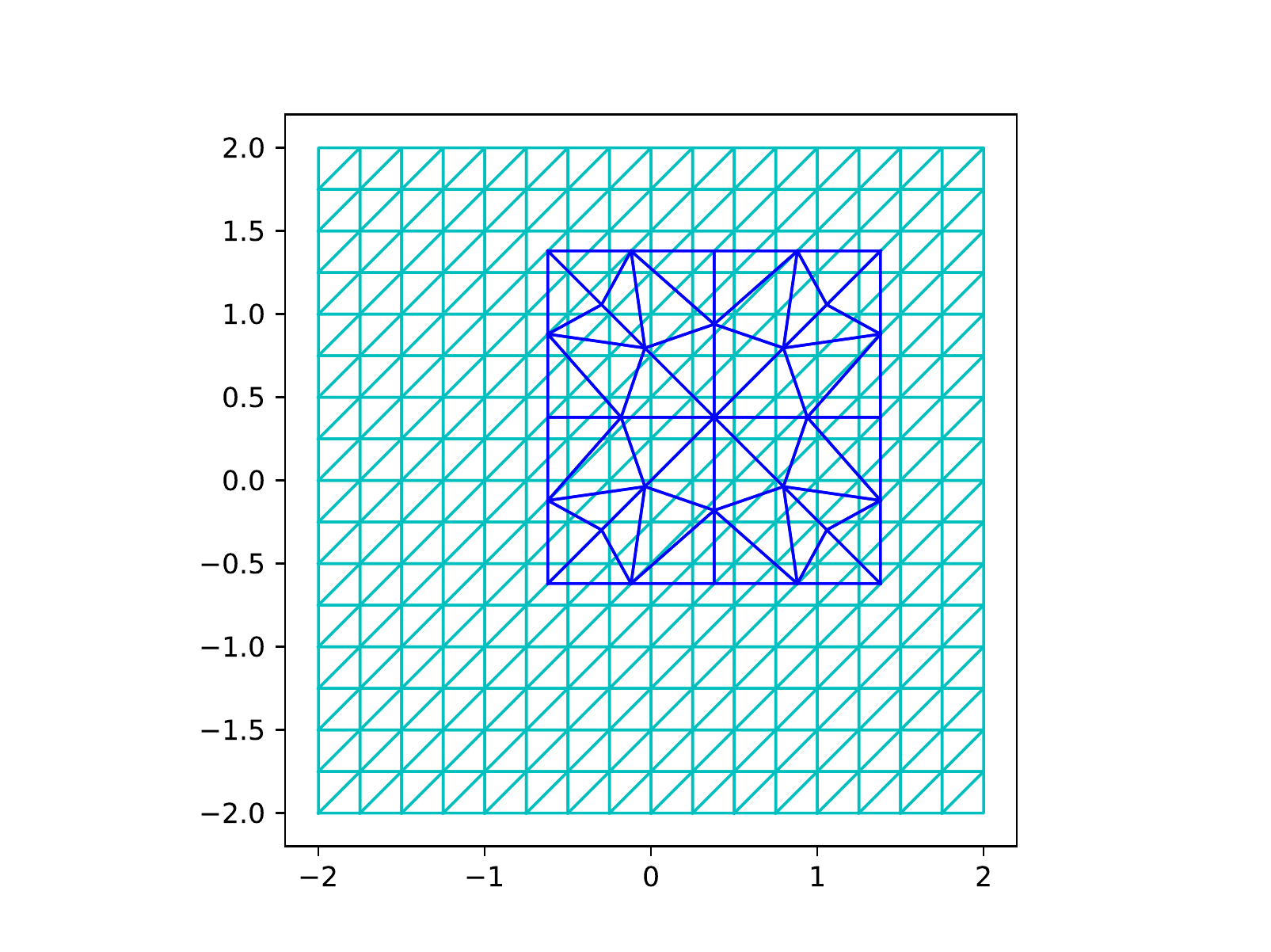}}
			\subfloat[Test 5]{\includegraphics[width=5cm]{figures/test_2_5-eps-converted-to}}
			\subfloat[Test 6]{\includegraphics[width=5cm]{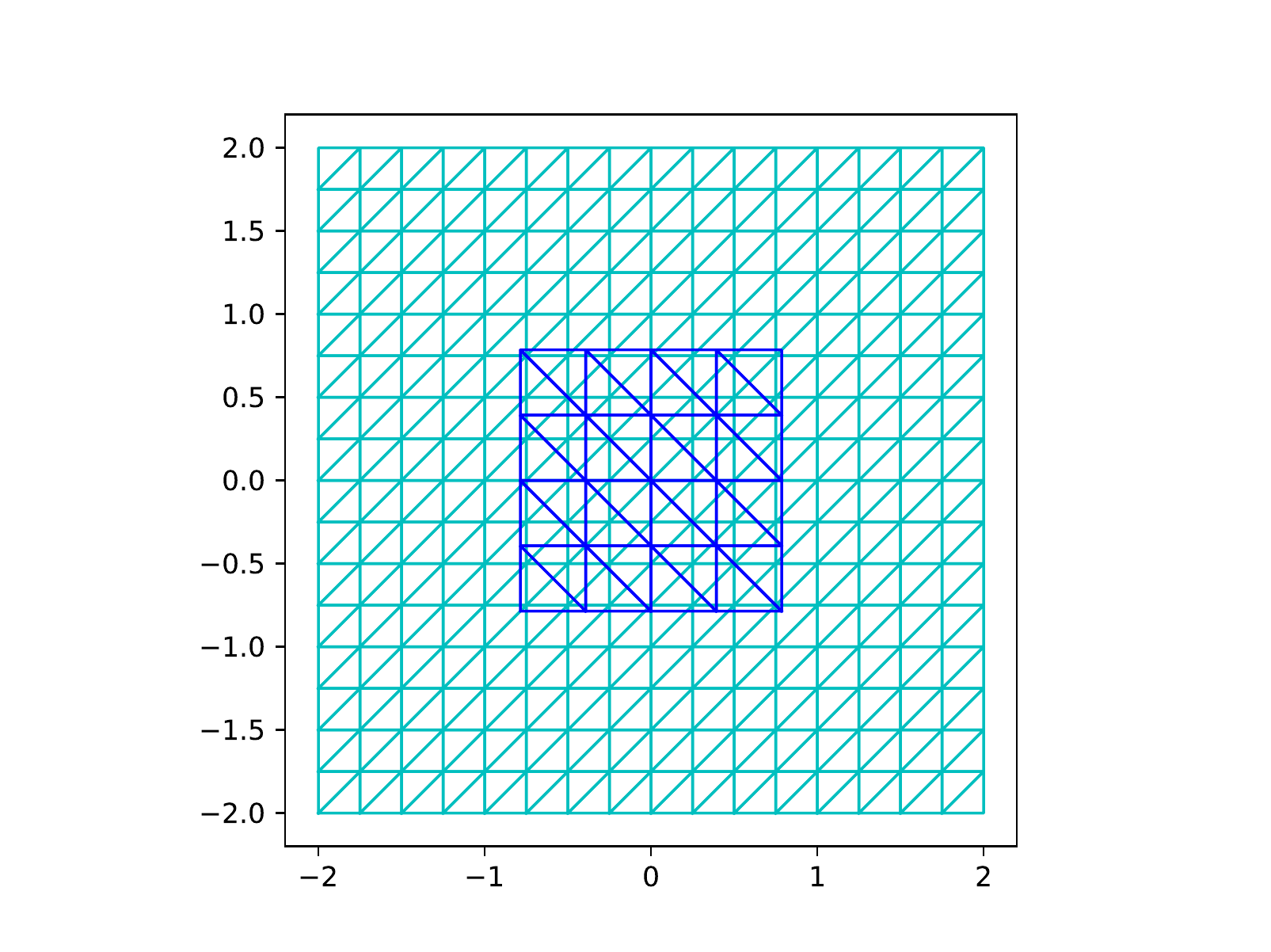}}\\
			\subfloat[Test 7]{\includegraphics[width=5cm]{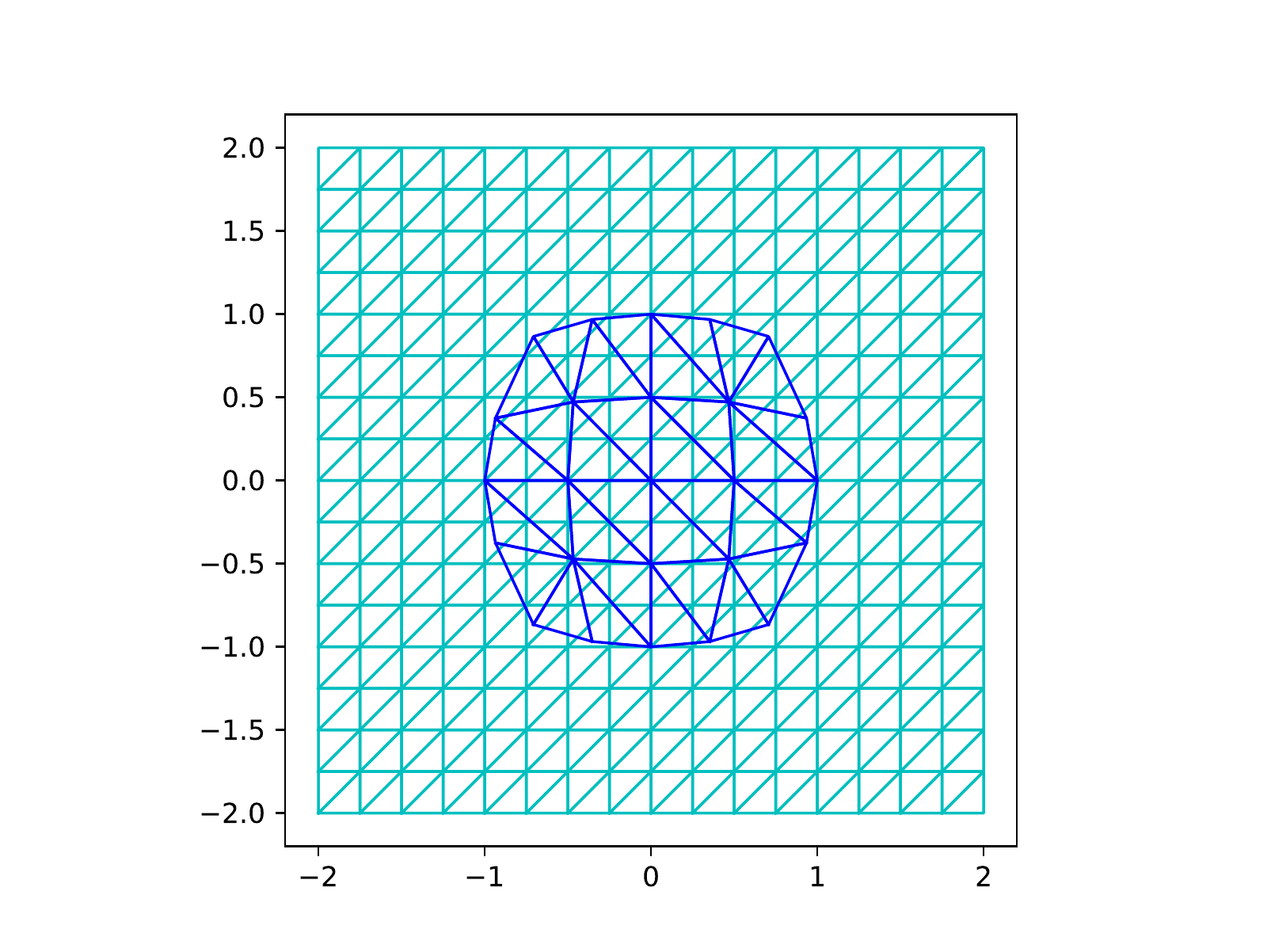}}
			\subfloat[Test 8]{\includegraphics[width=5cm]{figures/test3_geo-eps-converted-to}}\\
			\caption{All the geometric configurations of our tests with coarse meshes. In particular, for the first six tests, we have that the reference and actual configuration of the solid coincide. \blue We point out that the fluid mesh represented is the one related to the velocity. Moreover, each geometric configuration remains unchanged when we refine the two meshes\noblue.}
			\label{fig:geo_configs}
		\end{center}		
	\end{figure}

	The geometric configurations of our tests are summarized in Figure \ref{fig:geo_configs}: we can see the reciprocal position of fluid and solid meshes in the case of coarse spatial discretizations.
	
	\subsection{Quadrature rules for $\c_h$}
	Now we focus our attention on the quadrature rules used for the assembling of the interface matrix: as mentioned above, we set $\c$ as in \eqref{eq:c_h1}. With the chosen finite element spaces, the basis functions both for the velocity and the Lagrange multiplier are piecewise linear and this implies that we need to compute integrals with second order precision for the $L^2$ contributions and with first order precision for the gradients terms. In particular, in order to have exact computations, we used Gauss quadrature rules with the required accuracy for both approaches. Furthermore, in order to better understand the behavior of the system when the interface matrix is assembled without mesh intersection, we implemented, only for this case, also a Gauss quadrature rule with third order precision.
	
	For the sake of precision, we report the definition of these formulas on the reference triangle $\hatT$.
	
	\begin{defi}[Gauss rule, order 1]
		A polynomial of degree one $f$ can be integrated exactly by the following one-point formula
		\begin{equation}
			\int_{\hatT} f(\x)\,d\x \approx |\hatT| \, f(\x_{\hatT})
		\end{equation}
		where $\x_{\hatT}$ is the barycenter of $\hatT$.
	\end{defi}
	
	\begin{defi}[Gauss rule, order 2]
		A polynomial of degree two $f$ can be integrated exactly by the following three-points formula
		\begin{equation}
			\int_{\hatT} f(\x)\,d\x \approx \frac{|\hatT|}{3}\sum_{k=1}^3f\big(\x^{(k)}\big)
		\end{equation}
		where the quadrature nodes, represented in barycentric coordinates, are $\x^{(1)}=(2/3,1/6,1/6)$ and its permutation.
		
	\end{defi}
	
	\begin{defi}[Gauss rule, order 3]
		A polynomial of degree three $f$ can be integrated exactly by the following four-points formula
		\begin{equation}
			\int_{\hatT} f(\x)\,d\x \approx |\hatT| \bigg( \frac{25}{48}\sum_{k=1}^3f\big(\x^{(k)}\big)  - \frac{9}{16}f(\x_{\hatT})\bigg)
		\end{equation}
		where the quadrature nodes, represented in barycentric coordinates, are $\x^{(1)}=(3/5,1/5,1/5)$ with its permutation and the barycenter $\x_{\hatT}$ of $\hatT$.
	\end{defi}
	
	\subsection{Test 1. Sanity check}\label{sub:test1}
	Let us consider Problem~\ref{pro:stationary_tests} on the fluid square domain $\Omega=[-2,2]^2$ and with $\B=\Os=[-1,1]^2$. This assumption implies that the initial datum $\Xbar$ is simply the identity, hence $\d=\mathbf{0}$. We compute the right hand sides $\f$, $\g$, and $\d$ in order to obtain the following solutions:
	\begin{equation}\label{eq:standard_sol}
		\begin{aligned}
			&\u(x,y)=\curl\big((4-x^2)^2(4-y^2)^2\big)\\
			&\u_{\mid\partial\Omega}=\mathbf{0}\\
			&p(x,y)=150\sin(x)\\
			&\X(x,y)=\u(x,y)\\
			&\llambda(x,y)=\big(e^x,e^y\big).
		\end{aligned}
	\end{equation}
	In particular, we discretize both domains using perfectly matching uniform right-oriented meshes \blue(Figures \ref{fig:right1} and \ref{fig:right2})\noblue: we can see that the two methods for the assembling of the interface matrix are equivalent, indeed, the convergence rates confirm what expected. The errors and the convergence rates are collected in Tables ~\ref{table:test1_p1_fluid}, \ref{table:test1_p1_solid}, \ref{table:test1_p1p0_fluid}, \ref{table:test1_p1p0_solid}: in particular in the first two tables, we have the data related to the approximation with $\Pcal_1-iso-\Pcal_2/\Pcal_1/\Pcal_1/\Pcal_1$ element, while in Table~\ref{table:test1_p1p0_fluid} and \ref{table:test1_p1p0_solid}, we have the results related to the $\Pcal_1-iso-\Pcal_2/\Pcal_1+\Pcal_0/\Pcal_1/\Pcal_1$ approximation.
	
	\lg In Subsection \ref{sec:BP}, we discussed the main features of the Bercovier-Pironneau element, which is a low order Stokes pair. Therefore, we expect that the errors $\|p-p_h\|_{0}$, $\|\u-\u_h\|_{1}$ decay with order one, while $\|\u-\u_h\|_{0}$ decays with order two. Moreover, since for both displacement and multiplier we use $\Pcal_1$ elements, we expect that they converge with order two in the $L^2$ norm and with order one in the $H^1$ norm. From the results, we can see that, in the case of the classical Bercovier-Pironneau element, we get a superconvergence for pressure and Lagrange multiplier.\gl
	
	\begin{table}[h!]
		\begin{center}
			\begin{tabular}{p{1.3cm}|p{2.3cm}p{1.3cm}|p{2.3cm}p{1.3cm}|p{2.3cm}p{1.3cm}}
				\hline
				\multicolumn{7}{c}{\textbf{Errors and convergence rates for Test 1 $\bullet$ $\mathbf{\Pcal_1-iso-\Pcal_2/\Pcal_1/\Pcal_1/\Pcal_1}$}}\\
				\hline
				$h_\mathcal{T}$ & \multicolumn{2}{c|}{$\|p-p_h\|_{0}$} & \multicolumn{2}{c|}{$\|\u-\u_h\|_{0}$} & \multicolumn{2}{c}{$\|\u-\u_h\|_{1}$} \\
				& Error & Rate & Error & Rate & Error & Rate \\
				\hline
				\multicolumn{7}{c}{\textit{Coupling with mesh intersection}}\\
				\hline
				1/4&2.102e-02&-   &9.622e-03&-   &7.684e-02&-   \\
				1/8&6.251e-03&1.75&2.408e-03&2.00&3.834e-02&1.00\\
				1/16&1.977e-03&1.66&6.017e-04&2.00&1.915e-02&1.00\\
				1/32&6.546e-04&1.59&1.504e-04&2.00&9.572e-03&1.00\\
				1/64&2.233e-04&1.55&3.758e-05&2.00&4.785e-03&1.00\\
				1/128&7.745e-05&1.53&9.394e-06&2.00&2.392e-03&1.00\\
				\hline
				\multicolumn{7}{c}{\textit{Coupling without mesh intersection, quad. rule of order 2}}\\
				\hline
				1/4&2.102e-02&-   &9.622e-03&-   &7.684e-02&-   \\
				1/8&6.251e-03&1.75&2.408e-03&2.00&3.834e-02&1.00\\
				1/16&1.977e-03&1.66&6.017e-04&2.00&1.915e-02&1.00\\
				1/32&6.546e-04&1.59&1.504e-04&2.00&9.572e-03&1.00\\
				1/64&2.233e-04&1.55&3.758e-05&2.00&4.785e-03&1.00\\
				1/128&7.745e-05&1.53&9.394e-06&2.00&2.392e-03&1.00\\
				\hline
				\multicolumn{7}{c}{\textit{Coupling without mesh intersection, quad. rule of order 3}}\\
				\hline
				1/4&2.102e-02&-   &9.622e-03&-   &7.684e-02&-   \\
				1/8&6.251e-03&1.75&2.408e-03&2.00&3.834e-02&1.00\\
				1/16&1.977e-03&1.66&6.017e-04&2.00&1.915e-02&1.00\\
				1/32&6.546e-04&1.59&1.504e-04&2.00&9.572e-03&1.00\\
				1/64&2.233e-04&1.55&3.758e-05&2.00&4.785e-03&1.00\\
				1/128&7.745e-05&1.53&9.394e-06&2.00&2.392e-03&1.00\\
				\hline
			\end{tabular}
			\caption{Errors and convergence rates for the fluid variables of Test 1 discretized with $\Pcal_1-iso-\Pcal_2/\Pcal_1/\Pcal_1/\Pcal_1$}
			\label{table:test1_p1_fluid}
		\end{center}
	\end{table}

	\begin{table}[h!]
		\begin{center}
			\begin{tabular}{p{1cm}|p{1.8cm}p{0.8cm}|p{1.8cm}p{.8cm}|p{1.8cm}p{.8cm}|p{1.8cm}p{.8cm}}
				\hline
				\multicolumn{9}{c}{\textbf{\textbf{Errors and convergence rates for Test 1 $\bullet$ $\mathbf{\Pcal_1-iso-\Pcal_2/\Pcal_1/\Pcal_1/\Pcal_1}$}}}\\
				\hline
				$h_\mathcal{S}$ & \multicolumn{2}{c|}{$\|\X-\X_h\|_{0,\B}$} & \multicolumn{2}{c|}{$\|\X-\X_h\|_{1,\B}$} & \multicolumn{2}{c|}{$\|\llambda-\llambda_h\|_{0,\B}$}& \multicolumn{2}{c}{$\|\llambda-\llambda_h\|_{1,\B}$} \\
				& Error & Rate & Error & Rate & Error & Rate & Error & Rate\\
				\hline
				\multicolumn{9}{c}{\textit{Coupling with mesh intersection}}\\
				\hline
				1/8&8.011e-03&-   &4.972e-02&-   &1.908e-01&-   &4.166e-01&-   \\
				1/16&2.011e-03&1.99&2.479e-02&1.00&4.786e-02&2.00&1.111e-01&1.91\\
				1/32&5.032e-04&2.00&1.239e-02&1.00&1.197e-02&2.00&2.963e-02&1.91\\
				1/64&1.258e-04&2.00&6.193e-03&1.00&2.991e-03&2.00&8.185e-03&1.86\\
				1/128&3.146e-05&2.00&3.096e-03&1.00&7.476e-04&2.00&2.524e-03&1.70\\
				1/256&7.864e-06&2.00&1.548e-03&1.00&1.869e-04&2.00&9.445e-04&1.42\\
				\hline
				\multicolumn{9}{c}{\textit{Coupling without mesh intersection, quad. rule of order 2}}\\
				\hline
				1/8&8.011e-03&-   &4.972e-02&-   &1.908e-01&-   &4.166e-01&-   \\
				1/16&2.011e-03&1.99&2.479e-02&1.00&4.786e-02&2.00&1.111e-01&1.91\\
				1/32&5.032e-04&2.00&1.239e-02&1.00&1.197e-02&2.00&2.963e-02&1.91\\
				1/64&1.258e-04&2.00&6.193e-03&1.00&2.991e-03&2.00&8.185e-03&1.86\\
				1/128&3.146e-05&2.00&3.096e-03&1.00&7.476e-04&2.00&2.524e-03&1.70\\
				1/256&7.864e-06&2.00&1.548e-03&1.00&1.869e-04&2.00&9.445e-04&1.42\\
				\hline
				\multicolumn{9}{c}{\textit{Coupling without mesh intersection, quad. rule of order 3}}\\
				\hline
				1/8&8.011e-03&-   &4.972e-02&-   &1.908e-01&-   &4.166e-01&-   \\
				1/16&2.011e-03&1.99&2.479e-02&1.00&4.786e-02&2.00&1.111e-01&1.91\\
				1/32&5.032e-04&2.00&1.239e-02&1.00&1.197e-02&2.00&2.963e-02&1.91\\
				1/64&1.258e-04&2.00&6.193e-03&1.00&2.991e-03&2.00&8.185e-03&1.86\\
				1/128&3.146e-05&2.00&3.096e-03&1.00&7.476e-04&2.00&2.524e-03&1.70\\
				1/256&7.864e-06&2.00&1.548e-03&1.00&1.869e-04&2.00&9.445e-04&1.42\\
				\hline
			\end{tabular}
			\caption{Errors and convergence rates for the solid variables of Test 1 discretized with $\Pcal_1-iso-\Pcal_2/\Pcal_1/\Pcal_1/\Pcal_1$}
			\label{table:test1_p1_solid}
		\end{center}
	\end{table}

	\begin{table}[h!]
		\begin{center}
			\begin{tabular}{p{1.3cm}|p{2.3cm}p{1.3cm}|p{2.3cm}p{1.3cm}|p{2.3cm}p{1.3cm}}
				\hline
				\multicolumn{7}{c}{\textbf{Errors and convergence rates for Test 1 $\bullet$ $\mathbf{\Pcal_1-iso-\Pcal_2/\Pcal_1+\Pcal_0/\Pcal_1/\Pcal_1}$}}\\
				\hline
				$h_\mathcal{T}$ & \multicolumn{2}{c|}{$\|p-p_h\|_{0}$} & \multicolumn{2}{c|}{$\|\u-\u_h\|_{0}$} & \multicolumn{2}{c}{$\|\u-\u_h\|_{1}$} \\
				& Error & Rate & Error & Rate & Error & Rate \\
				\hline
				\multicolumn{7}{c}{\textit{Coupling with mesh intersection}}\\
				\hline
				1/4&7.981e-02&-   &1.043e-02&-   &8.042e-02&-   \\
				1/8&3.939e-02&1.02&2.617e-03&1.99&4.017e-02&1.00\\
				1/16&1.957e-02&1.01&6.549e-04&2.00&2.008e-02&1.00\\
				1/32&9.749e-03&1.00&1.637e-04&2.00&1.004e-02&1.00\\
				1/64&4.866e-03&1.00&4.093e-05&2.00&5.018e-03&1.00\\
				1/128&2.431e-03&1.00&1.023e-05&2.00&2.509e-03&1.00\\
				\hline
				\multicolumn{7}{c}{\textit{Coupling without mesh intersection, quad. rule of order 2}}\\
				\hline
				1/4&7.981e-02&-   &1.043e-02&-   &8.042e-02&-   \\
				1/8&3.939e-02&1.02&2.617e-03&1.99&4.017e-02&1.00\\
				1/16&1.957e-02&1.01&6.549e-04&2.00&2.008e-02&1.00\\
				1/32&9.749e-03&1.00&1.637e-04&2.00&1.004e-02&1.00\\
				1/64&4.866e-03&1.00&4.093e-05&2.00&5.018e-03&1.00\\
				1/128&2.431e-03&1.00&1.023e-05&2.00&2.509e-03&1.00\\
				\hline
				\multicolumn{7}{c}{\textit{Coupling without mesh intersection, quad. rule of order 3}}\\
				\hline
				1/4&7.981e-02&-   &1.043e-02&-   &8.042e-02&-   \\
				1/8&3.939e-02&1.02&2.617e-03&1.99&4.017e-02&1.00\\
				1/16&1.957e-02&1.01&6.549e-04&2.00&2.008e-02&1.00\\
				1/32&9.749e-03&1.00&1.637e-04&2.00&1.004e-02&1.00\\
				1/64&4.866e-03&1.00&4.093e-05&2.00&5.018e-03&1.00\\
				1/128&2.431e-03&1.00&1.023e-05&2.00&2.509e-03&1.00\\
				\hline
			\end{tabular}
			\caption{Errors and convergence rates for the fluid variables of Test 1 discretized with $\Pcal_1-iso-\Pcal_2/\Pcal_1+\Pcal_0/\Pcal_1/\Pcal_1$}
			\label{table:test1_p1p0_fluid}
		\end{center}
	\end{table}
	
	\begin{table}[h!]
		\begin{center}
			\begin{tabular}{p{1cm}|p{1.8cm}p{0.8cm}|p{1.8cm}p{.8cm}|p{1.8cm}p{.8cm}|p{1.8cm}p{.8cm}}
				\hline
				\multicolumn{9}{c}{\textbf{\textbf{Errors and convergence rates for Test 1 $\bullet$ $\mathbf{\Pcal_1-iso-\Pcal_2/\Pcal_1+\Pcal_0/\Pcal_1/\Pcal_1}$}}}\\
				\hline
				$h_\mathcal{S}$ & \multicolumn{2}{c|}{$\|\X-\X_h\|_{0,\B}$} & \multicolumn{2}{c|}{$\|\X-\X_h\|_{1,\B}$} & \multicolumn{2}{c|}{$\|\llambda-\llambda_h\|_{0,\B}$}& \multicolumn{2}{c}{$\|\llambda-\llambda_h\|_{1,\B}$} \\
				& Error & Rate & Error & Rate & Error & Rate & Error & Rate\\
				\hline
				\multicolumn{9}{c}{\textit{Coupling with mesh intersection}}\\
				\hline
				1/8&8.854e-03&-   &5.239e-02&-   &2.300e-01&-   &1.861e+00&-   \\
				1/16&2.228e-03&1.99&2.616e-02&1.00&5.802e-02&1.99&9.338e-01&0.99\\
				1/32&5.578e-04&2.00&1.308e-02&1.00&1.453e-02&2.00&4.672e-01&1.00\\
				1/64&1.395e-04&2.00&6.539e-03&1.00&3.632e-03&2.00&2.336e-01&1.00\\
				1/128&3.487e-05&2.00&3.269e-03&1.00&9.078e-04&2.00&1.168e-01&1.00\\
				1/256&8.718e-06&2.00&1.635e-03&1.00&2.269e-04&2.00&5.840e-02&1.00\\
				\hline
				\multicolumn{9}{c}{\textit{Coupling without mesh intersection, quad. rule of order 2}}\\
				\hline
				1/8&8.854e-03&-   &5.239e-02&-   &2.300e-01&-   &1.861e+00&-   \\
				1/16&2.228e-03&1.99&2.616e-02&1.00&5.802e-02&1.99&9.338e-01&0.99\\
				1/32&5.578e-04&2.00&1.308e-02&1.00&1.453e-02&2.00&4.672e-01&1.00\\
				1/64&1.395e-04&2.00&6.539e-03&1.00&3.632e-03&2.00&2.336e-01&1.00\\
				1/128&3.487e-05&2.00&3.269e-03&1.00&9.078e-04&2.00&1.168e-01&1.00\\
				1/256&8.718e-06&2.00&1.635e-03&1.00&2.269e-04&2.00&5.840e-02&1.00\\
				\hline
				\multicolumn{9}{c}{\textit{Coupling without mesh intersection, quad. rule of order 3}}\\
				\hline
				1/8&8.854e-03&-   &5.239e-02&-   &2.300e-01&-   &1.861e+00&-   \\
				1/16&2.228e-03&1.99&2.616e-02&1.00&5.802e-02&1.99&9.338e-01&0.99\\
				1/32&5.578e-04&2.00&1.308e-02&1.00&1.453e-02&2.00&4.672e-01&1.00\\
				1/64&1.395e-04&2.00&6.539e-03&1.00&3.632e-03&2.00&2.336e-01&1.00\\
				1/128&3.487e-05&2.00&3.269e-03&1.00&9.078e-04&2.00&1.168e-01&1.00\\
				1/256&8.718e-06&2.00&1.635e-03&1.00&2.269e-04&2.00&5.840e-02&1.00\\
				\hline
			\end{tabular}
			\caption{Errors and convergence rates for the solid variables of Test 1 discretized with $\Pcal_1-iso-\Pcal_2/\Pcal_1+\Pcal_0/\Pcal_1/\Pcal_1$}
			\label{table:test1_p1p0_solid}
		\end{center}
	\end{table}
	
	\subsection{Test 2}\label{sub:test2}
	We consider the same problem as in Test 1, but changing the meshes: we discretize the fluid domain with a right-oriented uniform mesh \blue(Figures \ref{fig:right1} and \ref{fig:right2})\noblue, while the solid domain is discretized with a left-oriented uniform mesh \blue(Figure \ref{fig:left})\noblue.
	
\begin{figure}[!h]
		\centering
		\includegraphics[width=0.24\linewidth]{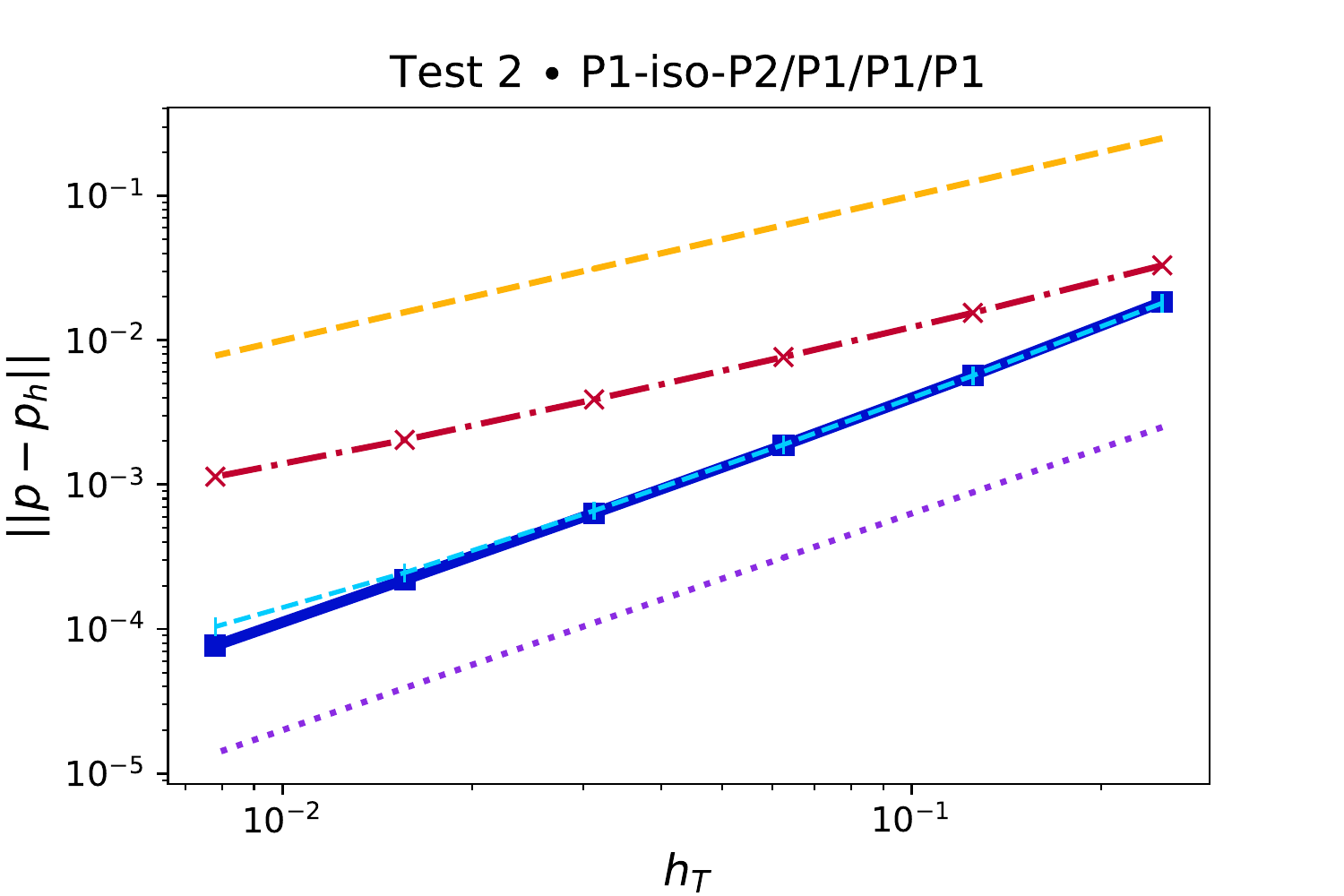}
		\includegraphics[width=0.24\linewidth]{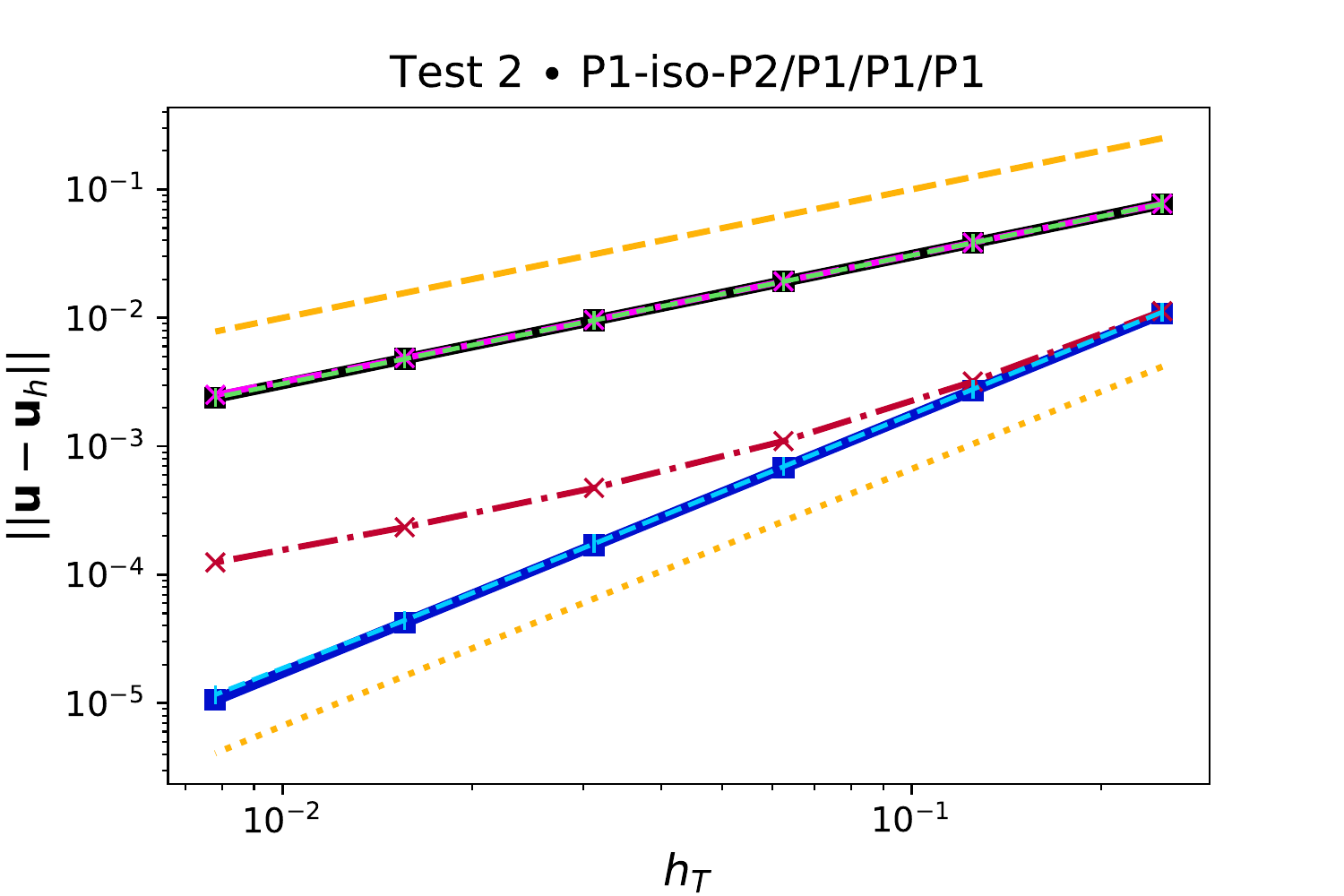}
		\includegraphics[width=0.24\linewidth]{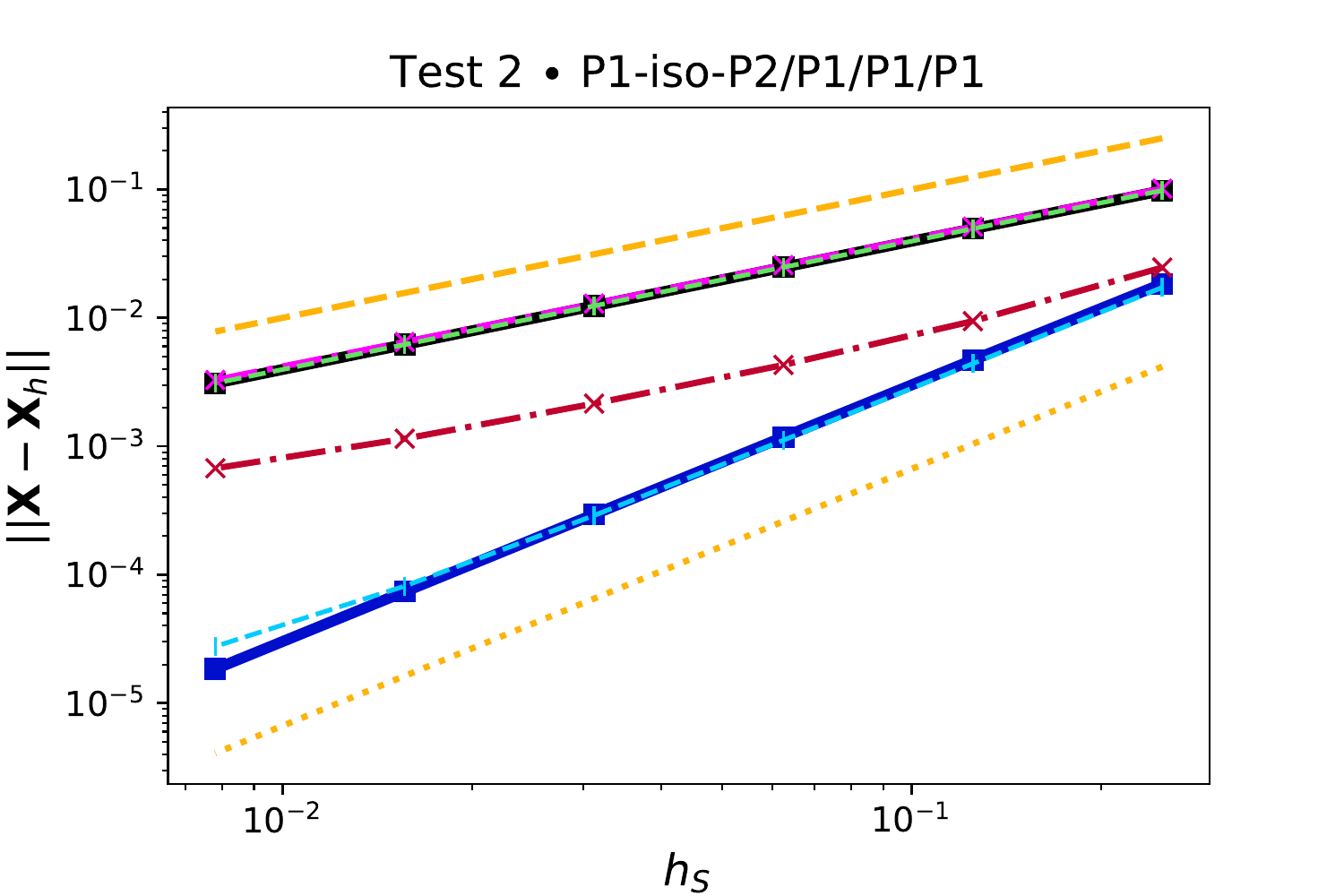}
		\includegraphics[width=0.24\linewidth]{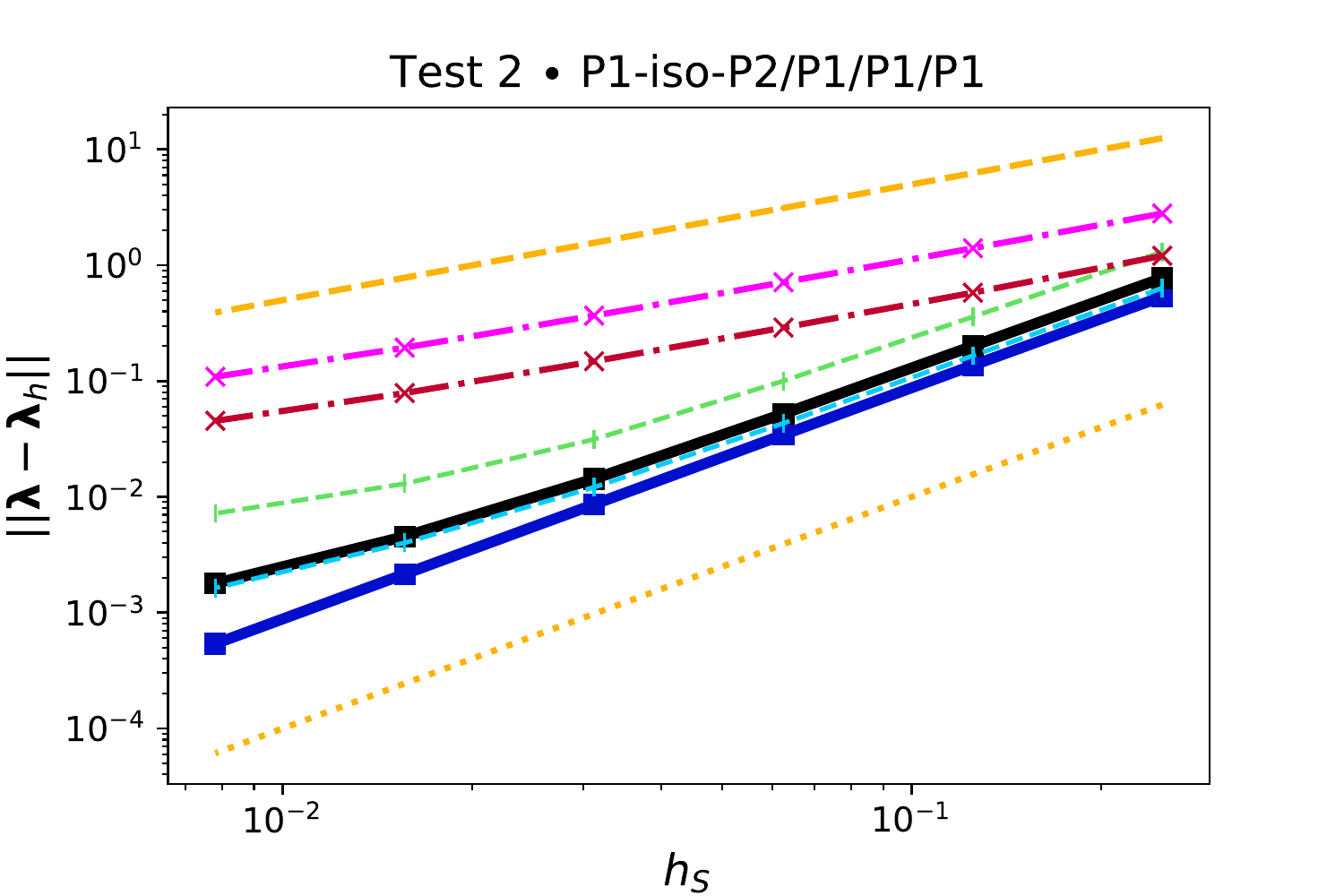}\\
		\includegraphics[trim=180 16 50 50,width=1.35\linewidth]{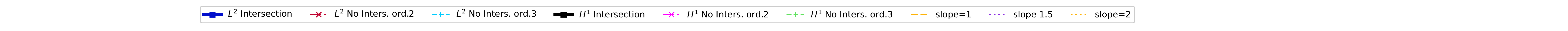}
		\caption{Convergence plots of Test 2 with $\Pcal_1-iso-\Pcal_2/\Pcal_1/\Pcal_1/\Pcal_1$}
		\label{fig:test2_p1}
	\end{figure}

\begin{figure}[!h]
		\centering
		\includegraphics[width=0.24\linewidth]{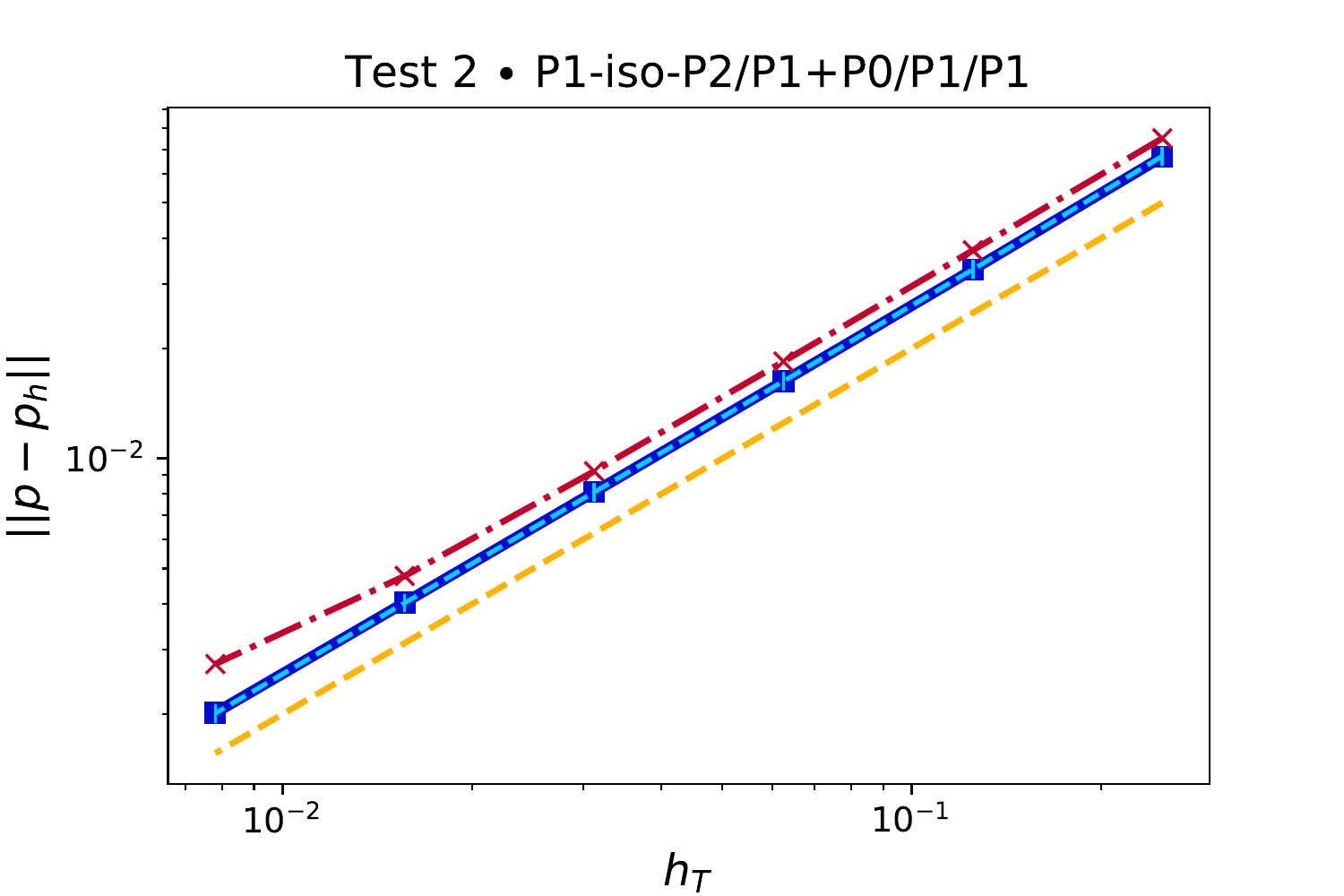}
		\includegraphics[width=0.24\linewidth]{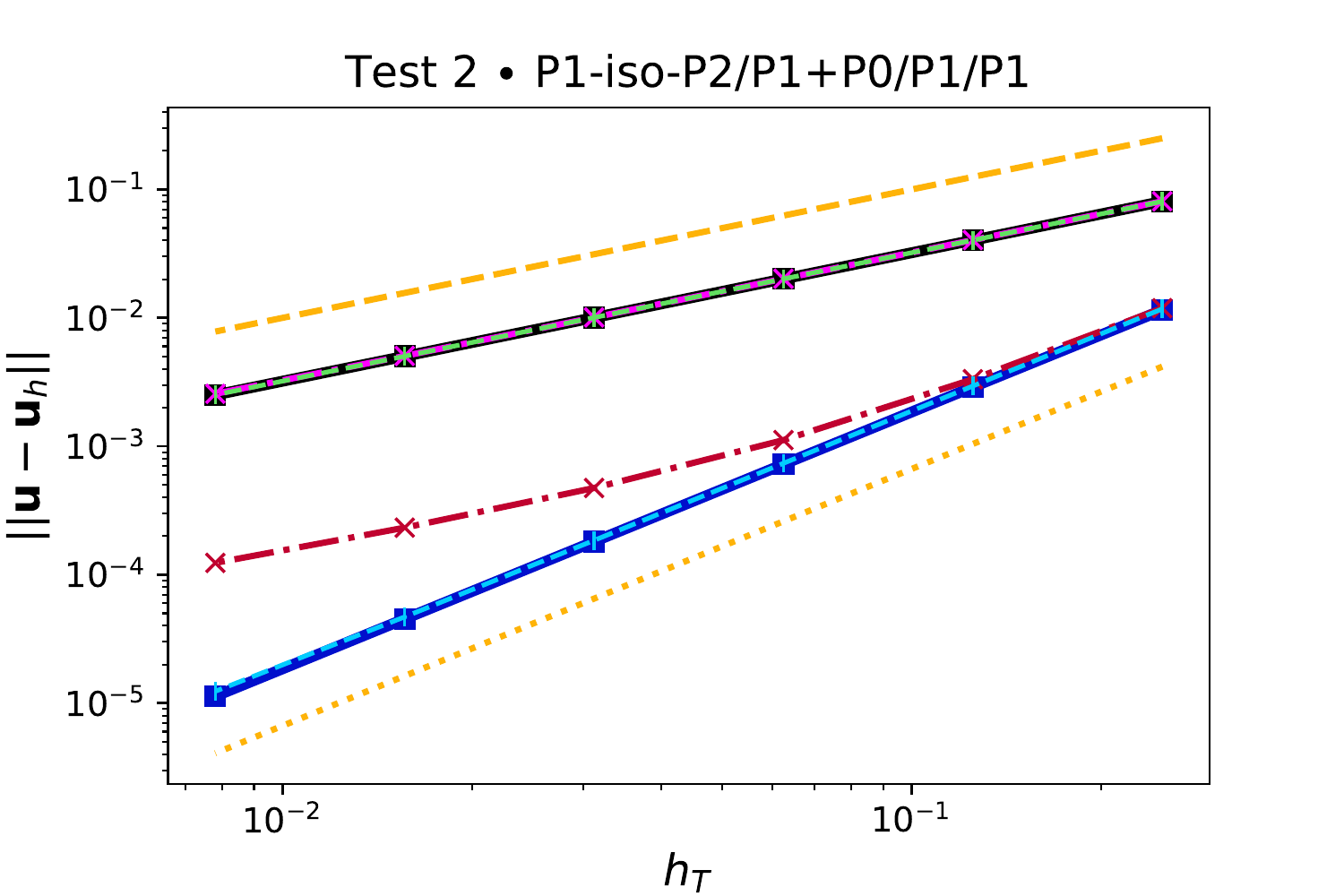}
		\includegraphics[width=0.24\linewidth]{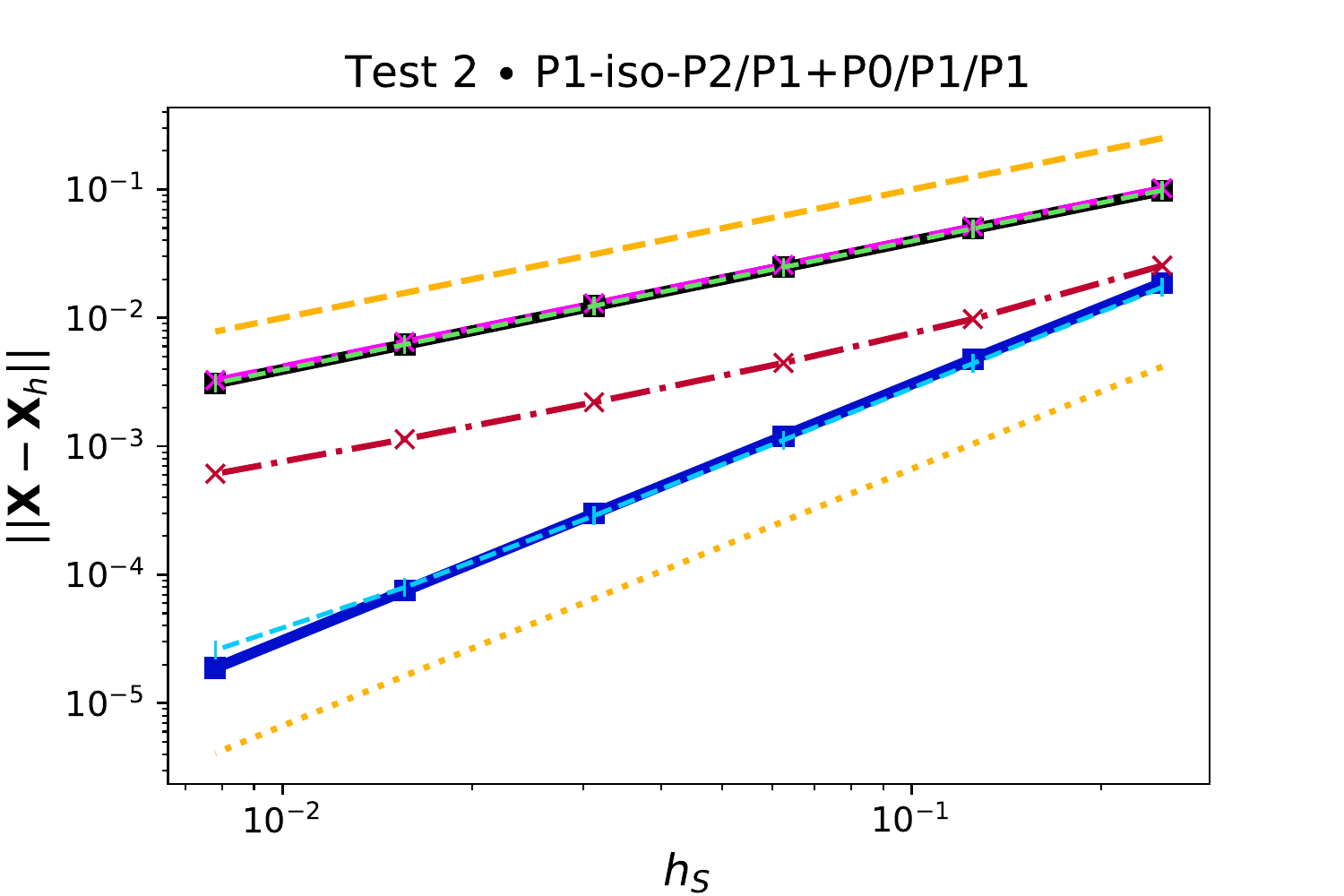}
		\includegraphics[width=0.24\linewidth]{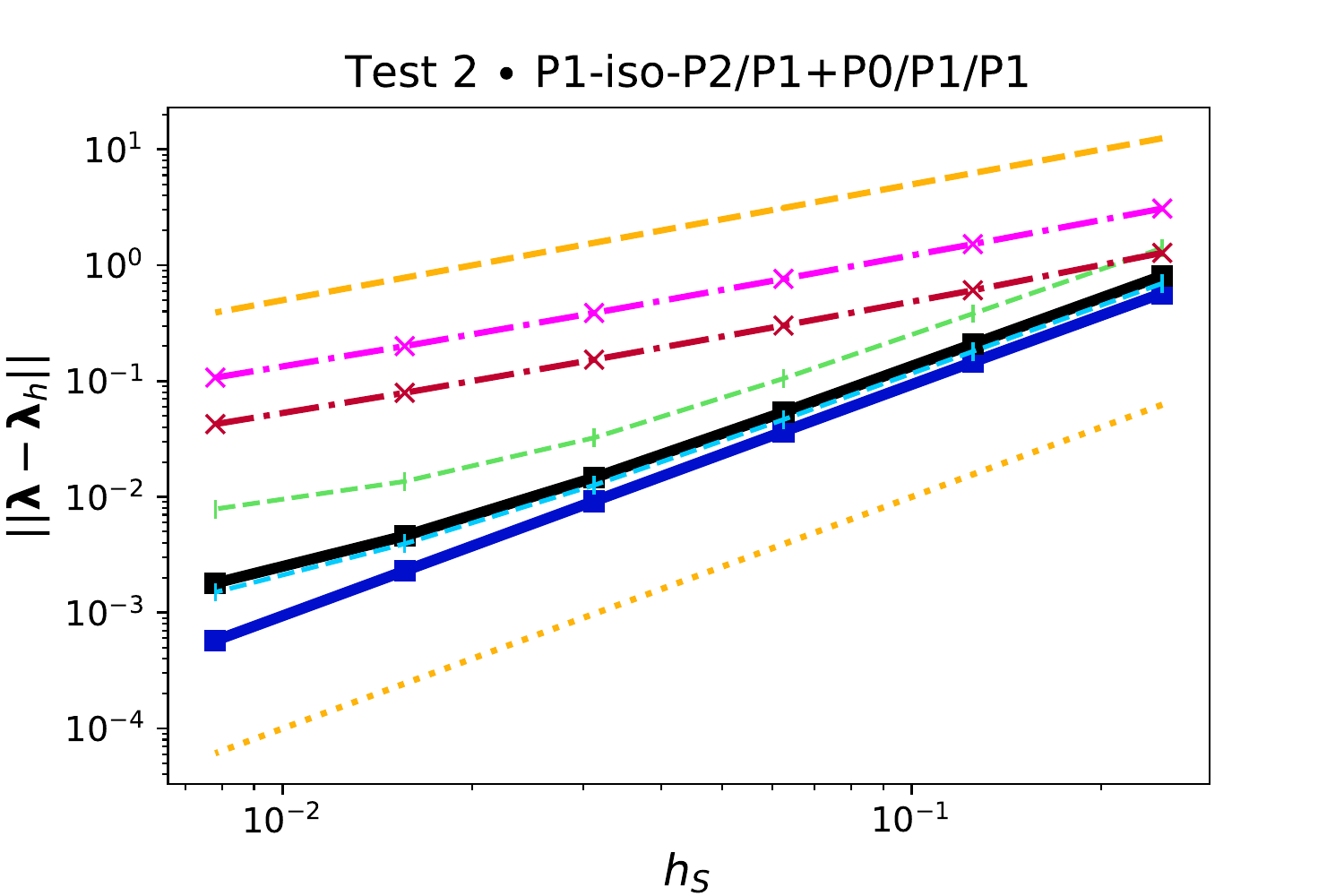}\\
		\includegraphics[trim=150 16 50 50,width=1.4\linewidth]{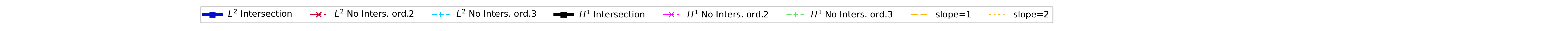}
		\caption{Convergence plots of Test 2 with $\Pcal_1-iso-\Pcal_2/\Pcal_1+P0/\Pcal_1/\Pcal_1$}
		\label{fig:test2_p1p0}
	\end{figure}

	In this test, we consider the case of meshes with coincident edges to make a first comparison between the methods in a situation where the intersection of the meshes does not have a complex geometry. From Figures \ref{fig:test2_p1} and \ref{fig:test2_p1p0}, we can, however, already notice that the convergence rates of the method without intersection with quadrature rule of order two are lower than the convergence rates of the method with intersection. On the other hand, the increase of the quadrature order for the method without intersection produces in this case a significant improvement, almost reaching the good behavior of the method with intersection.
	
	\subsection{Test 3}\label{sub:test3}
\begin{figure}[!h]
		\centering
		\includegraphics[width=0.24\linewidth]{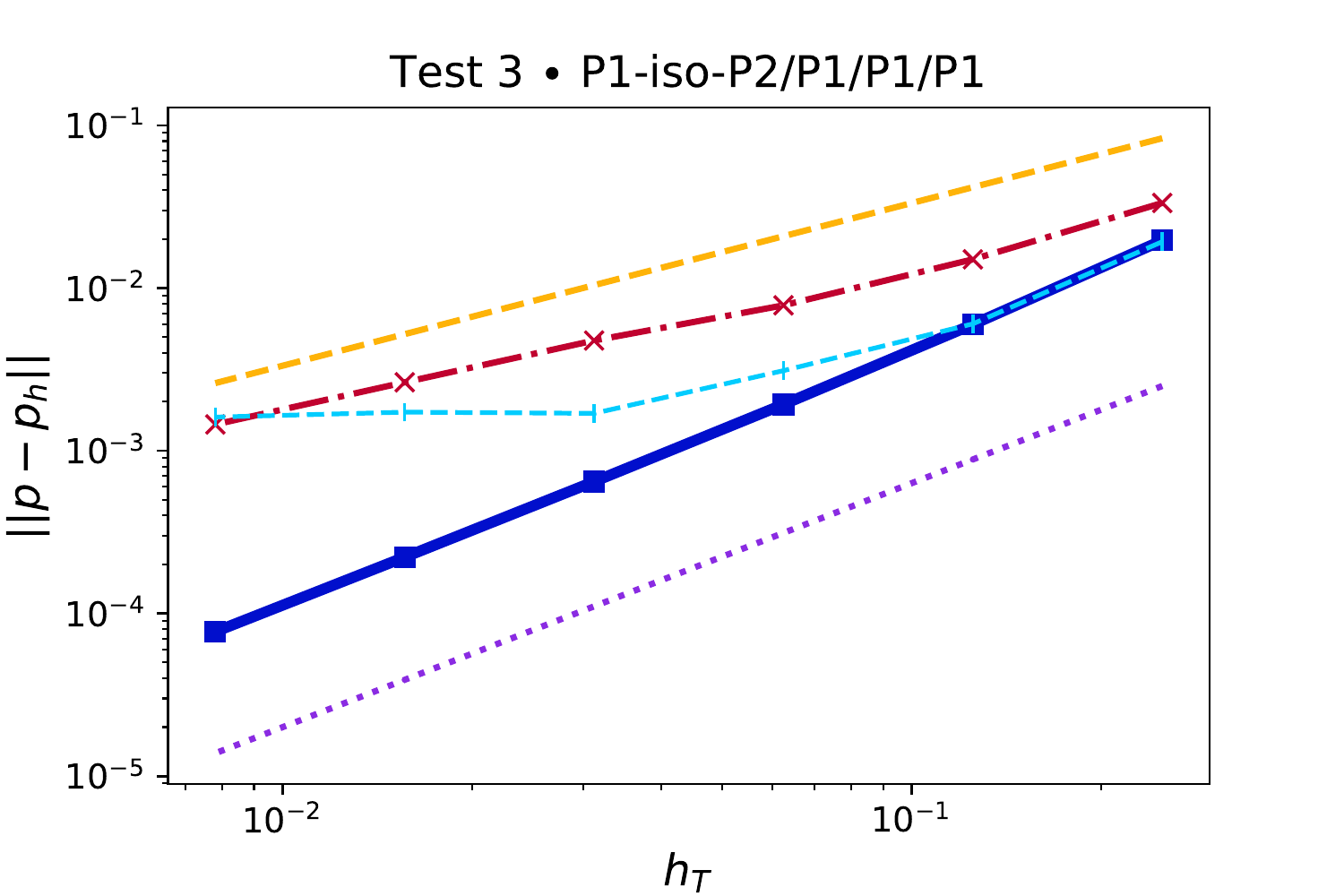}
		\includegraphics[width=0.24\linewidth]{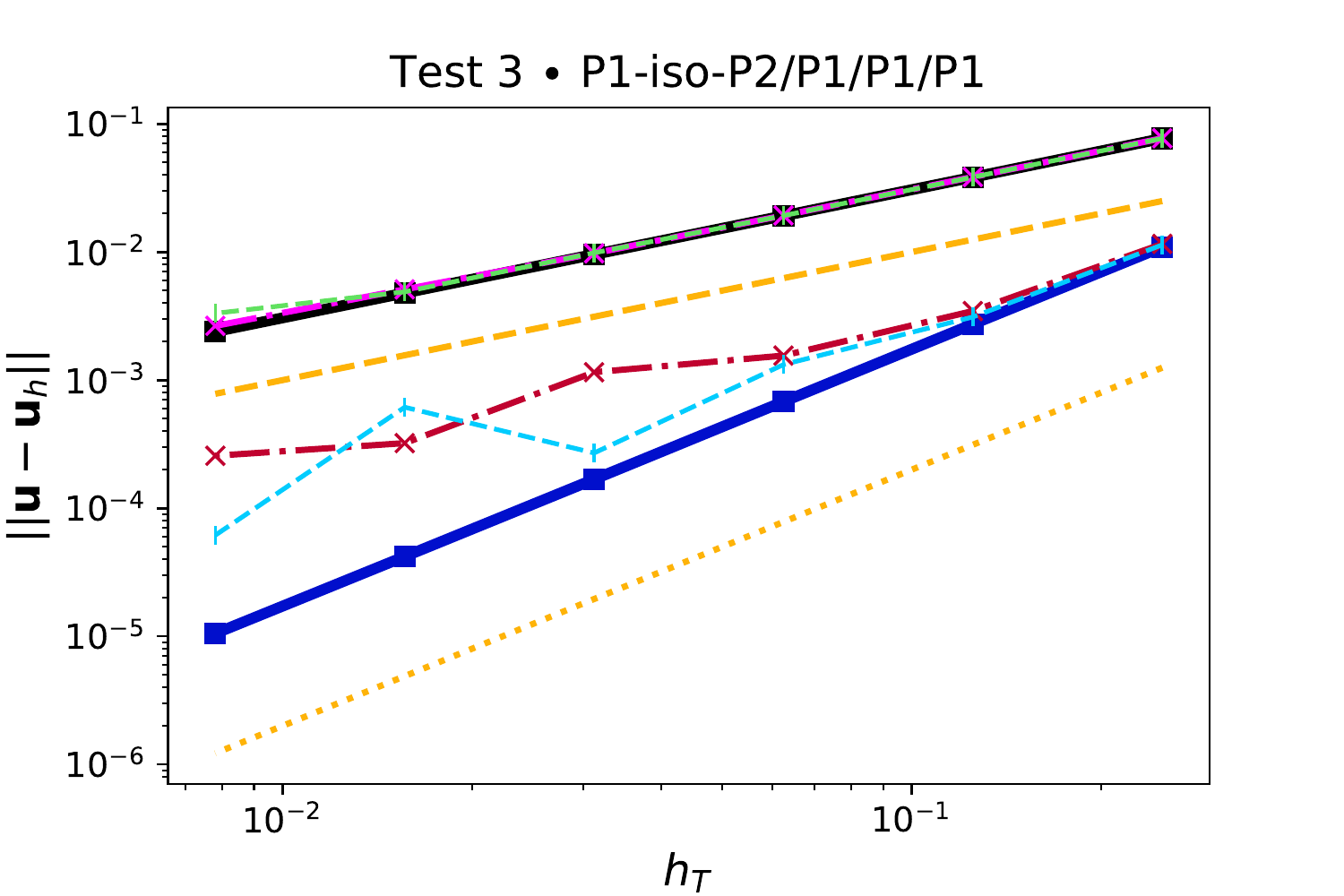}
		\includegraphics[width=0.24\linewidth]{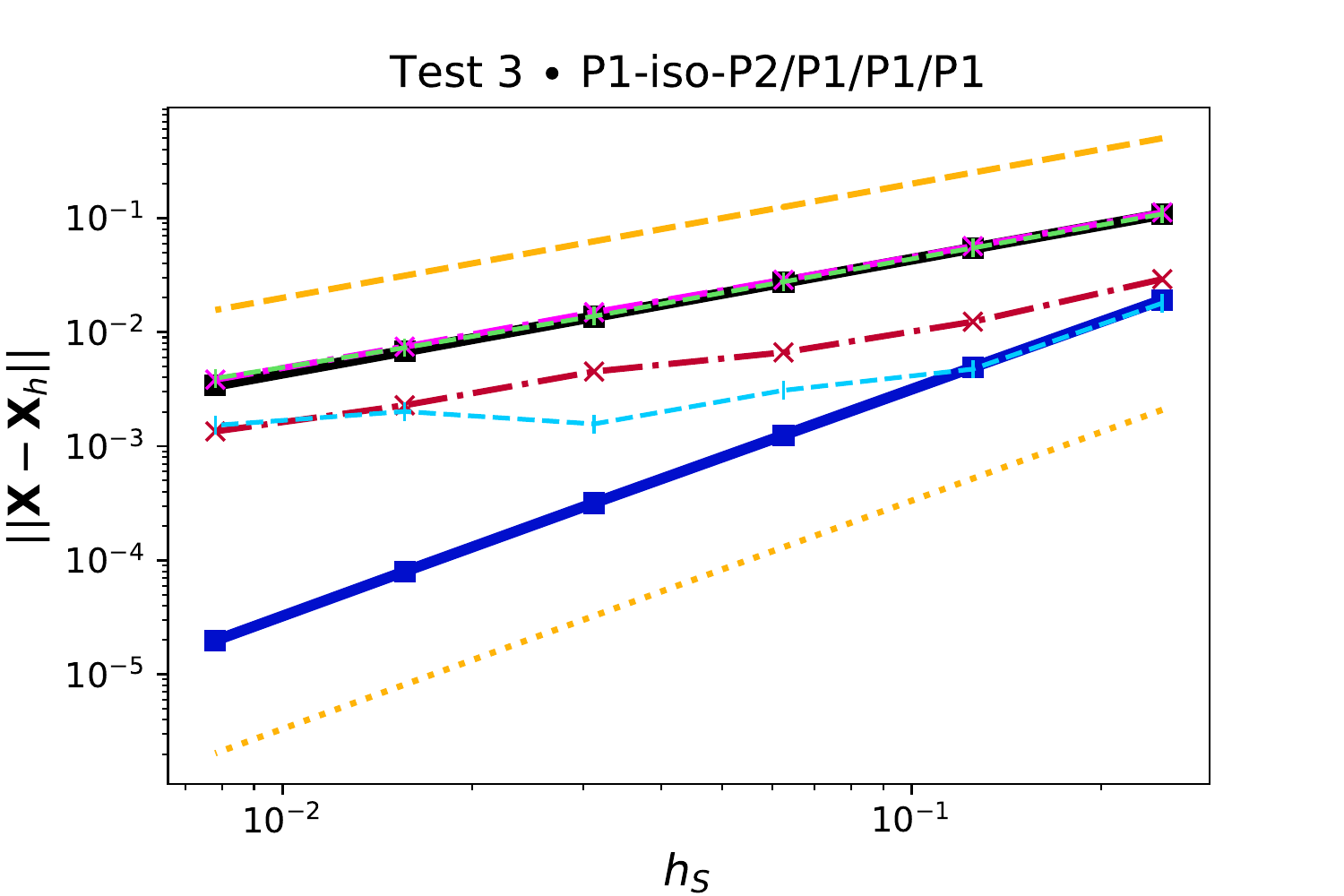}
		\includegraphics[width=0.24\linewidth]{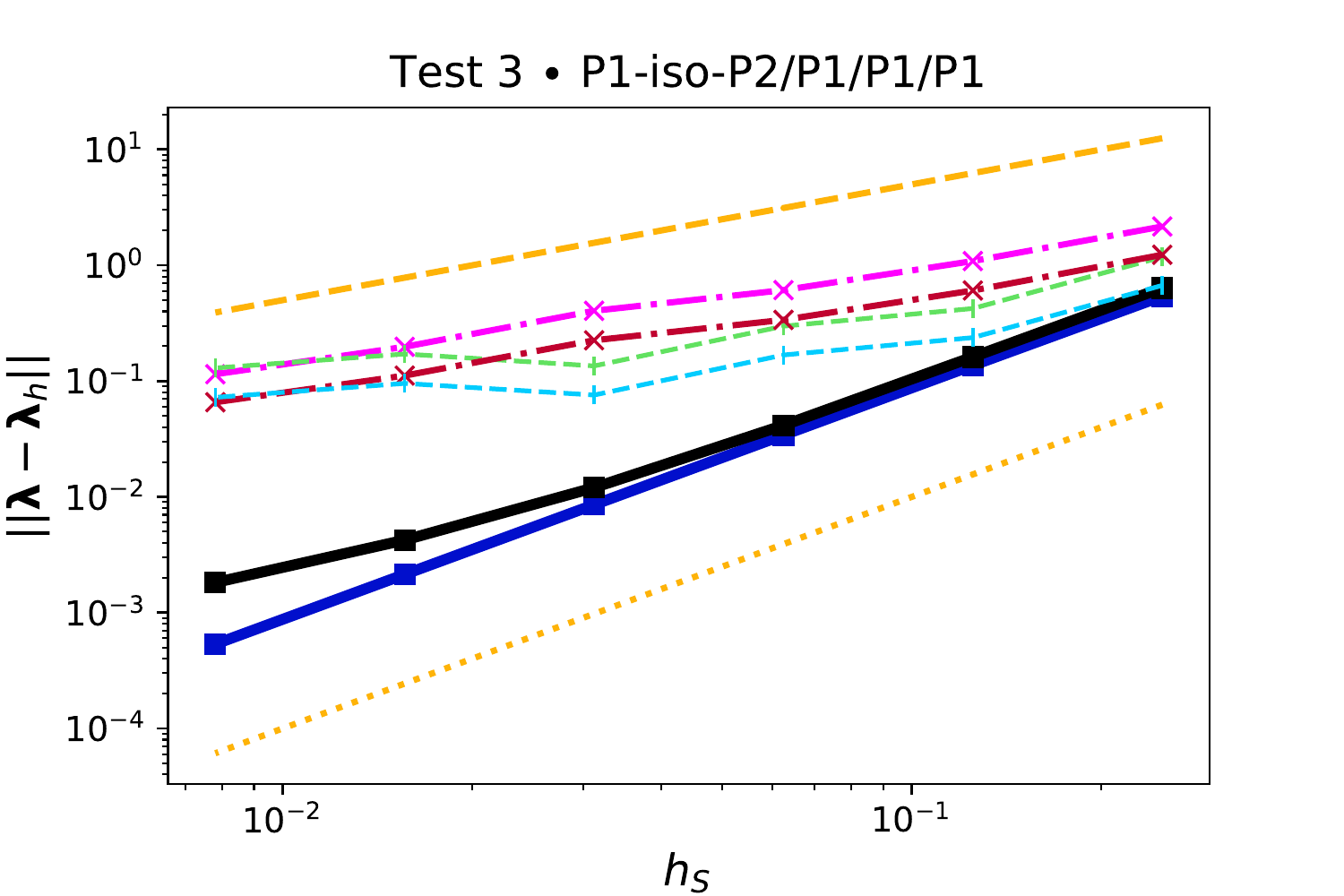}\\
		\includegraphics[trim=180 16 50 50,width=1.35\linewidth]{figures/label_p1-eps-converted-to}
		\caption{Convergence plots of Test 3 with $\Pcal_1-iso-\Pcal_2/\Pcal_1/\Pcal_1/\Pcal_1$}
		\label{fig:test3_p1}
	\end{figure}

\begin{figure}[!h]
		\centering
		\includegraphics[width=0.24\linewidth]{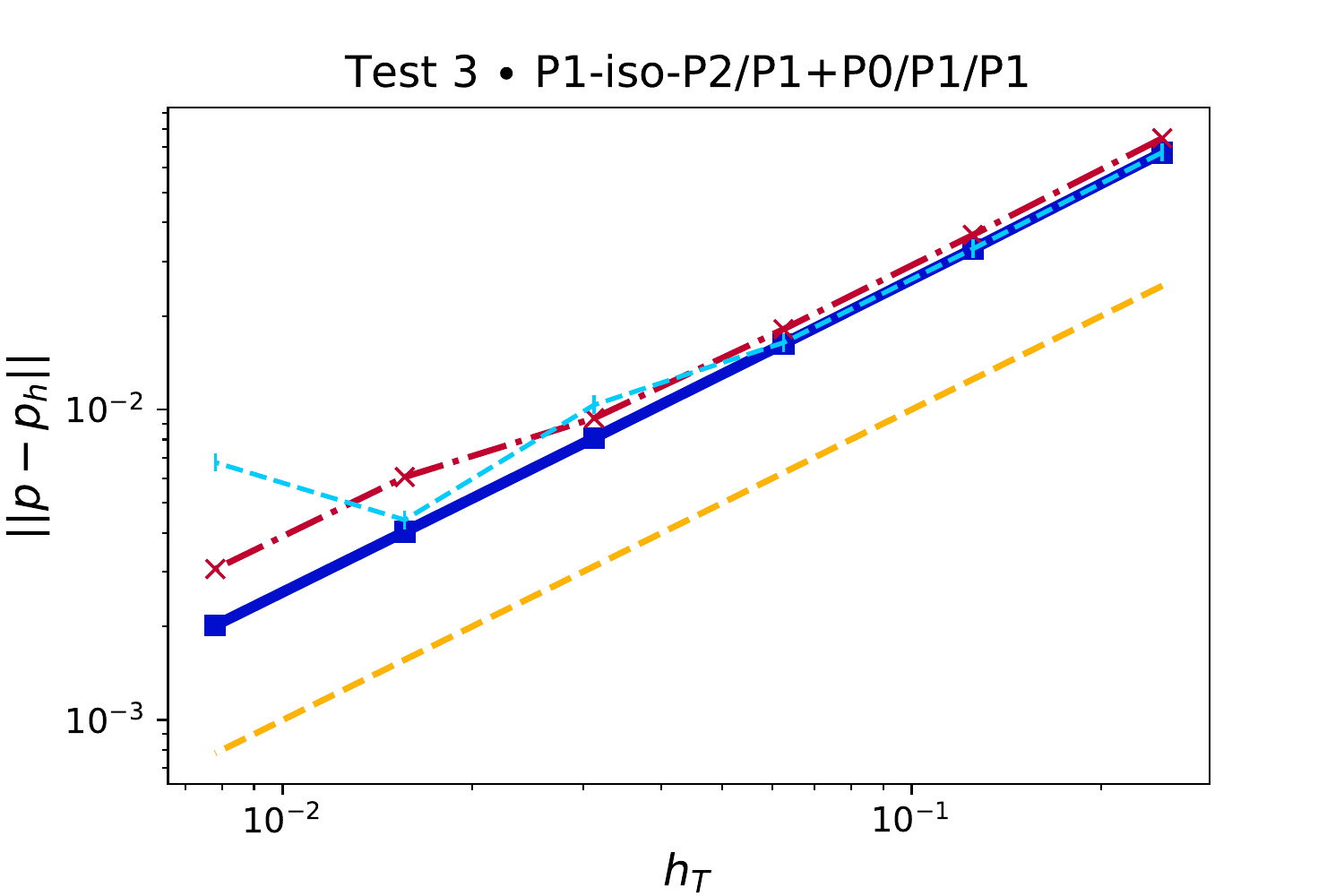}
		\includegraphics[width=0.24\linewidth]{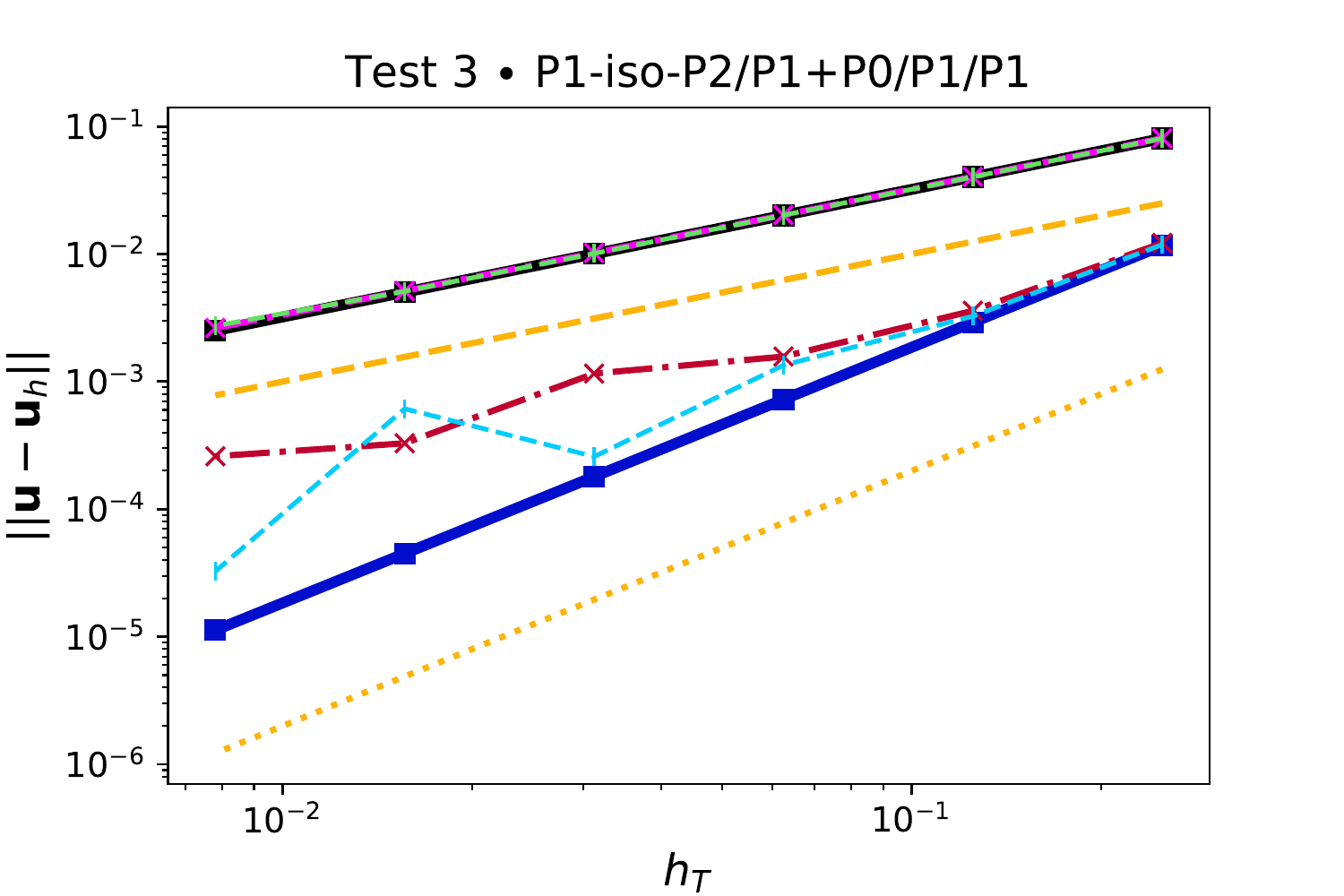}
		\includegraphics[width=0.24\linewidth]{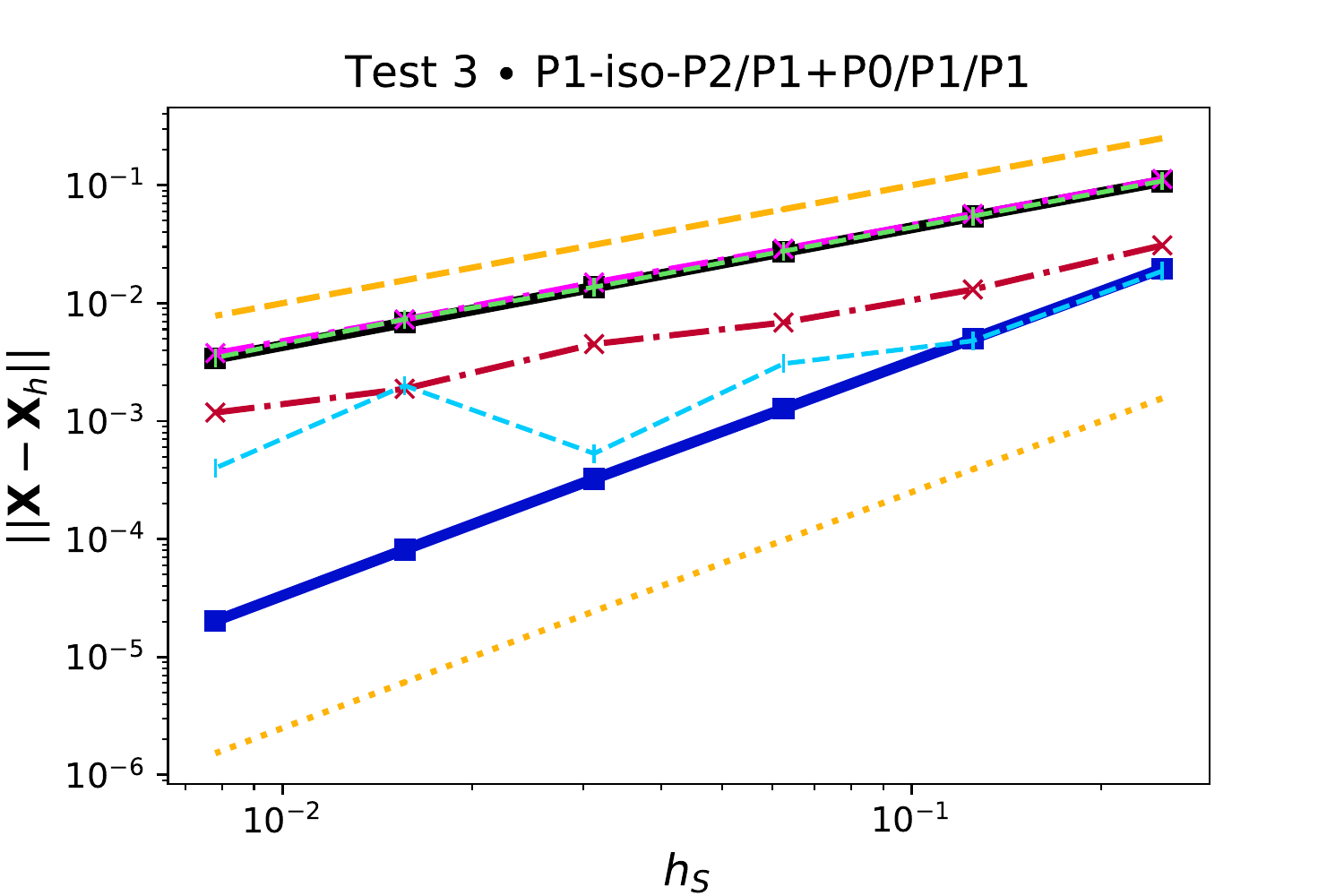}
		\includegraphics[width=0.24\linewidth]{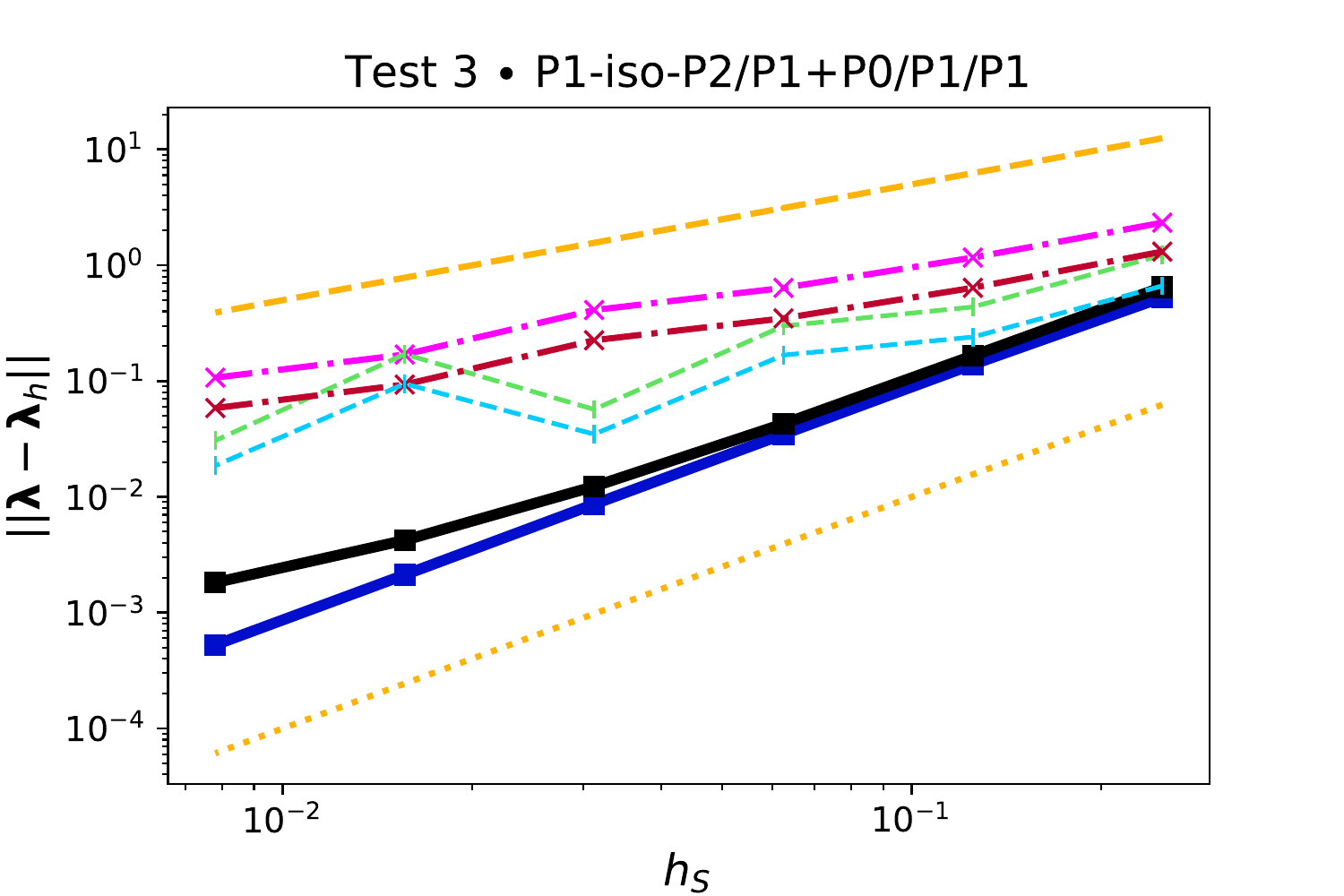}\\
		\includegraphics[trim=150 16 50 50,width=1.4\linewidth]{figures/label_p0-eps-converted-to}
		\caption{Convergence plots of Test 3 with $\Pcal_1-iso-\Pcal_2/\Pcal_1+\Pcal_0/\Pcal_1/\Pcal_1$}
		\label{fig:test3_p1p0}
	\end{figure}

	For this third test, we use again the solutions chosen in \eqref{eq:standard_sol}. On the other hand, we set $\Omega=[-2,2]^2$ and $\B=\Os=[-0.62,1.38]^2$ so that the interface of the solid does not match with the fluid mesh, discretizing them with uniform meshes as before. The results are reported in Figures \ref{fig:test3_p1} and \ref{fig:test3_p1p0}. Thanks to the new positioning of the structure, without matching boundaries, we get clearer results.  The only method capable of ensuring good convergence is the one that computes the intersection. In this case, the increase in precision of the quadrature rule used when the intersection is not calculated does not produce any improvement.
	
	\subsection{Test 4}\label{sub:test4}
\begin{figure}[!h]
		\centering
		\includegraphics[width=0.24\linewidth]{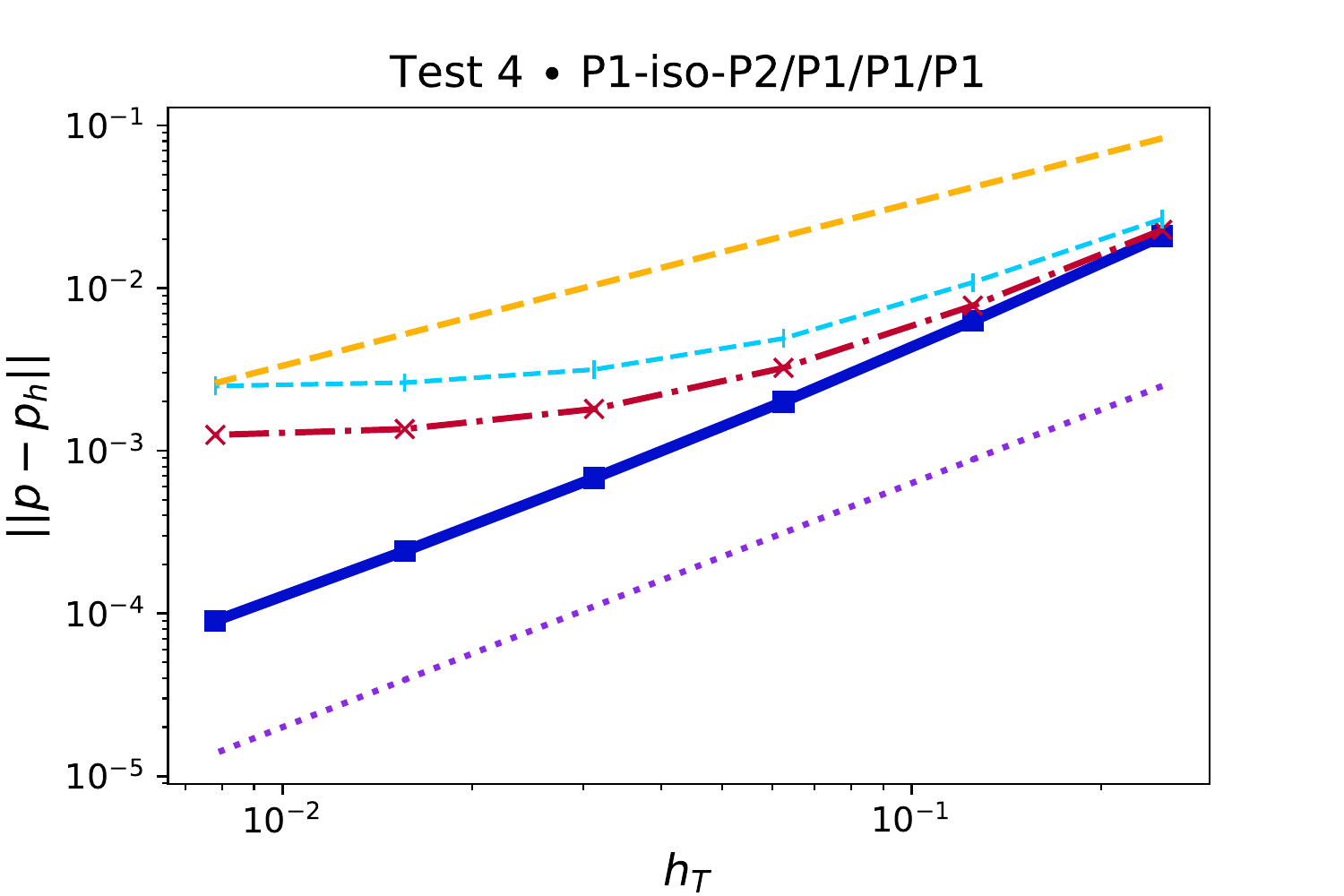}
		\includegraphics[width=0.24\linewidth]{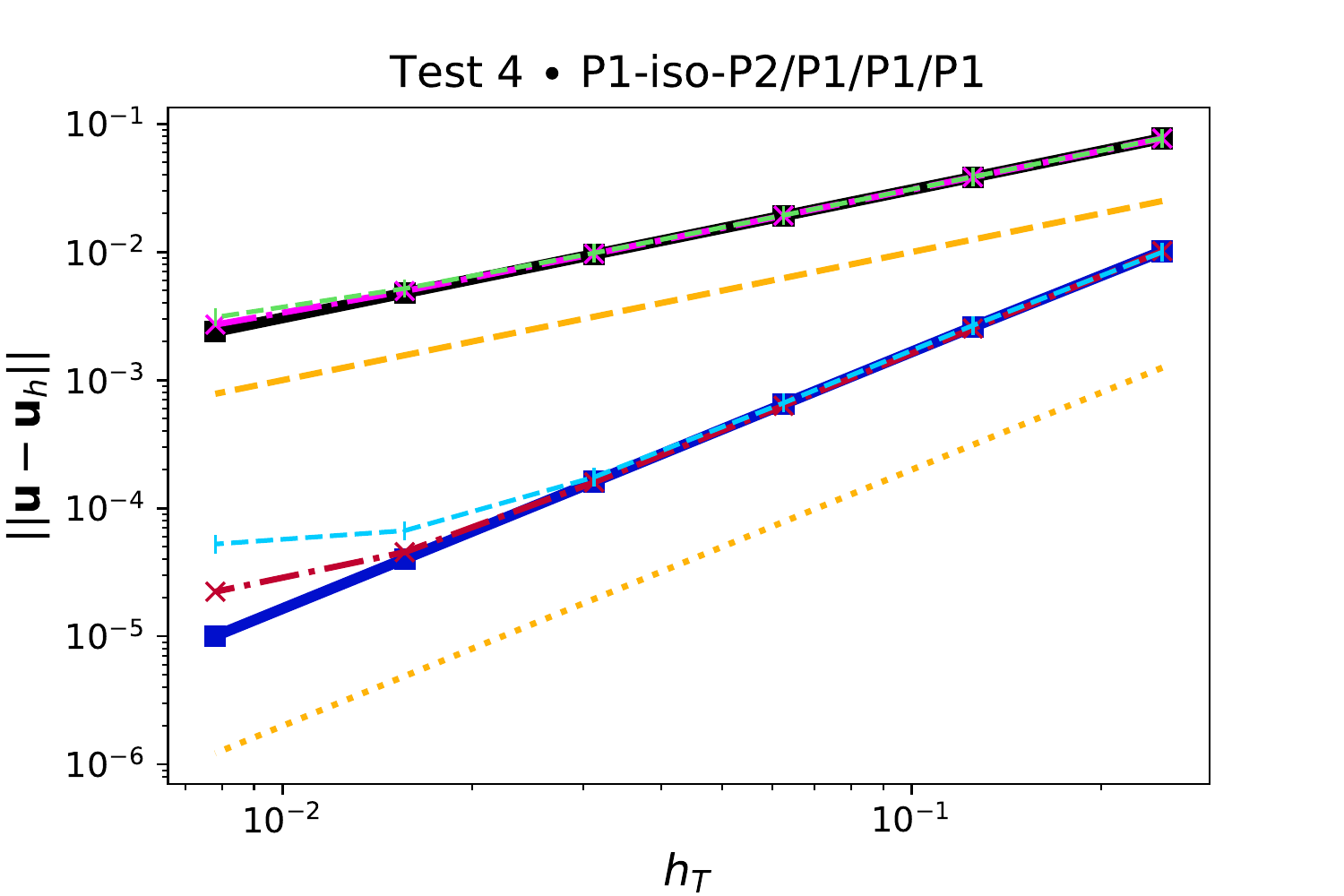}
		\includegraphics[width=0.24\linewidth]{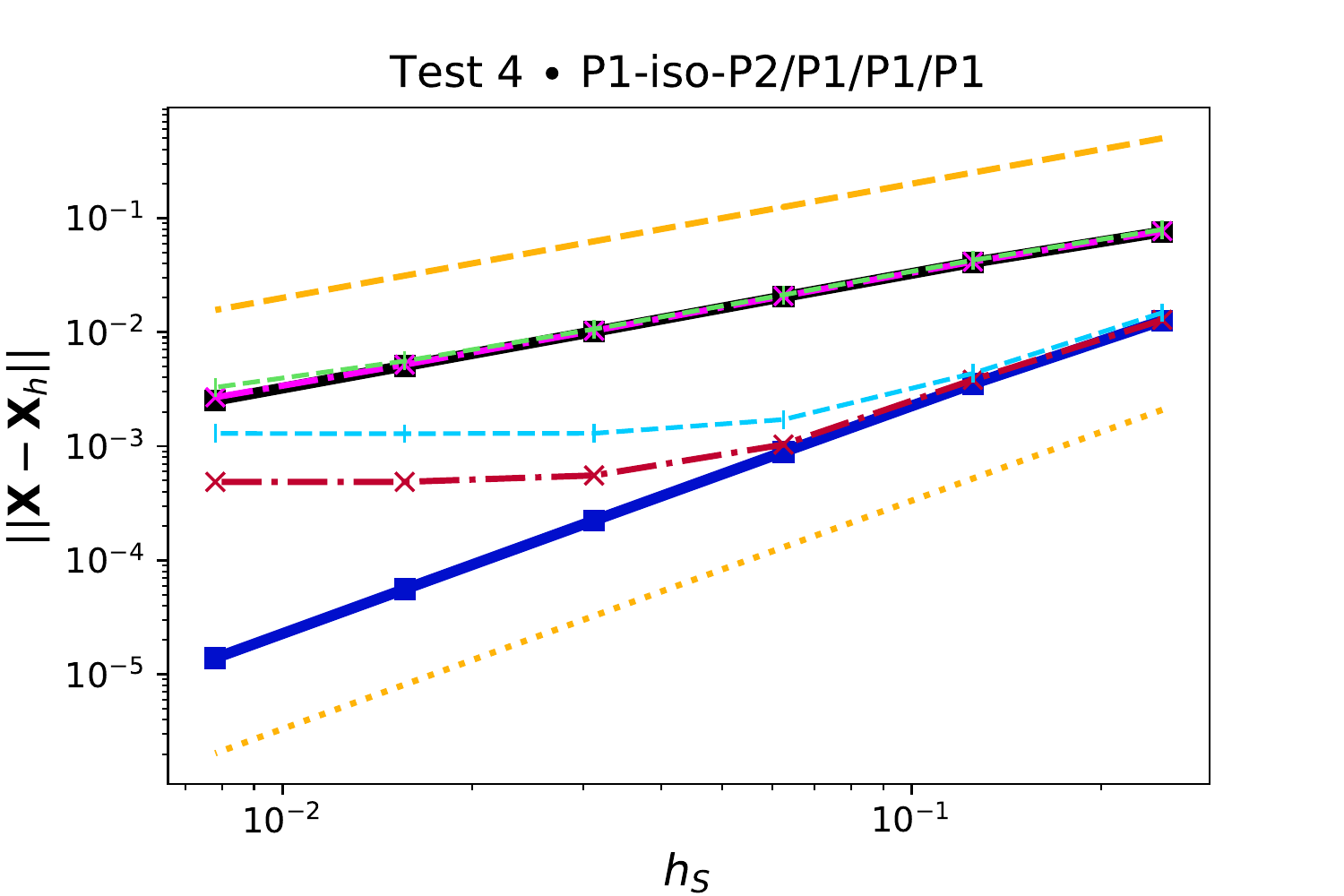}
		\includegraphics[width=0.24\linewidth]{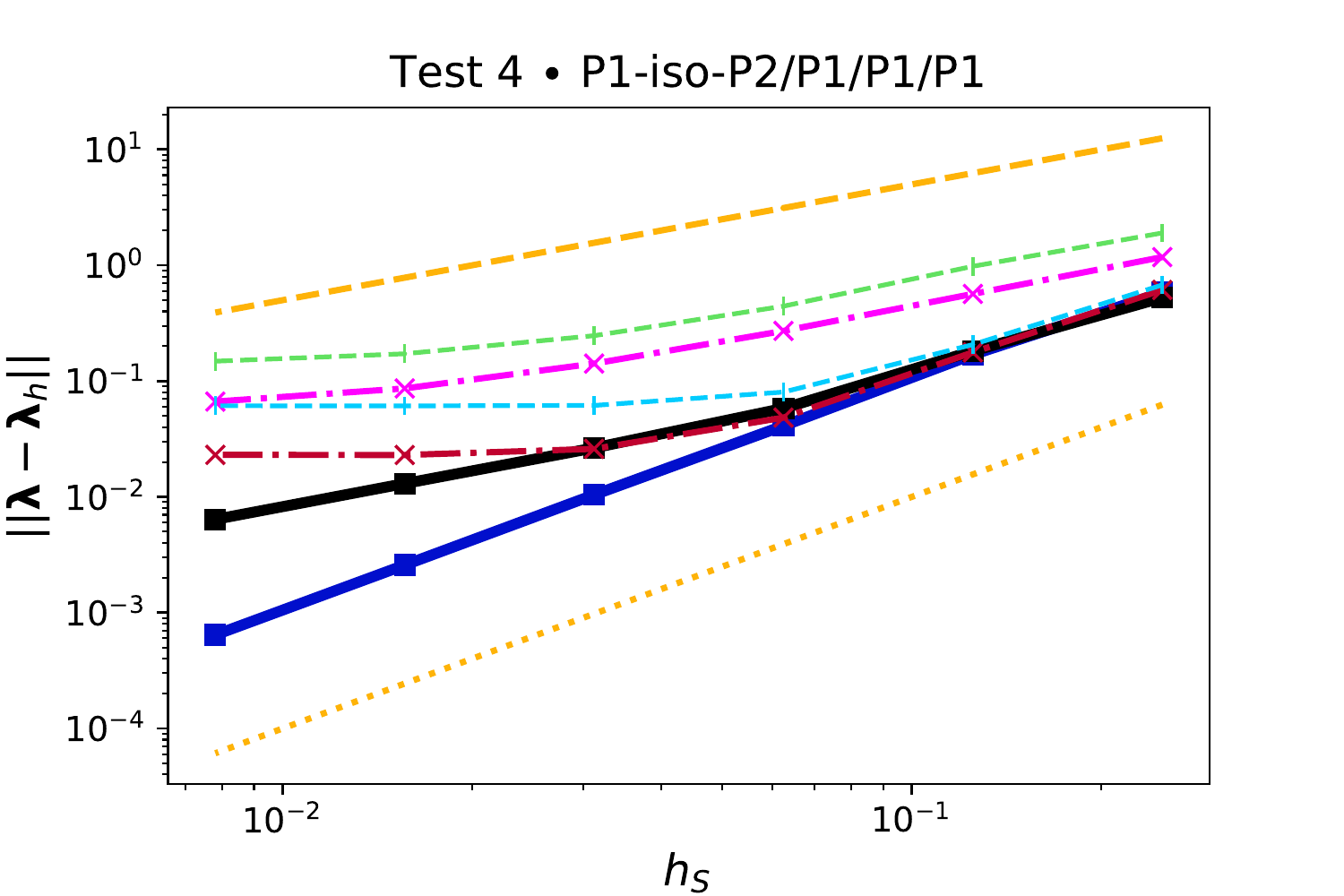}\\
		\includegraphics[trim=180 16 50 50,width=1.35\linewidth]{figures/label_p1-eps-converted-to}
		\caption{Convergence plots of Test 4 with $\Pcal_1-iso-\Pcal_2/\Pcal_1/\Pcal_1/\Pcal_1$}
		\label{fig:test4_p1}
	\end{figure}

\begin{figure}[!h]
		\centering
		\includegraphics[width=0.24\linewidth]{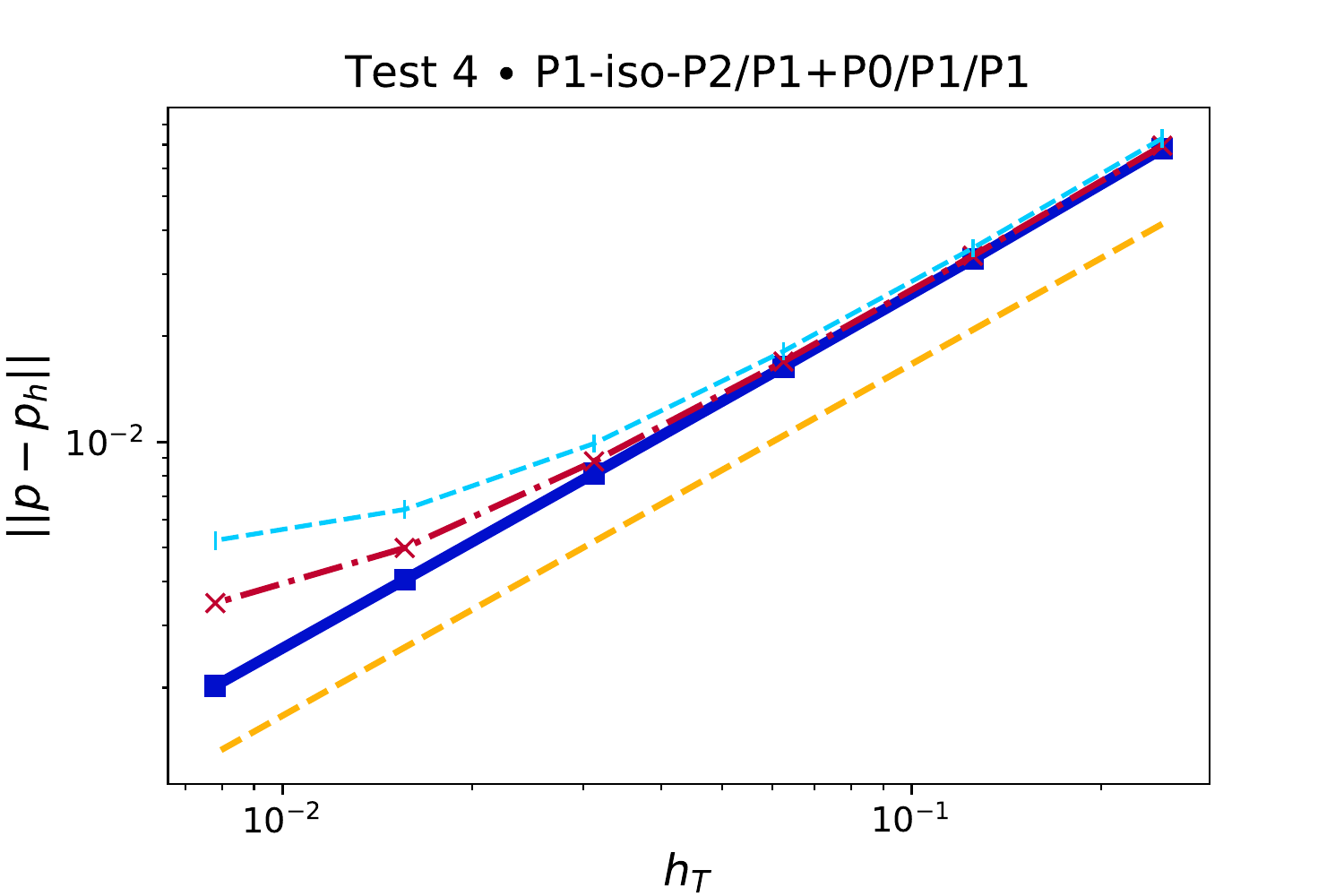}
		\includegraphics[width=0.24\linewidth]{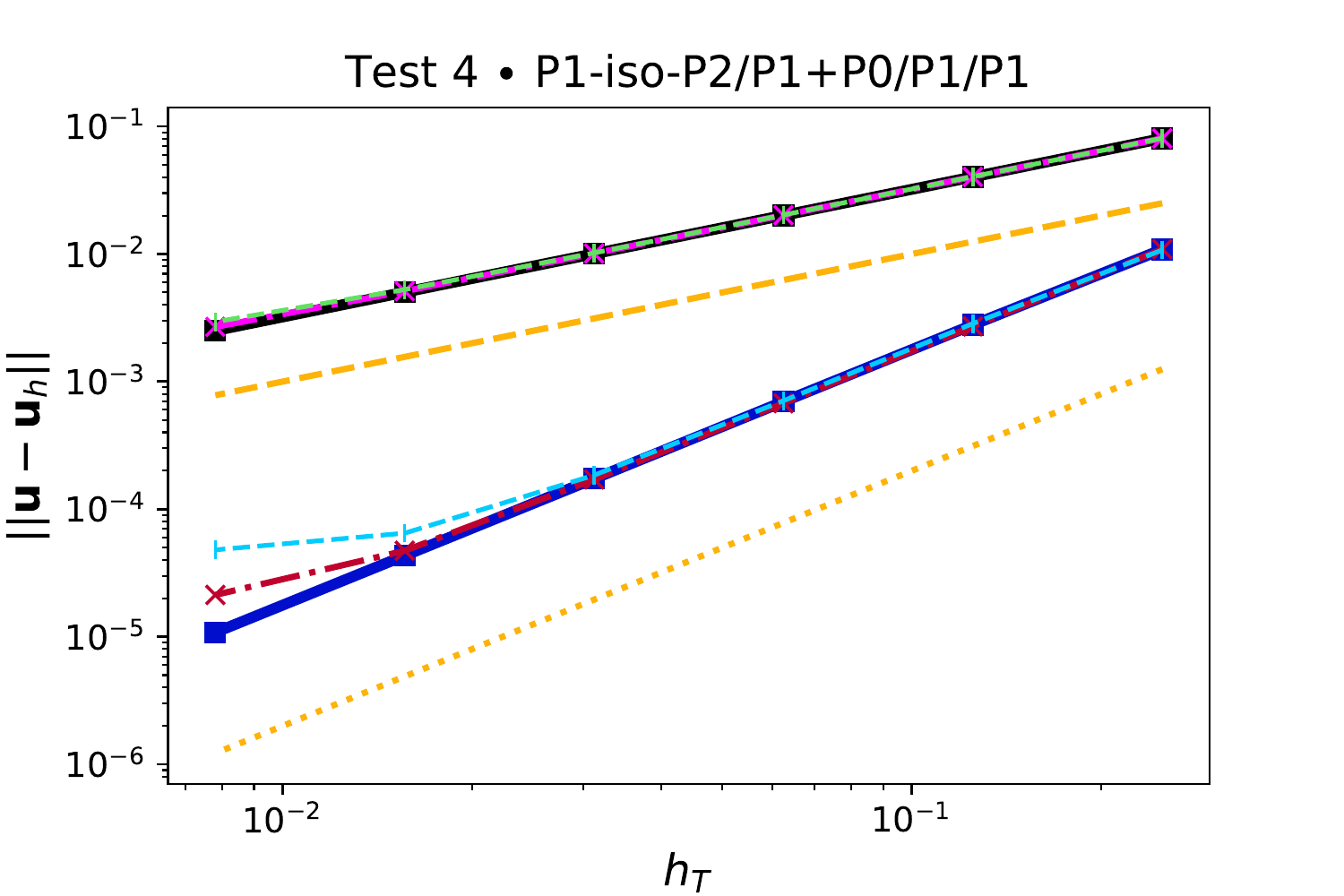}
		\includegraphics[width=0.24\linewidth]{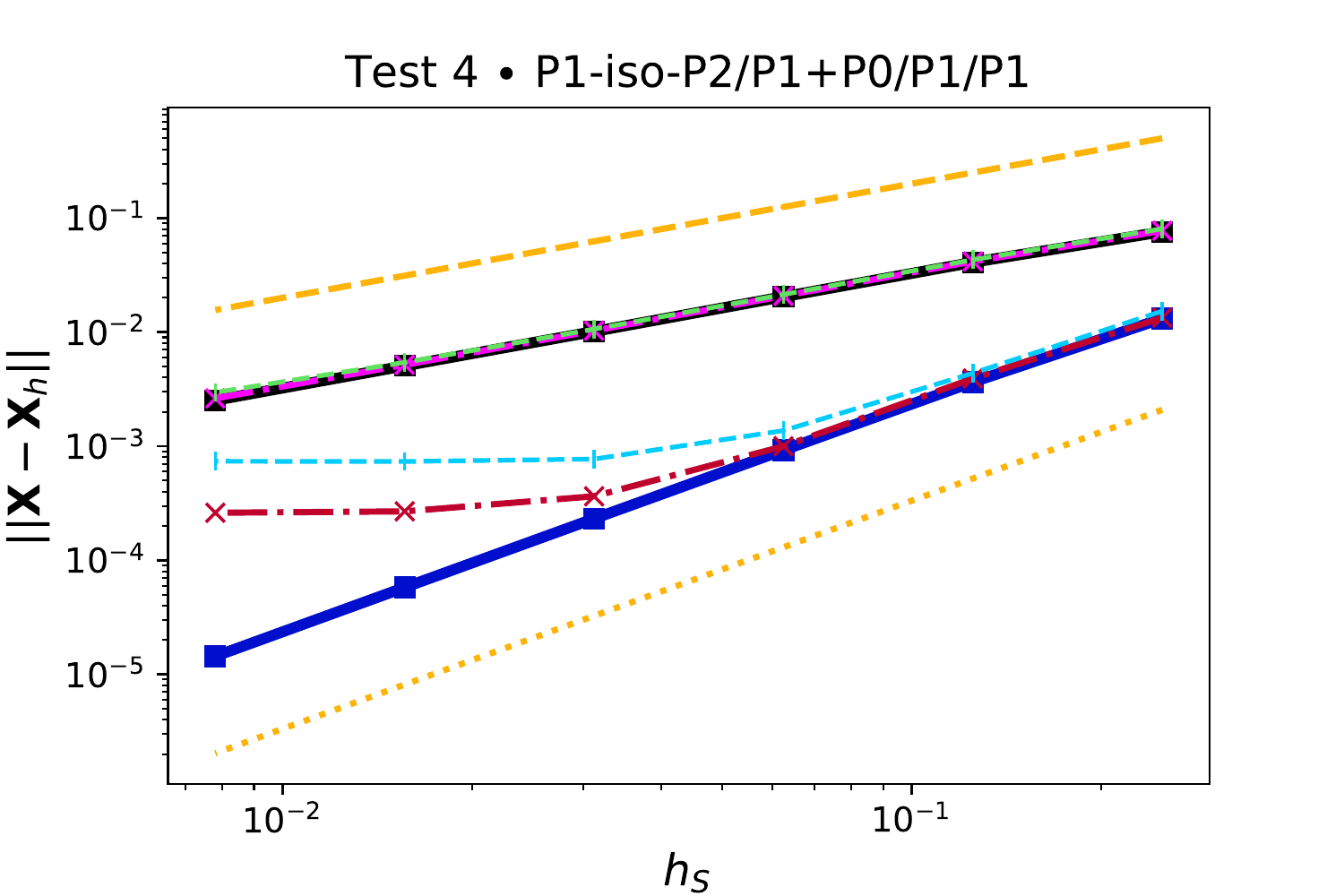}
		\includegraphics[width=0.24\linewidth]{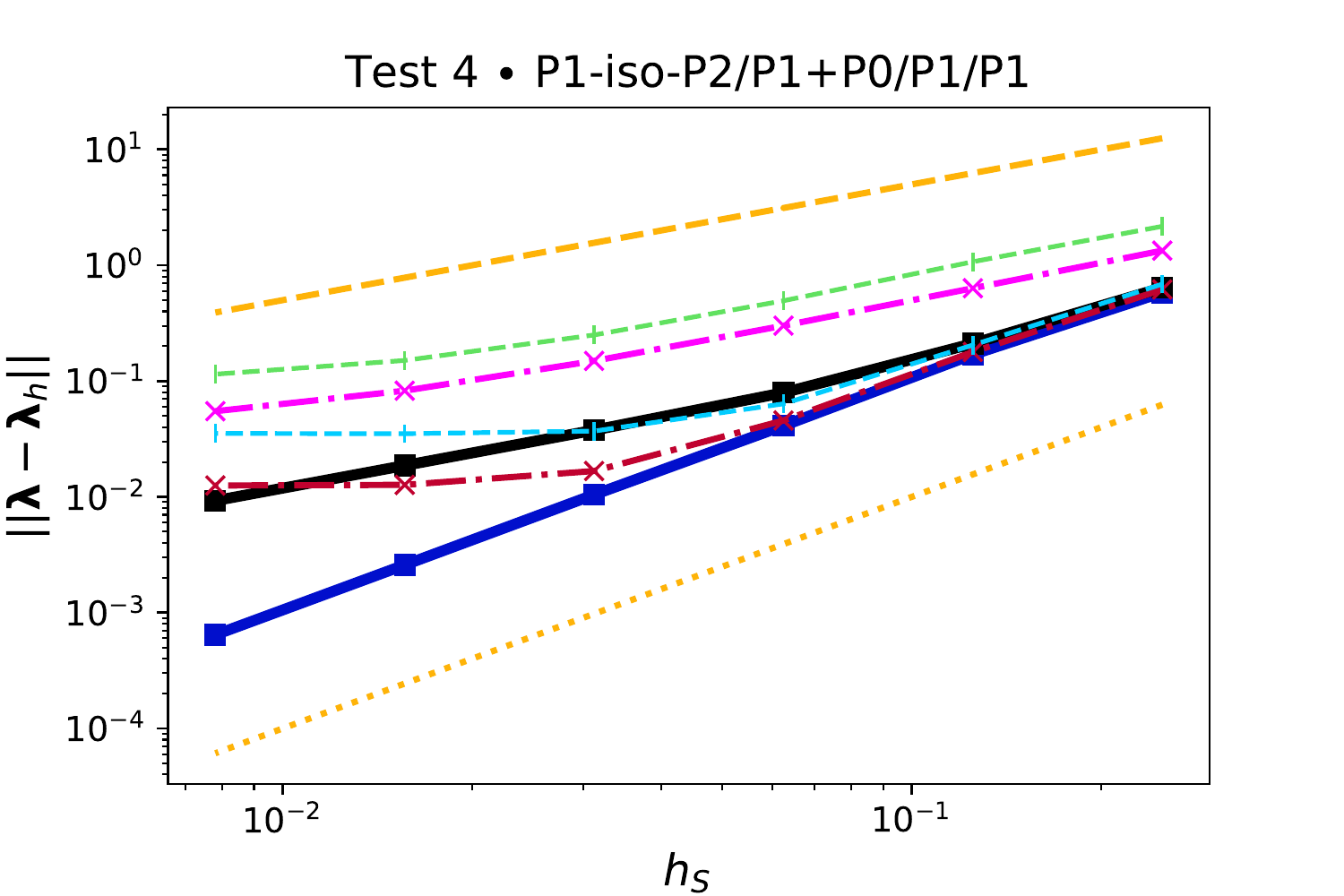}\\
		\includegraphics[trim=150 16 50 50,width=1.4\linewidth]{figures/label_p0-eps-converted-to}
		\caption{Convergence plots of Test 4 with $\Pcal_1-iso-\Pcal_2/\Pcal_1+\Pcal_0/\Pcal_1/\Pcal_1$}
		\label{fig:test4_p1p0}
	\end{figure}

	In the same setting of the previous test, we partition the solid domain making use of an unstructured grid \blue(Figure \ref{fig:unstr})\noblue. The use of this type of mesh for the discretization of the immersed body places us in an even more general setting. Figures \ref{fig:test4_p1} and \ref{fig:test4_p1p0} show that the bad performances seen in the previous tests for the methods without intersection are amplified; \blue in the case of the pressure, the non-intersection methods perform better when we use discontinuous elements\noblue. Furthermore, if in the previous tests $\|\llambda-\llambda_h\|_{1,\B}$ showed a convergence rate higher than $1.5$, this does not happen now, settling around $1$.
	
	\subsection{Test 5}\label{sub:test5}
	\lg The configuration of this test is the same we chose for Test~1, but, in this case, we \gl\blue consider a pressure which is discontinuous on $\partial\Os$\noblue.  We assume that $\Omega=[-2,2]^2$ and $\B=\Os=[-1,1]^2$: as before, for the \lg fluid domain \gl we choose a right-uniform mesh, while for the \lg structure \gl a left-uniform one. The fact that the reference and the actual solid domain coincide implies that $\Xbar=\mathbf{I}_\B$. We compute the right hand sides $\f$, $\g$ in order to obtain the following solutions:
	\begin{equation}
		\begin{aligned}
			&\u(x,y)=\curl\big((4-x^2)^2(4-y^2)^2\big)\\
			&\u_{\mid\partial\Omega}=\mathbf{0}\\
			&p(x,y) = 
			\begin{cases}
				150\sin(x)-\frac{50}{3}&\text{ in }\Omega^f\\
				150\sin(x)+50&\text{ in }\Os=\B.
			\end{cases}\\
			&\X(x,y)=\u(x,y)\\
			&\llambda(x,y)=\big(e^x,e^y\big).
		\end{aligned}
	\end{equation}
	Here, the choice of $p$ implies that we have to introduce a new weak term in the definition of the right hand side $\f$ taking into account the \blue discontinuity\noblue. Indeed, we have
	\begin{equation*}
		\int_\Omega p\, \div\v\,\dx = \int_{\Os} \v\cdot\grad p \,\dx - \int_{\partial\Os} (\v\cdot\n_s)\,p\,\dA + \int_{\Omega\smallsetminus\Os} \v\cdot\grad p \,\dx - \int_{\partial\Os} (\v\cdot\n_f)\,p\,\dA;
	\end{equation*}
	therefore, in order to compute the boundary integrals combining $\partial\Os$ with the basis for $\V_h$, we have to take care of the immersion of the solid boundary in the fluid mesh. \blue We particularly emphasize that we are in the situation where the discontinuity of the pressure not only matches the interface, but also the velocity mesh\noblue. In this situation, we can see how the choice of element for pressure affects the convergence of all the variables. As expected and presented in \cite{mass}, the enhanced space $\Pcal_1+\Pcal_0$ allows us to reach the optimal convergence without presenting the Gibbs phenomenon that we get using the classical $\Pcal_1$ pressure.
	
\begin{figure}[!h]
		\centering
		\includegraphics[width=0.24\linewidth]{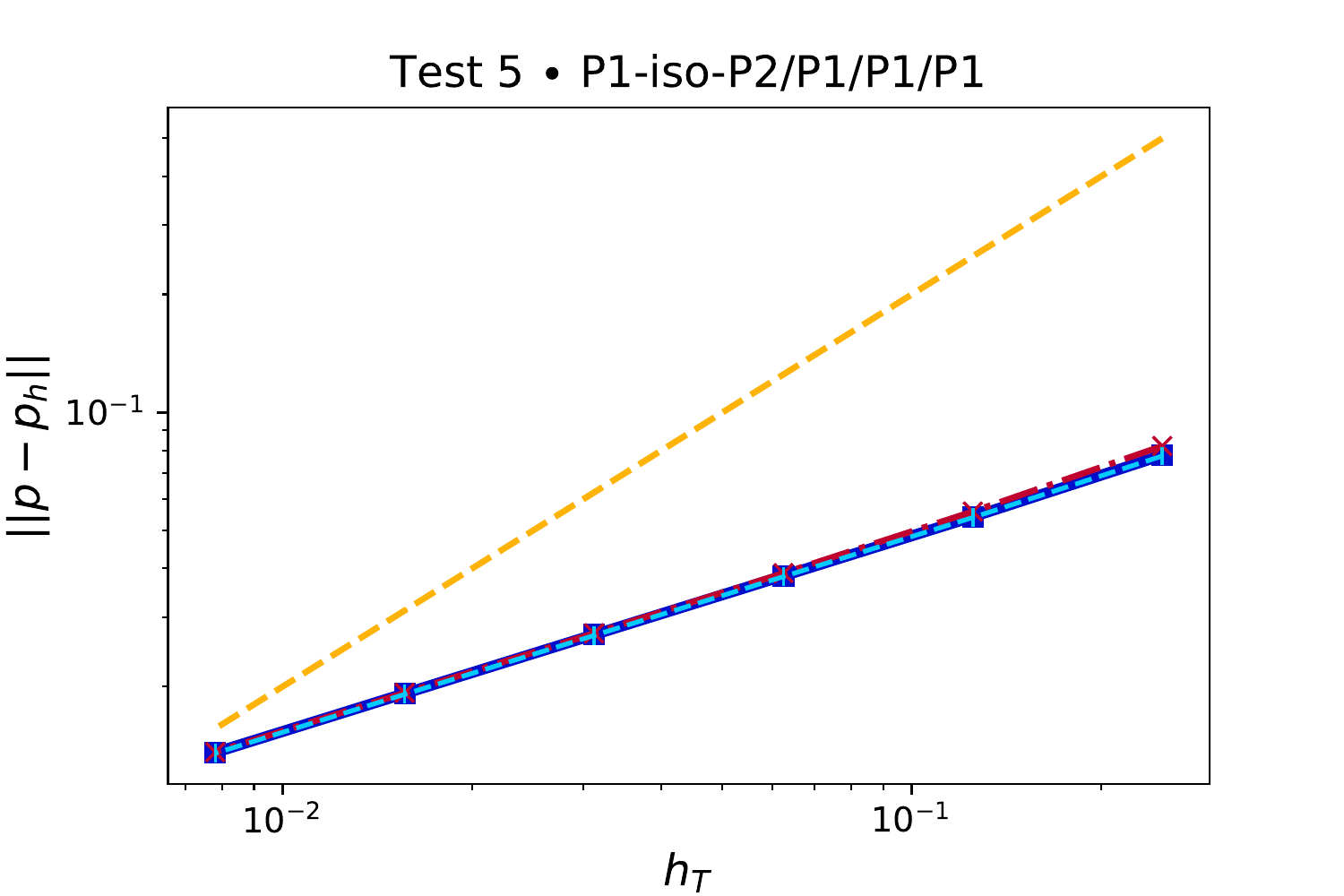}
		\includegraphics[width=0.24\linewidth]{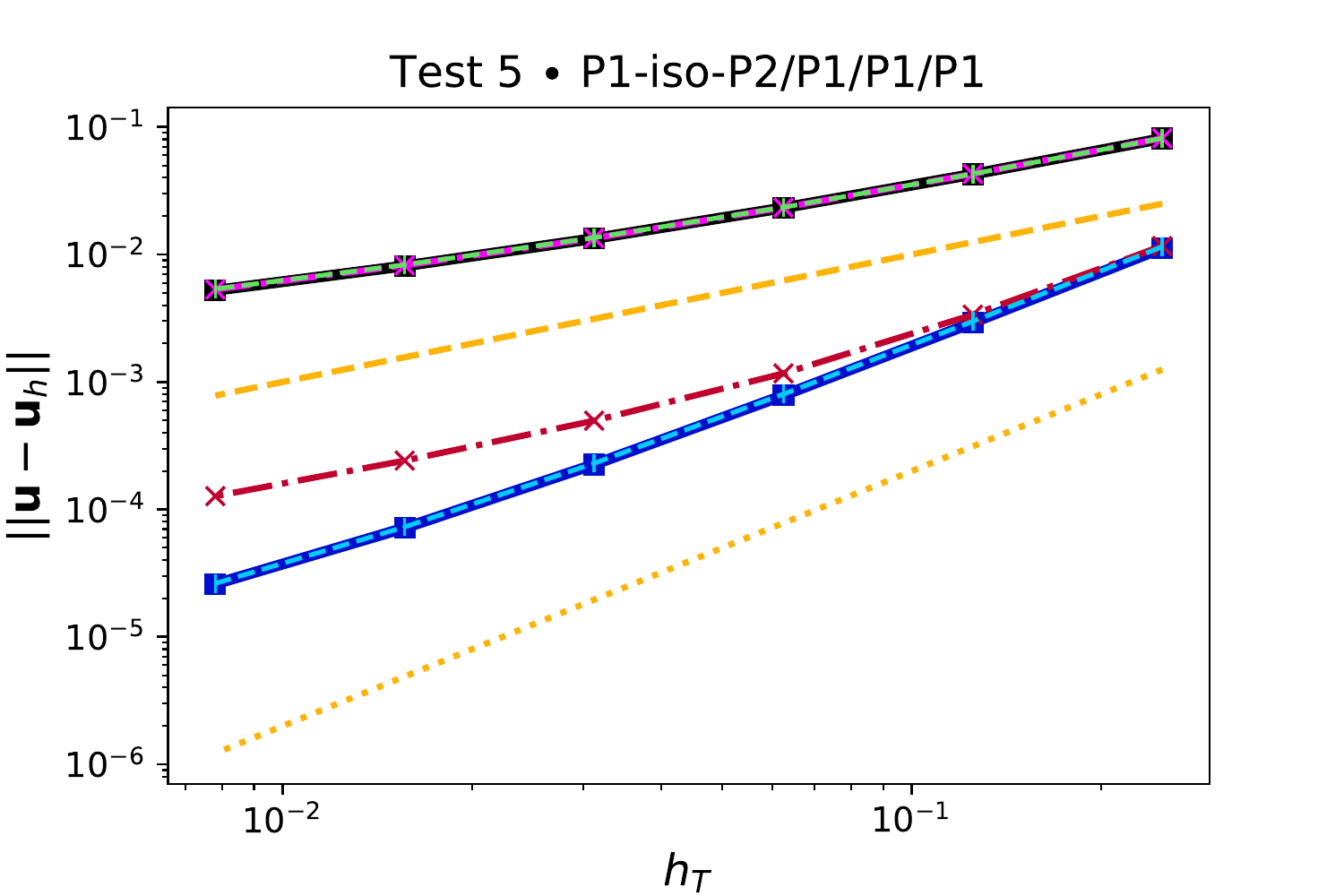}
		\includegraphics[width=0.24\linewidth]{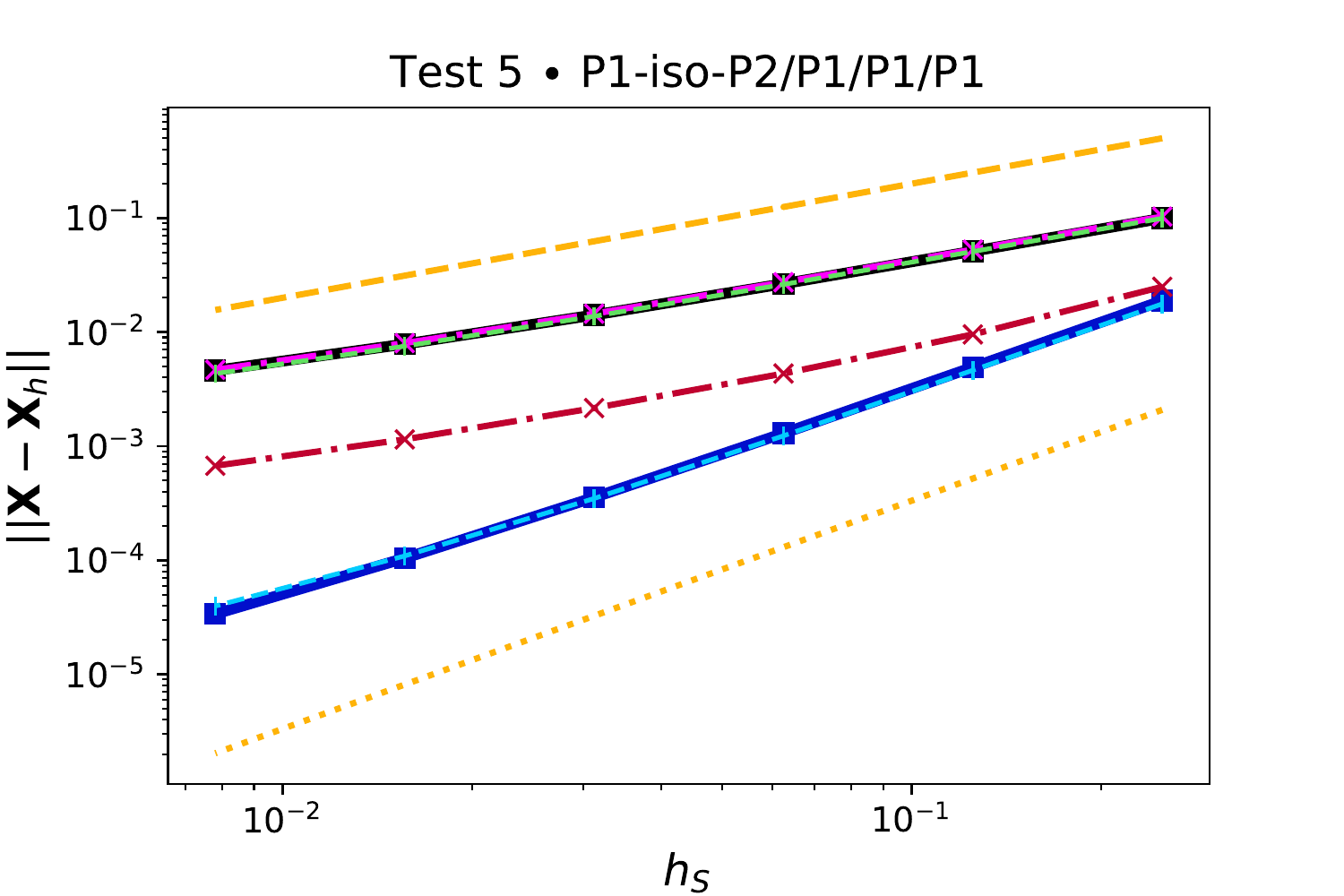}
		\includegraphics[width=0.24\linewidth]{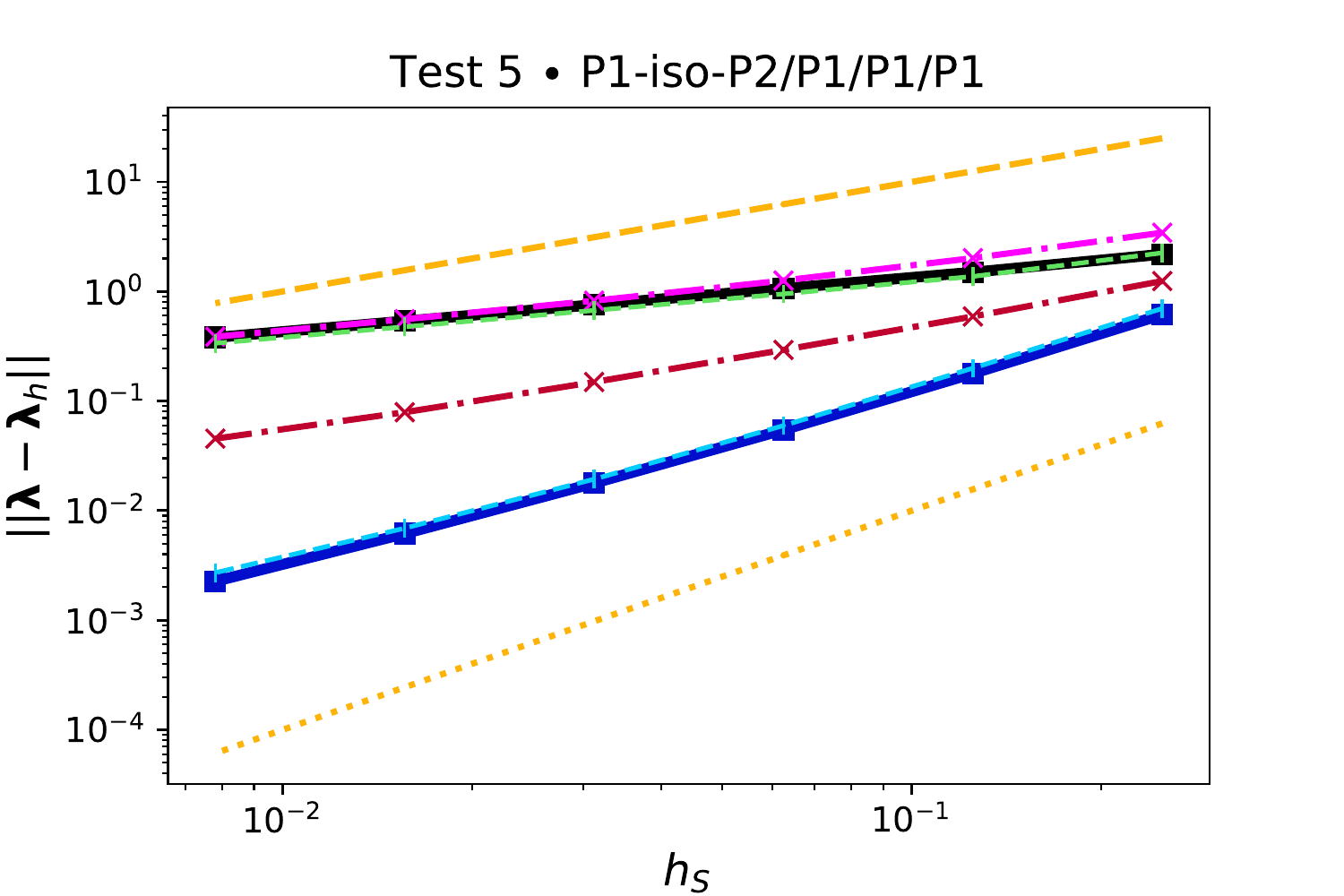}\\
		\includegraphics[trim=180 16 50 50,width=1.35\linewidth]{figures/label_p1-eps-converted-to}
		\caption{Convergence plots of Test 5 with $\Pcal_1-iso-\Pcal_2/\Pcal_1/\Pcal_1/\Pcal_1$}
		\label{fig:test5_p1}
	\end{figure}

\begin{figure}[!h]
		\centering
		\includegraphics[width=0.24\linewidth]{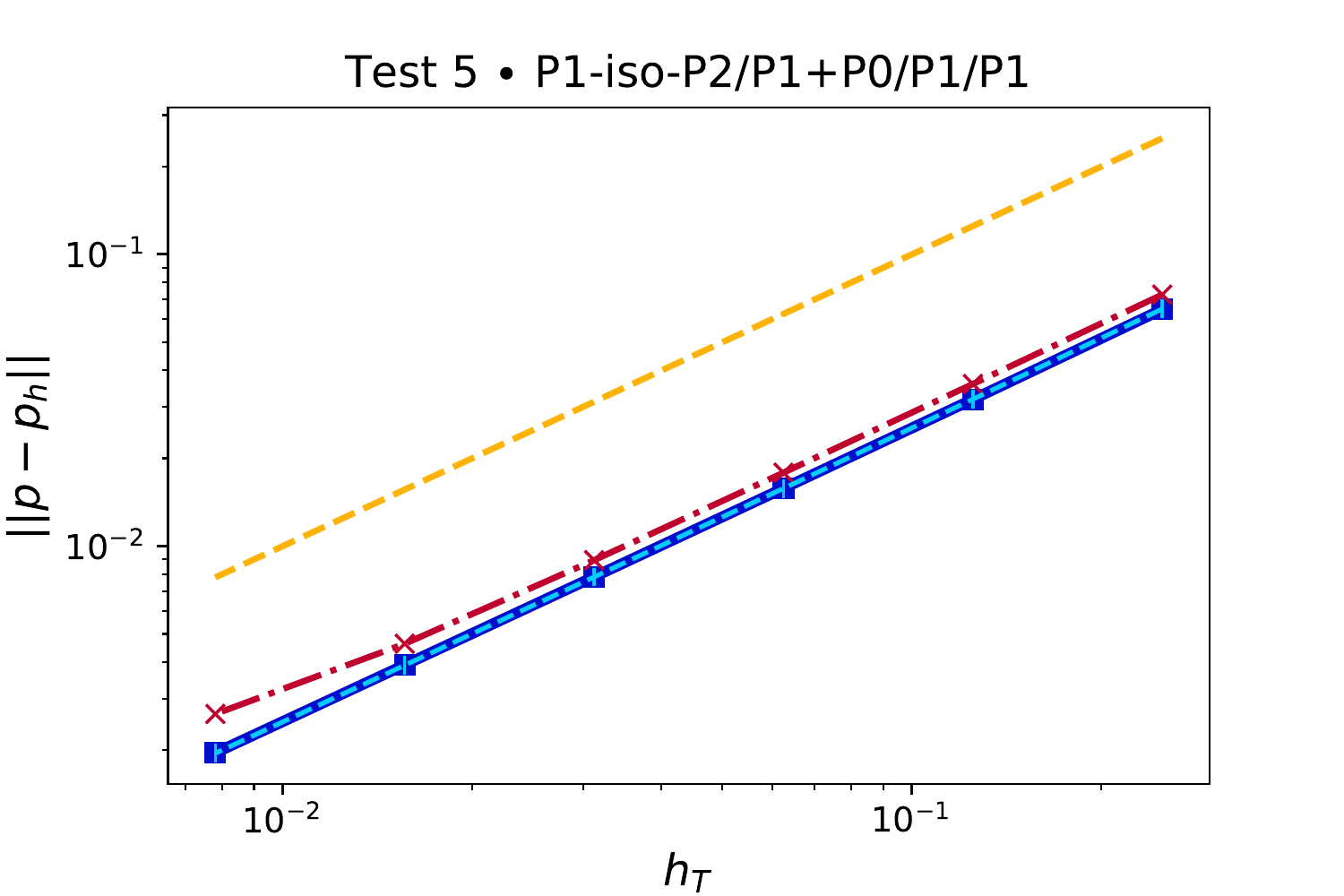}
		\includegraphics[width=0.24\linewidth]{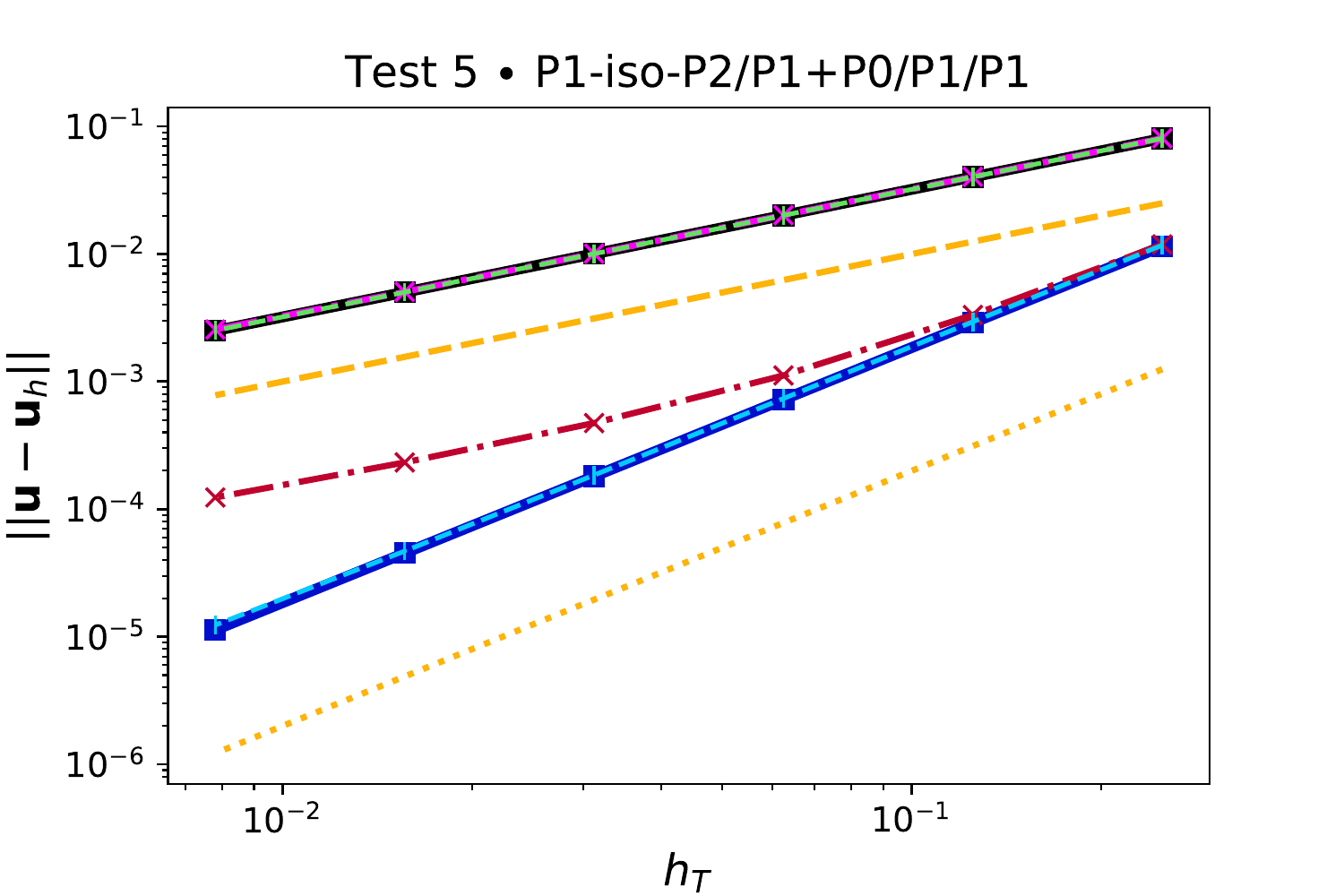}
		\includegraphics[width=0.24\linewidth]{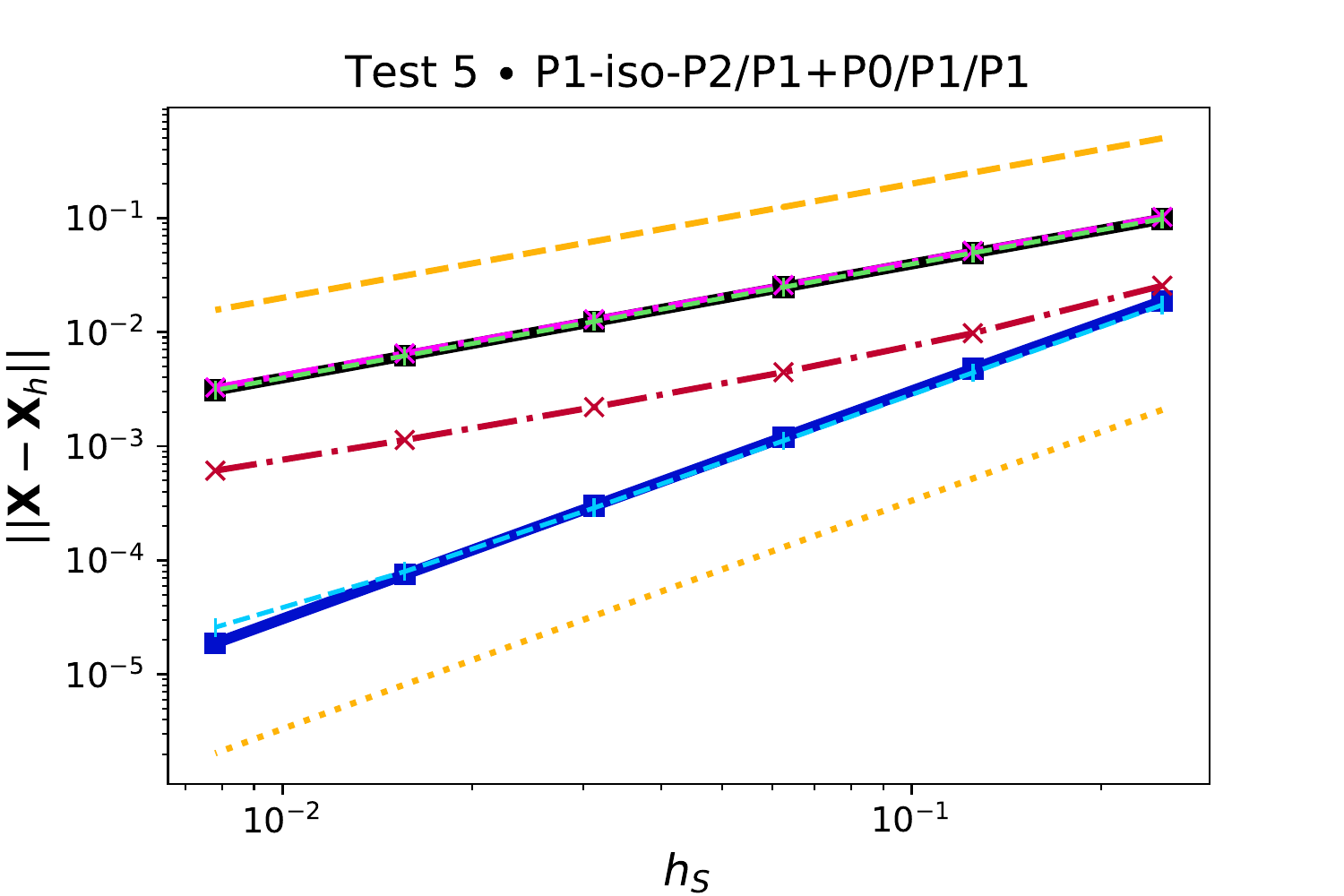}
		\includegraphics[width=0.24\linewidth]{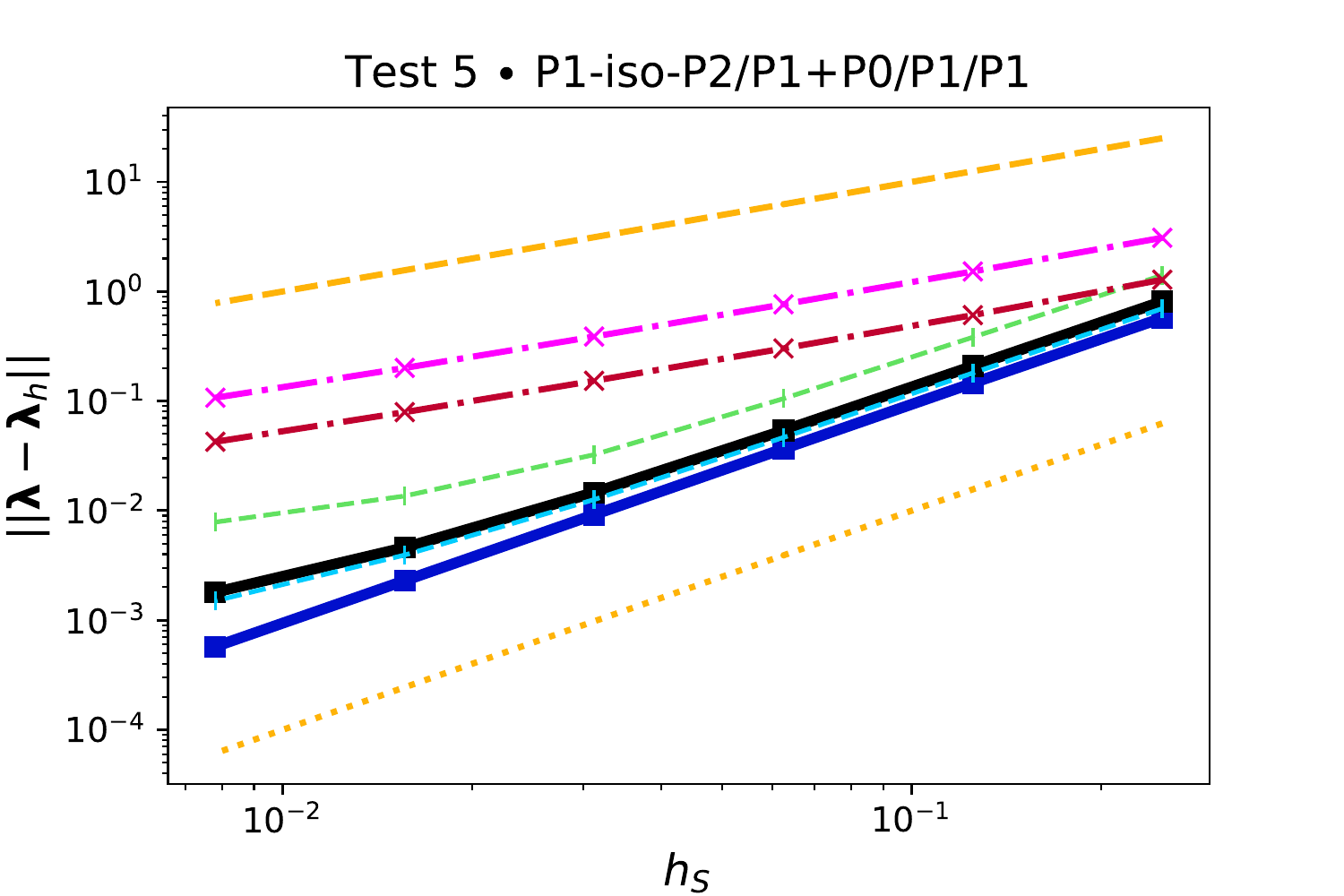}\\
		\includegraphics[trim=150 16 50 50,width=1.4\linewidth]{figures/label_p0-eps-converted-to}
		\caption{Convergence plots of Test 5 with $\Pcal_1-iso-\Pcal_2/\Pcal_1+\Pcal_0/\Pcal_1/\Pcal_1$}
		\label{fig:test5_p1p0}
	\end{figure}

	Regarding the computation of the intersection, we can see again that avoiding this procedure we do not obtain the optimal convergence rate. In particular, we can see that the behavior is similar to the one seen in the first test. The results are collected in Figures~\ref{fig:test5_p1} and \ref{fig:test5_p1p0}.
	
	\subsection{Test 6}\label{sub:test6}
\begin{figure}[!h]
		\centering
		\includegraphics[width=0.24\linewidth]{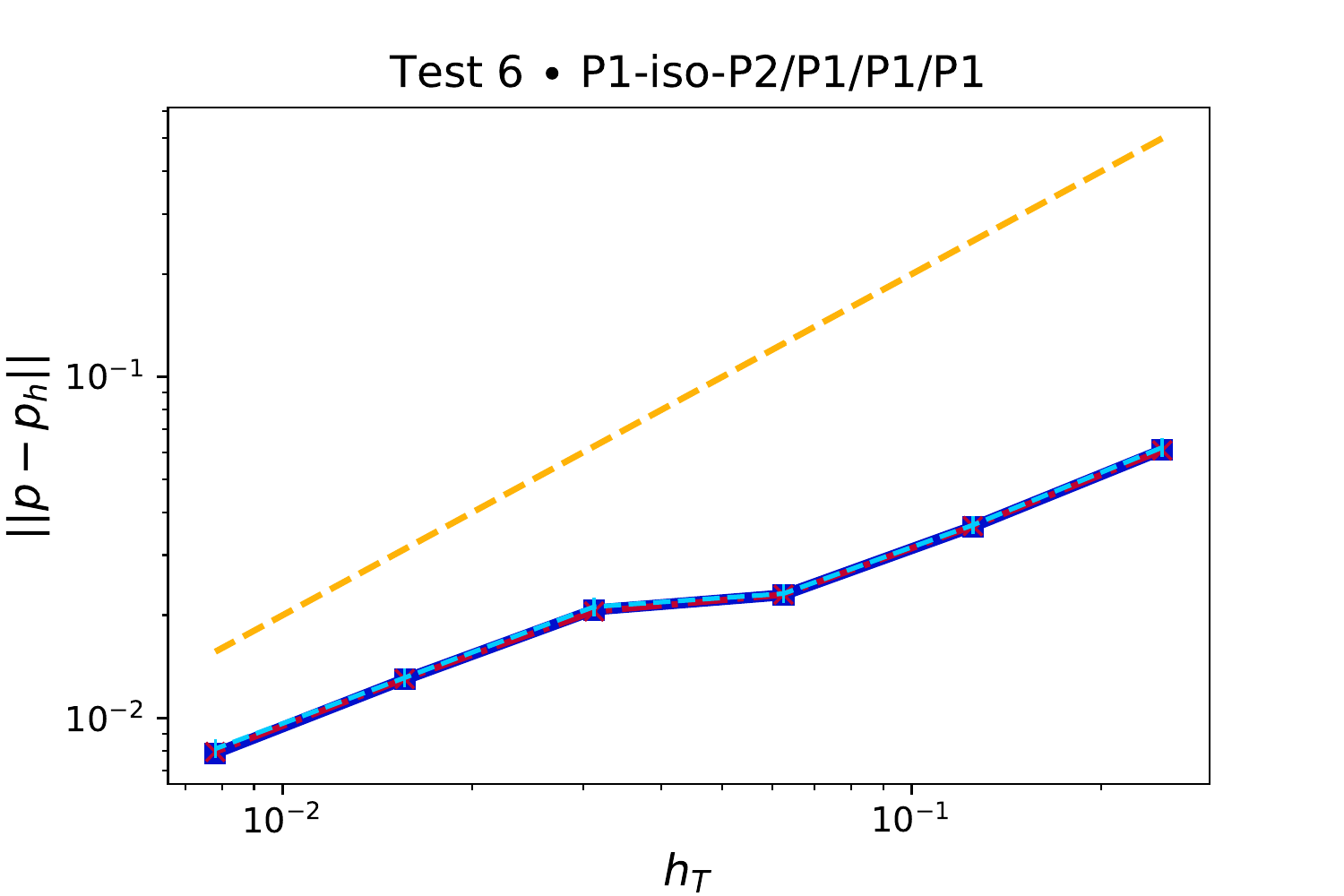}
		\includegraphics[width=0.24\linewidth]{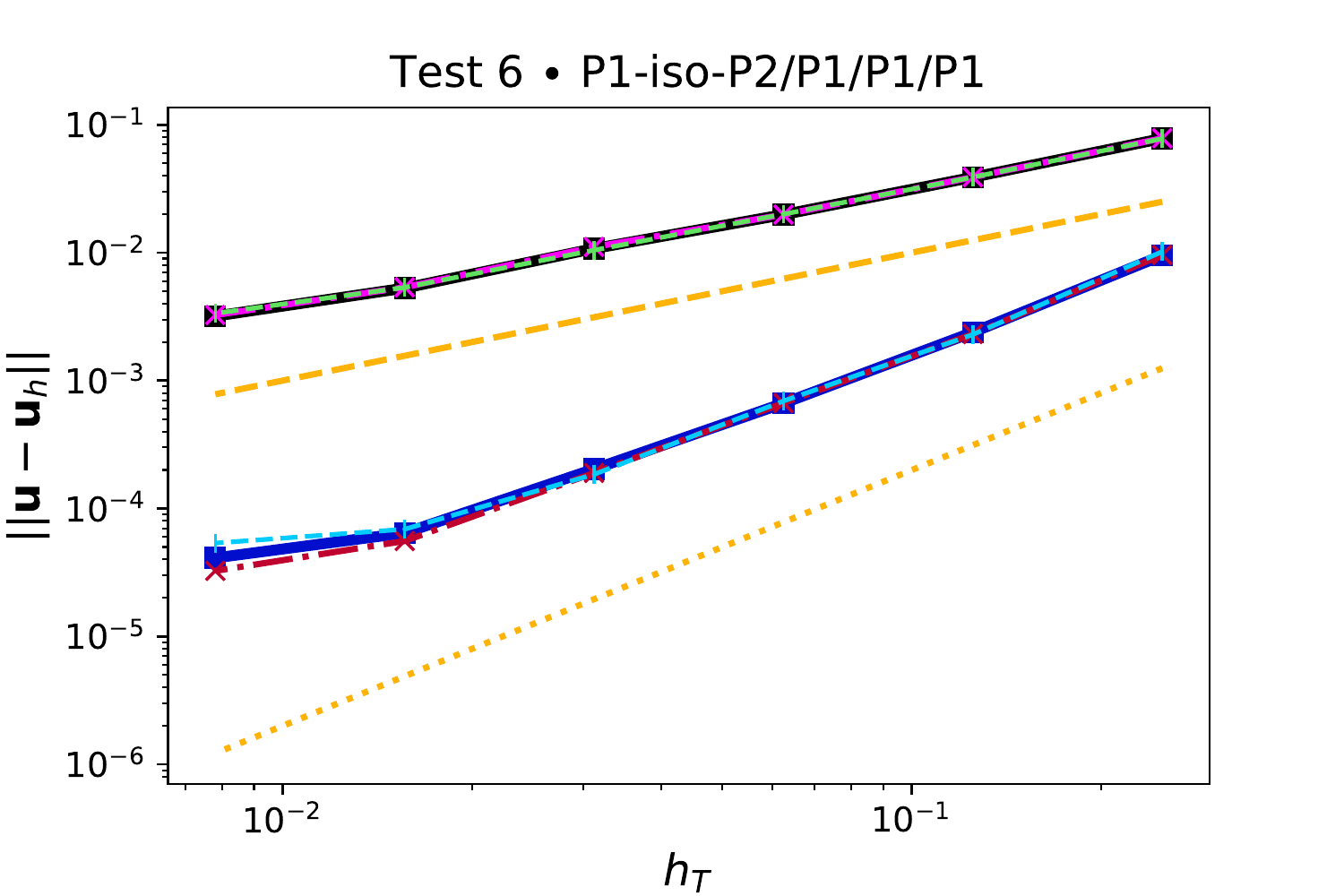}
		\includegraphics[width=0.24\linewidth]{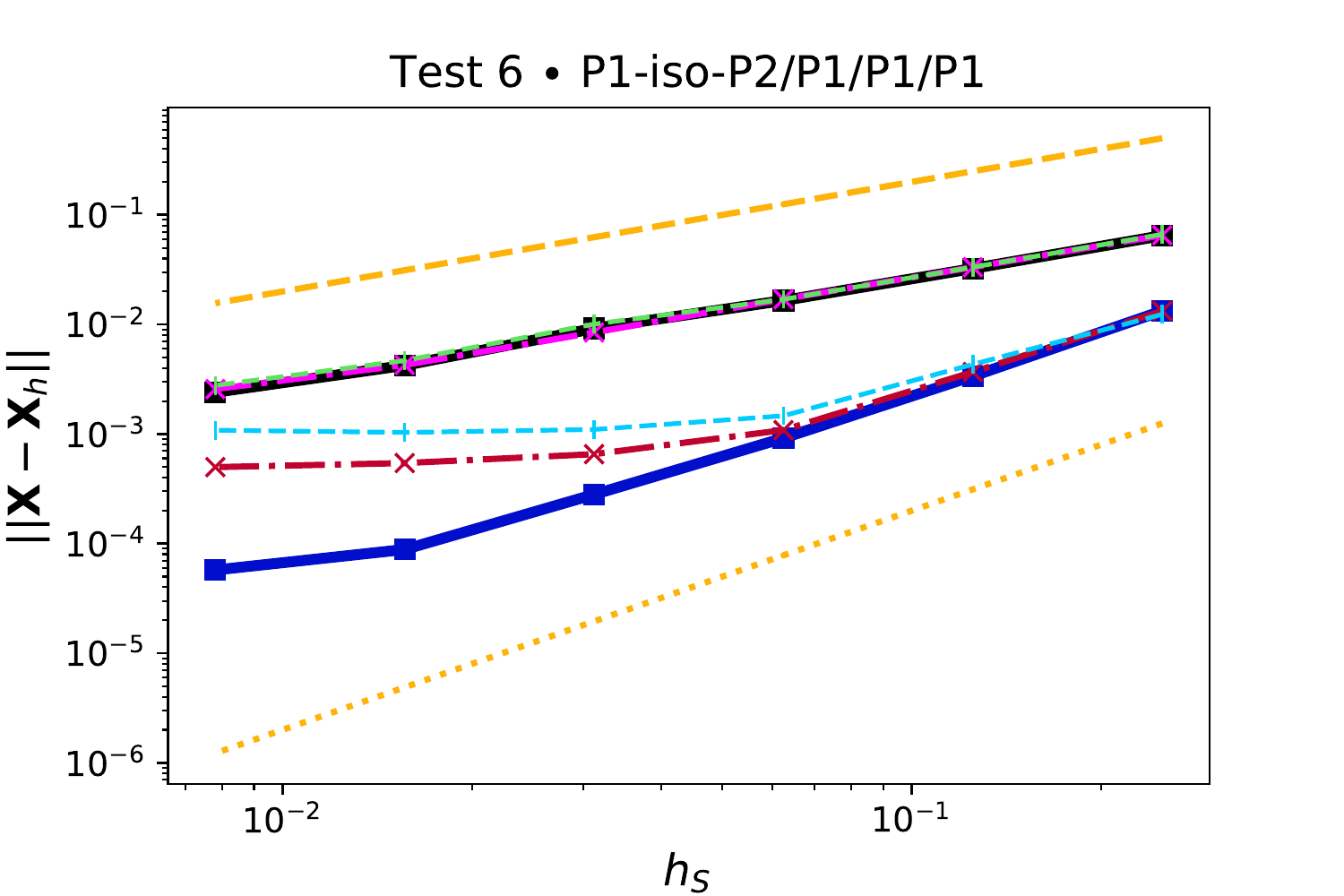}
		\includegraphics[width=0.24\linewidth]{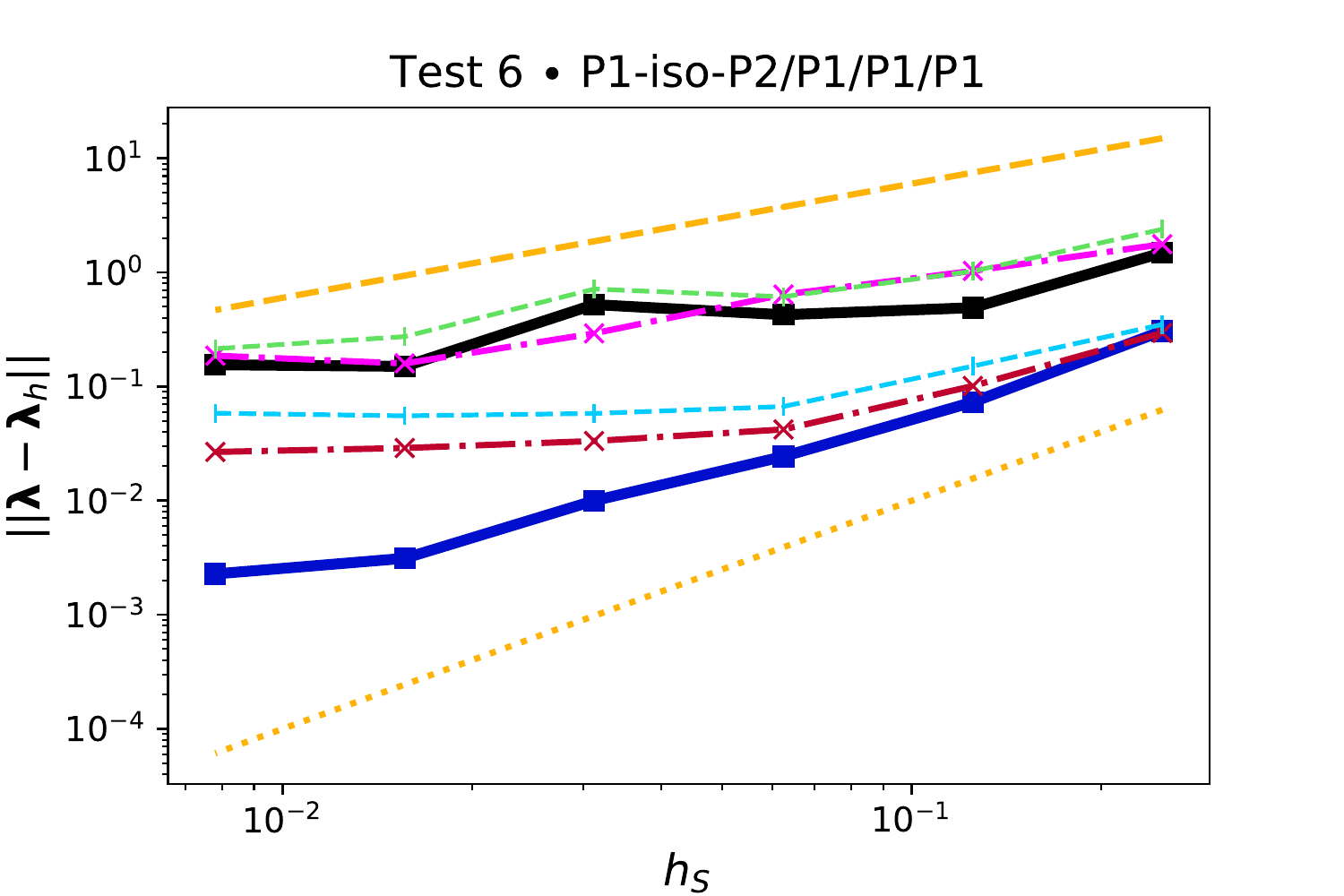}\\
		\includegraphics[trim=180 16 50 50,width=1.35\linewidth]{figures/label_p1-eps-converted-to}
		\caption{Convergence plots of Test 6 with $\Pcal_1-iso-\Pcal_2/\Pcal_1/\Pcal_1/\Pcal_1$}
		\label{fig:test6_p1}
	\end{figure}

\begin{figure}[!h]
		\centering
		\includegraphics[width=0.24\linewidth]{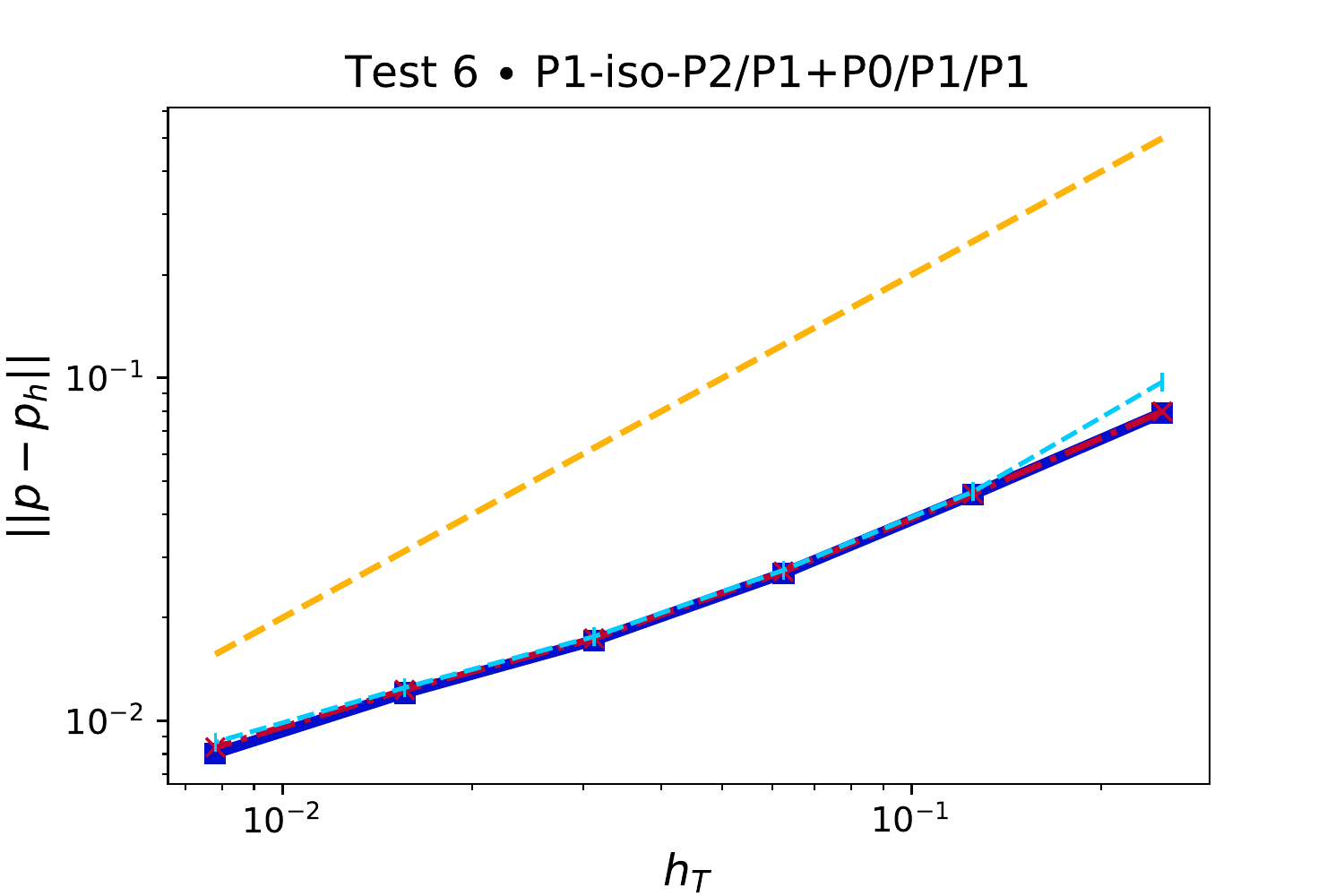}
		\includegraphics[width=0.24\linewidth]{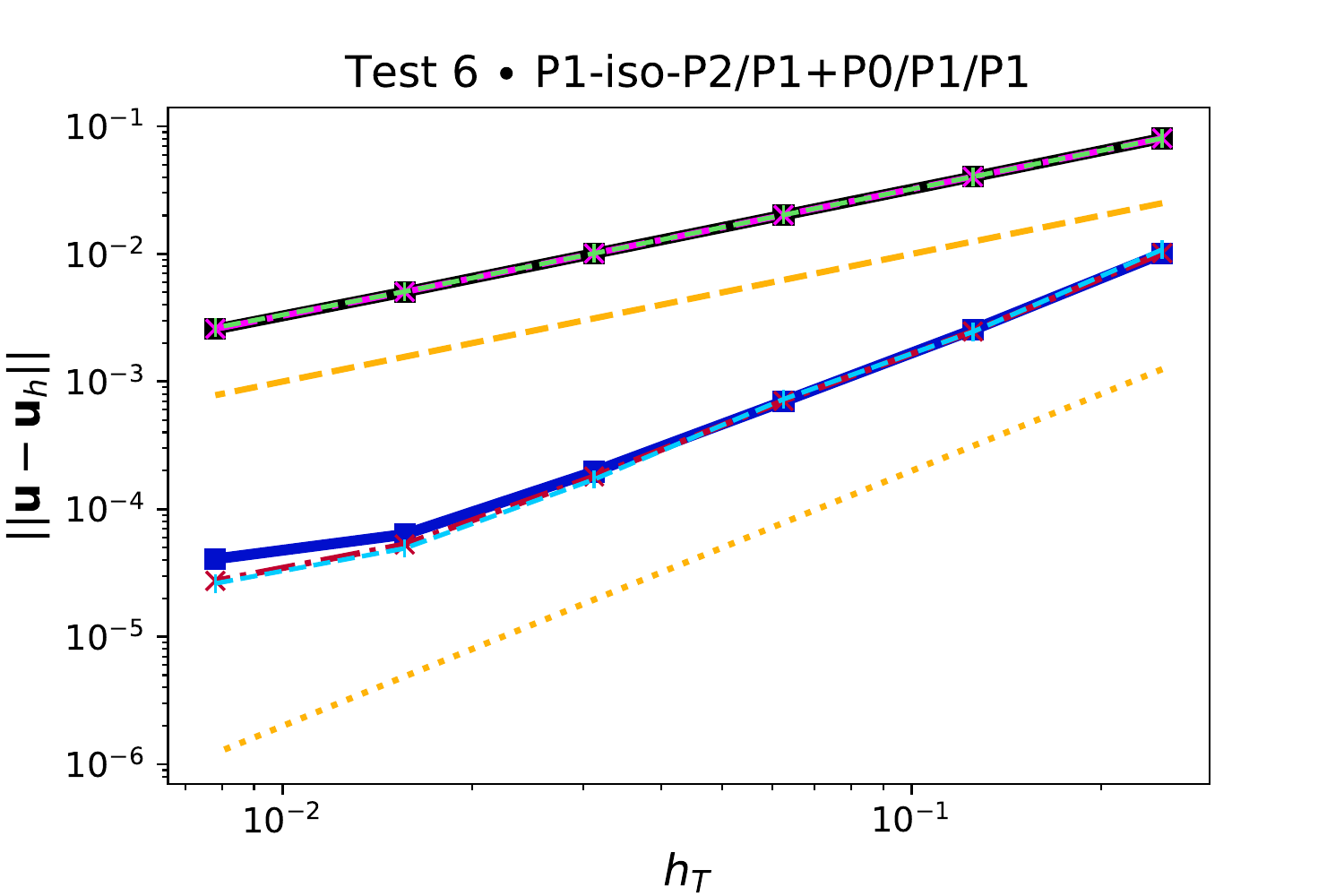}
		\includegraphics[width=0.24\linewidth]{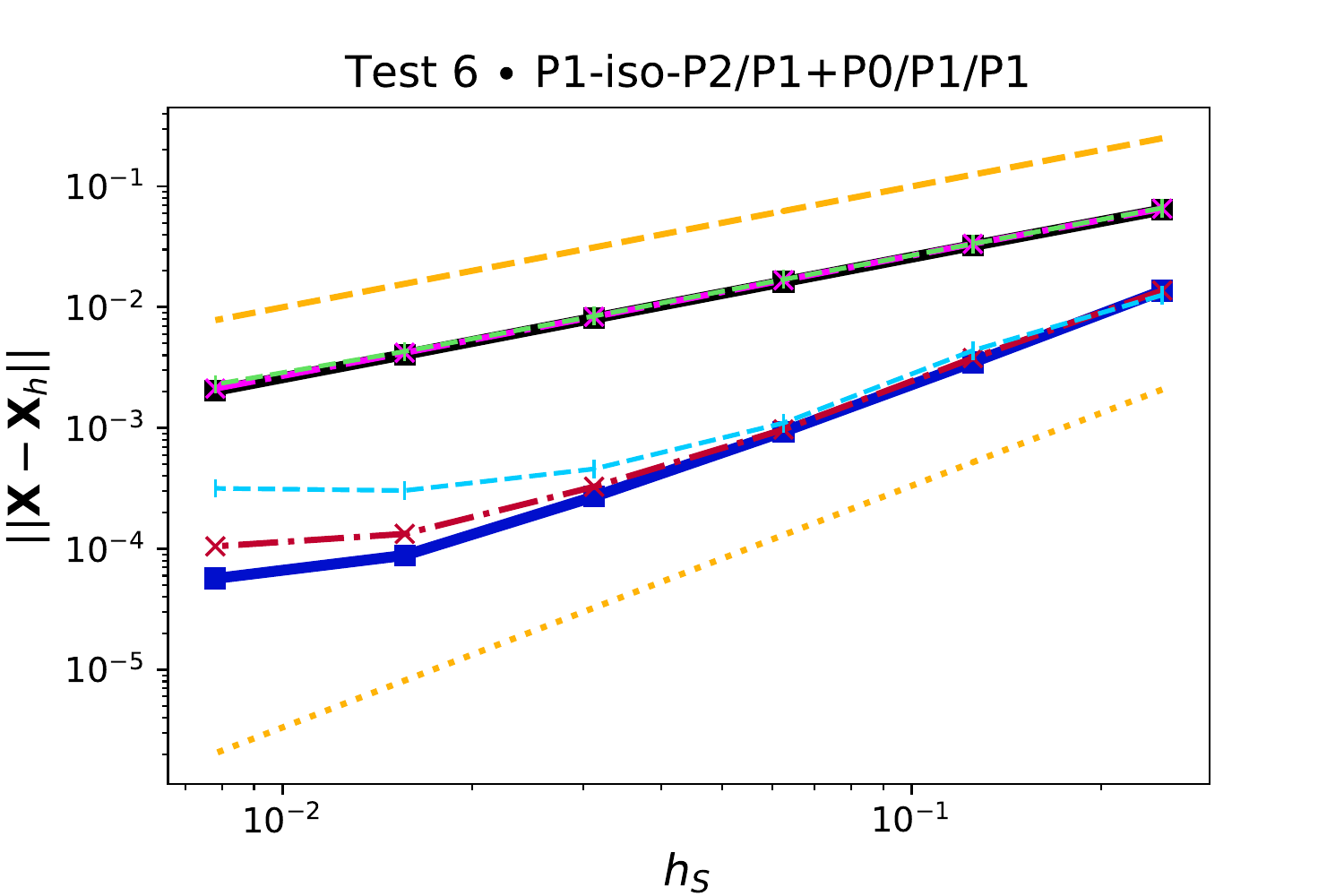}
		\includegraphics[width=0.24\linewidth]{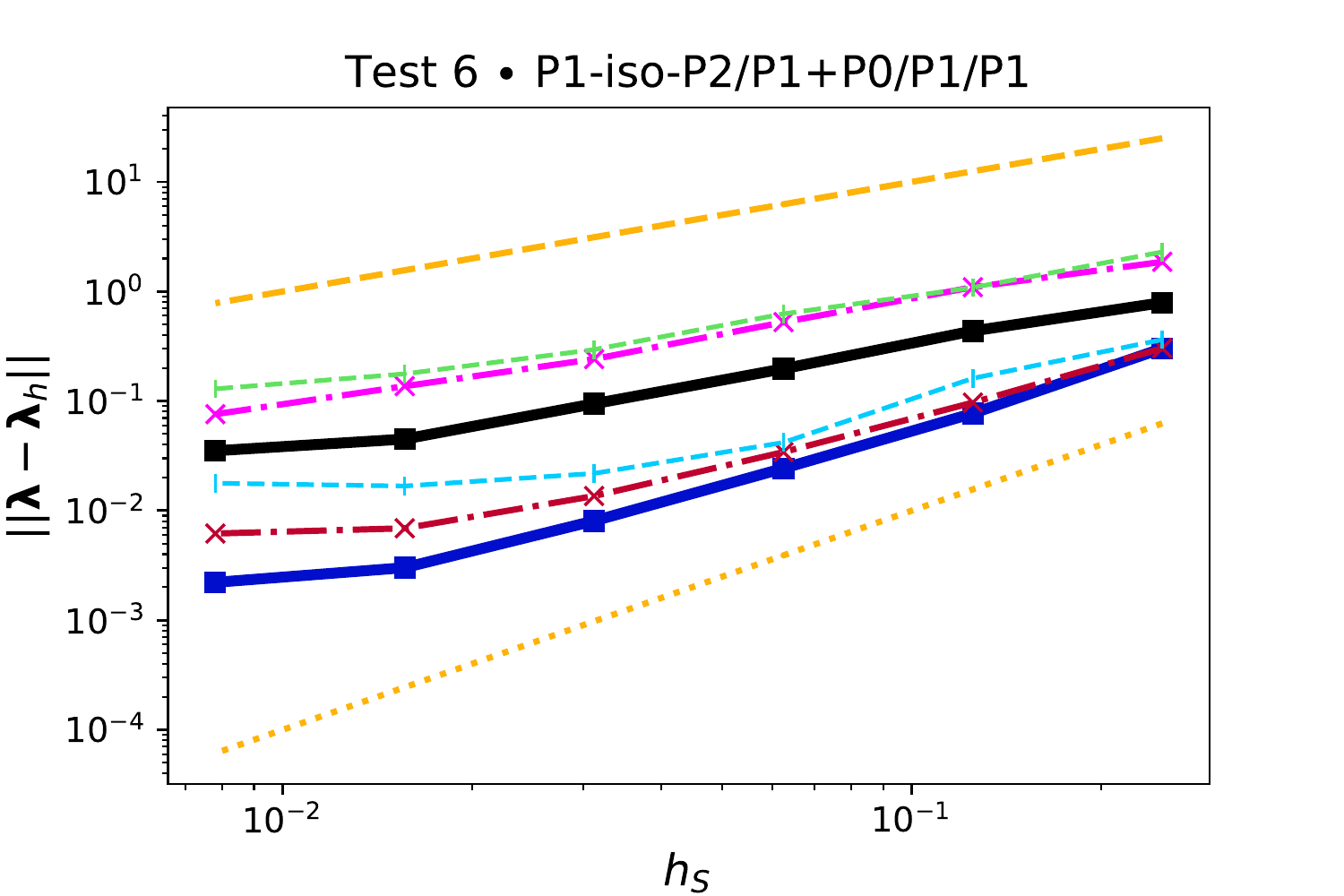}\\
		\includegraphics[trim=150 16 50 50,width=1.4\linewidth]{figures/label_p0-eps-converted-to}
		\caption{Convergence plots of Test 6 with $\Pcal_1-iso-\Pcal_2/\Pcal_1+\Pcal_0/\Pcal_1/\Pcal_1$}
		\label{fig:test6_p1p0}
	\end{figure}

	We consider again a discontinuous pressure \blue on the boundary of the structure but in the case where the fluid triangulation does not match it\noblue: for this reason, we choose ${\B=\Os=[-\pi/4, \pi/4]^2}$ and we discretize it with a uniform left mesh. Our choices for $\Omega$ and its discretization remain unchanged.
	We compute the right hand sides $\f$, $\g$, and $\d$ in order to obtain the following solutions:
	\begin{equation}
		\begin{aligned}
			&\u(x,y)=\curl\big((4-x^2)^2(4-y^2)^2\big)\\
			&\u_{\mid\partial\Omega}=\mathbf{0}\\
			&p(x,y) = 
			\begin{cases}
				150\sin(x)+50\frac{\pi^2}{4}\big(\frac{\pi^2}{4}-16\big)^{-1}&\text{ in }\Of\\
				150\sin(x)+50&\text{ in }\Os=\B.
			\end{cases}\\
			&\X(x,y)=\u(x,y)\\
			&\llambda(x,y)=\big(e^x,e^y\big).
		\end{aligned}
	\end{equation}
	Due to the non matching boundary with respect to the singularity, also the $\Pcal_1+\Pcal_0$ choice does not reach the optimal convergence rate presenting the Gibbs phenomenon. In this case, as we can see in Figures~\ref{fig:test6_p1} and \ref{fig:test6_p1p0}, the non-computation of the mesh intersection does not affect the fluid variables basically; on the other hand, for the solid ones, the behavior changes drastically.
	
	\subsection{Test 7, disk}\label{sub:test7}
	At this point, we consider a situation where the initial map $\Xbar$ is not trivial: indeed, we consider as reference domain for the solid the square $\B=[-1,1]^2$ (partitioned with a uniform left mesh), while the actual immersed body is the unit disk ${\Os=\{ \x\in\R^2:\norm{x}\le 1\}}$; hence we have
	\begin{equation*}
		\Xbar(x,y) = \bigg( x\sqrt{1-\frac{y^2}{2}} ,\,y\sqrt{1-\frac{x^2}{2}} \bigg).
	\end{equation*}
	On the other hand, the domain $\Omega$ is still the same as before, and the computations are done in order to obtain as solutions the ones chosen in \eqref{eq:standard_sol}. In addition, since the map $\Xbar$ is nontrivial, $\d=\u(\Xbar)-\X\neq\mathbf{0}$ and we have again to deal with the computation of $\c(\mmu,\u(\Xbar))$ \red intersecting \nored the two meshes as we do in order to compute $\f$.
	
\begin{figure}[!h]
		\centering
		\includegraphics[width=0.24\linewidth]{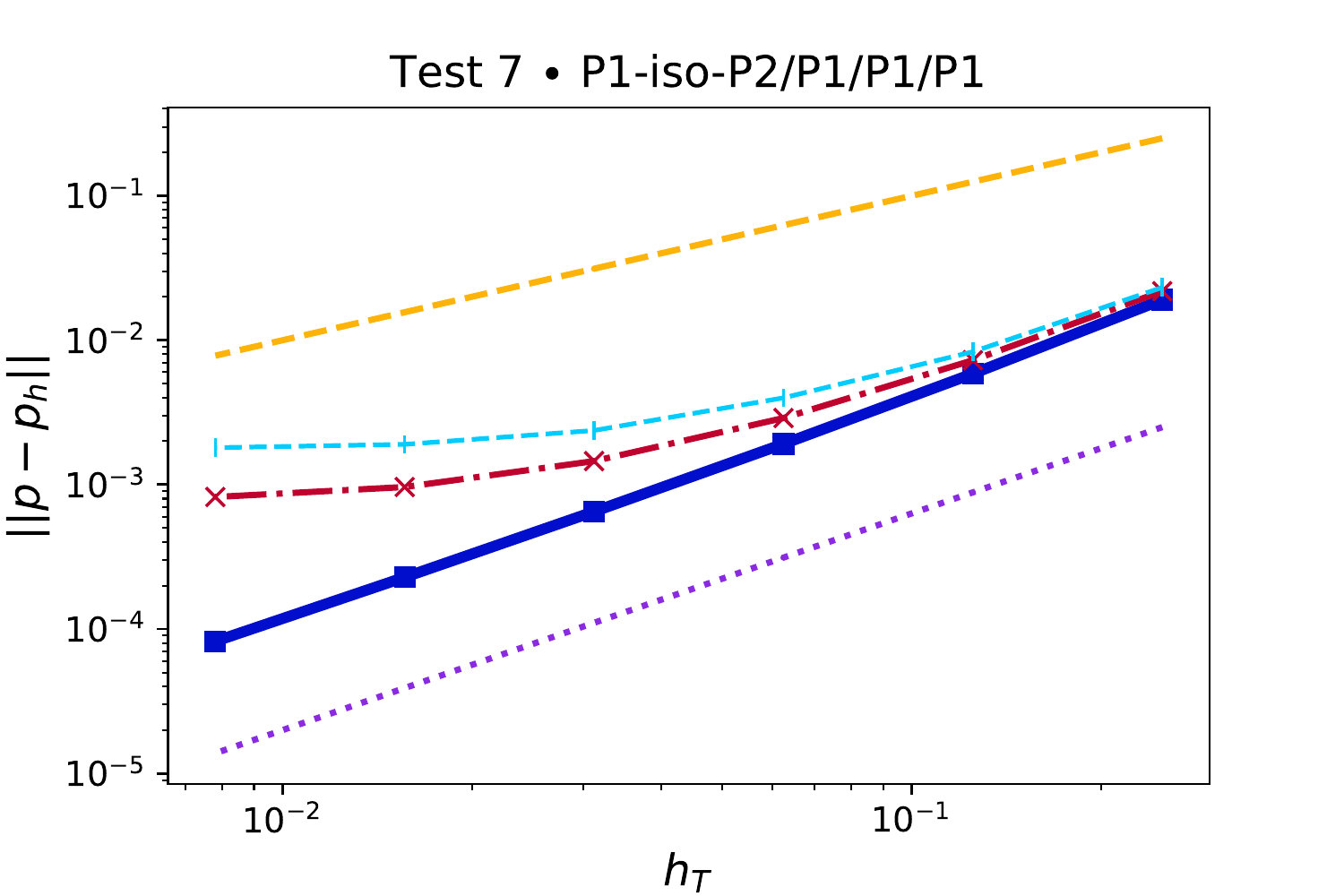}
		\includegraphics[width=0.24\linewidth]{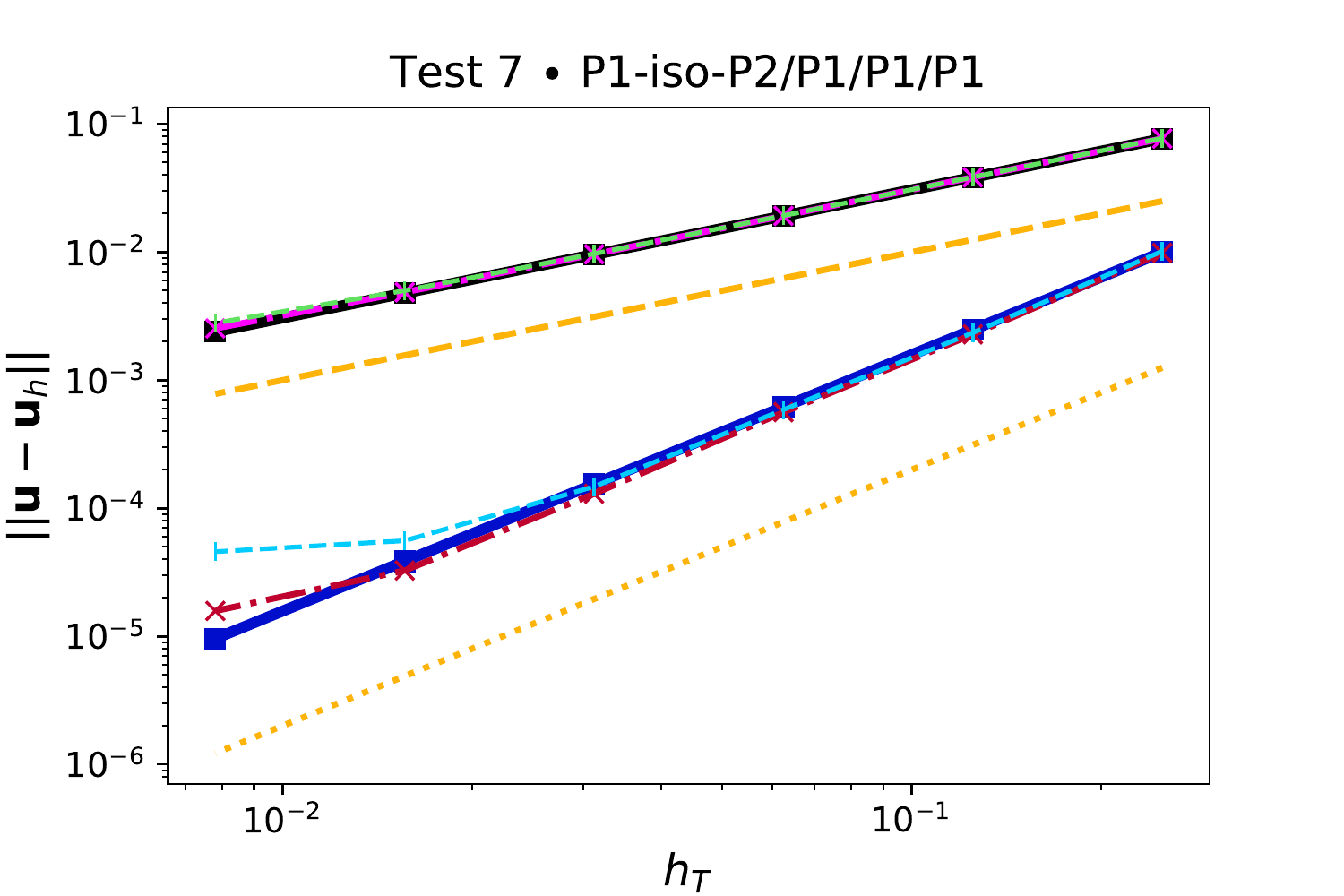}
		\includegraphics[width=0.24\linewidth]{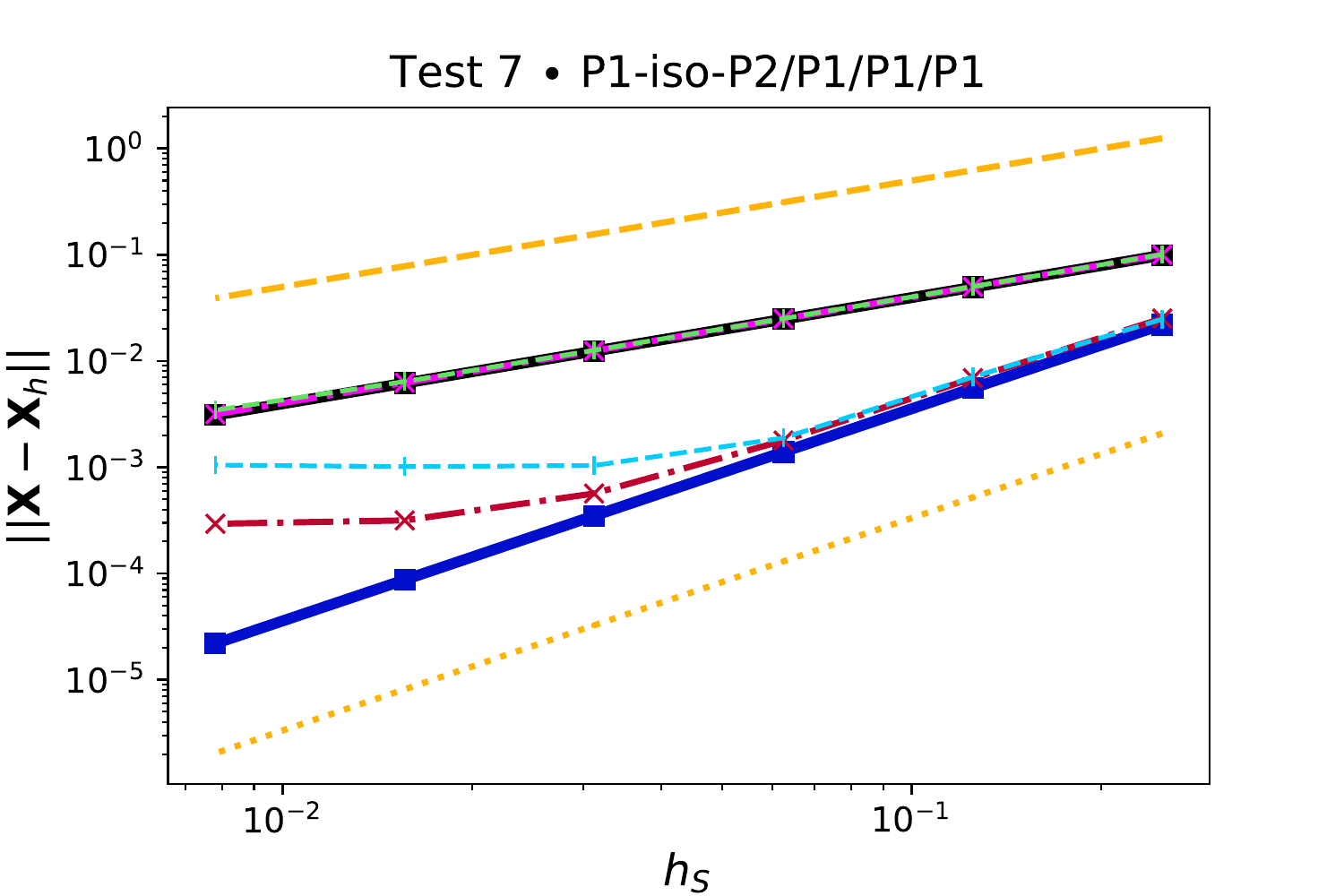}
		\includegraphics[width=0.24\linewidth]{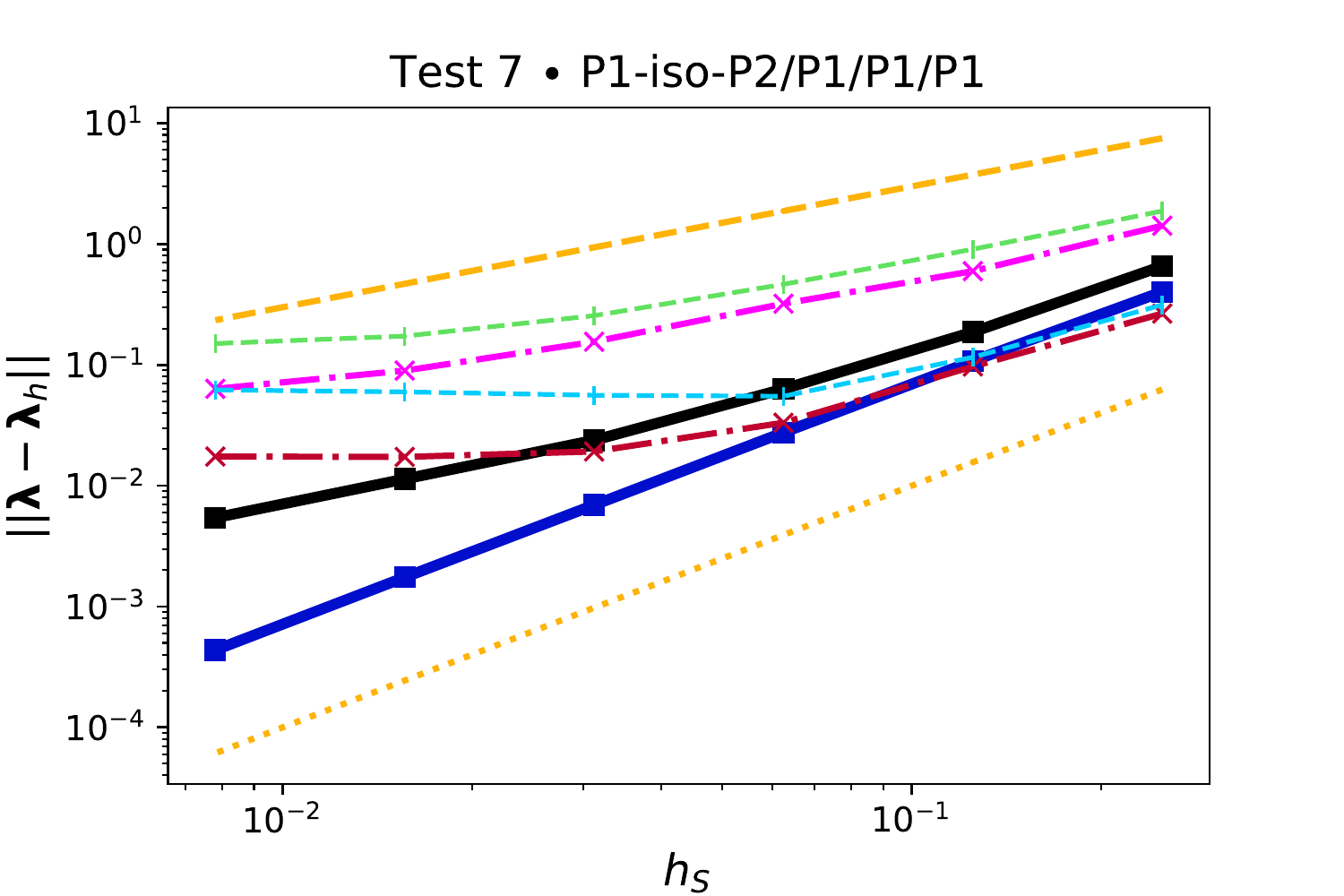}\\
		\includegraphics[trim=180 16 50 50,width=1.35\linewidth]{figures/label_p1-eps-converted-to}
		\caption{Convergence plots of Test 7 with $\Pcal_1-iso-\Pcal_2/\Pcal_1/\Pcal_1/\Pcal_1$}
		\label{fig:test7_p1}
	\end{figure}

\begin{figure}[!h]
		\centering
		\includegraphics[width=0.24\linewidth]{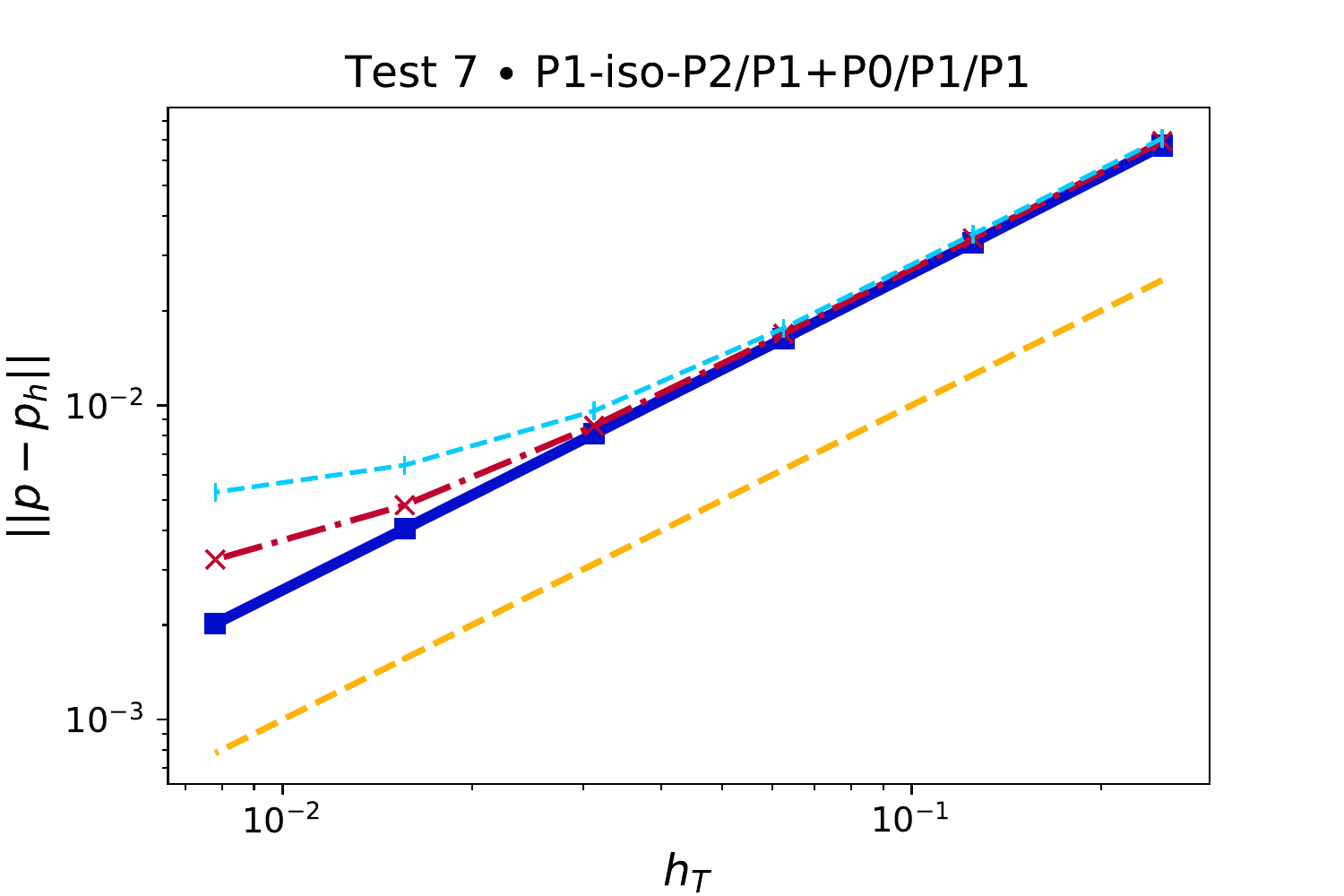}
		\includegraphics[width=0.24\linewidth]{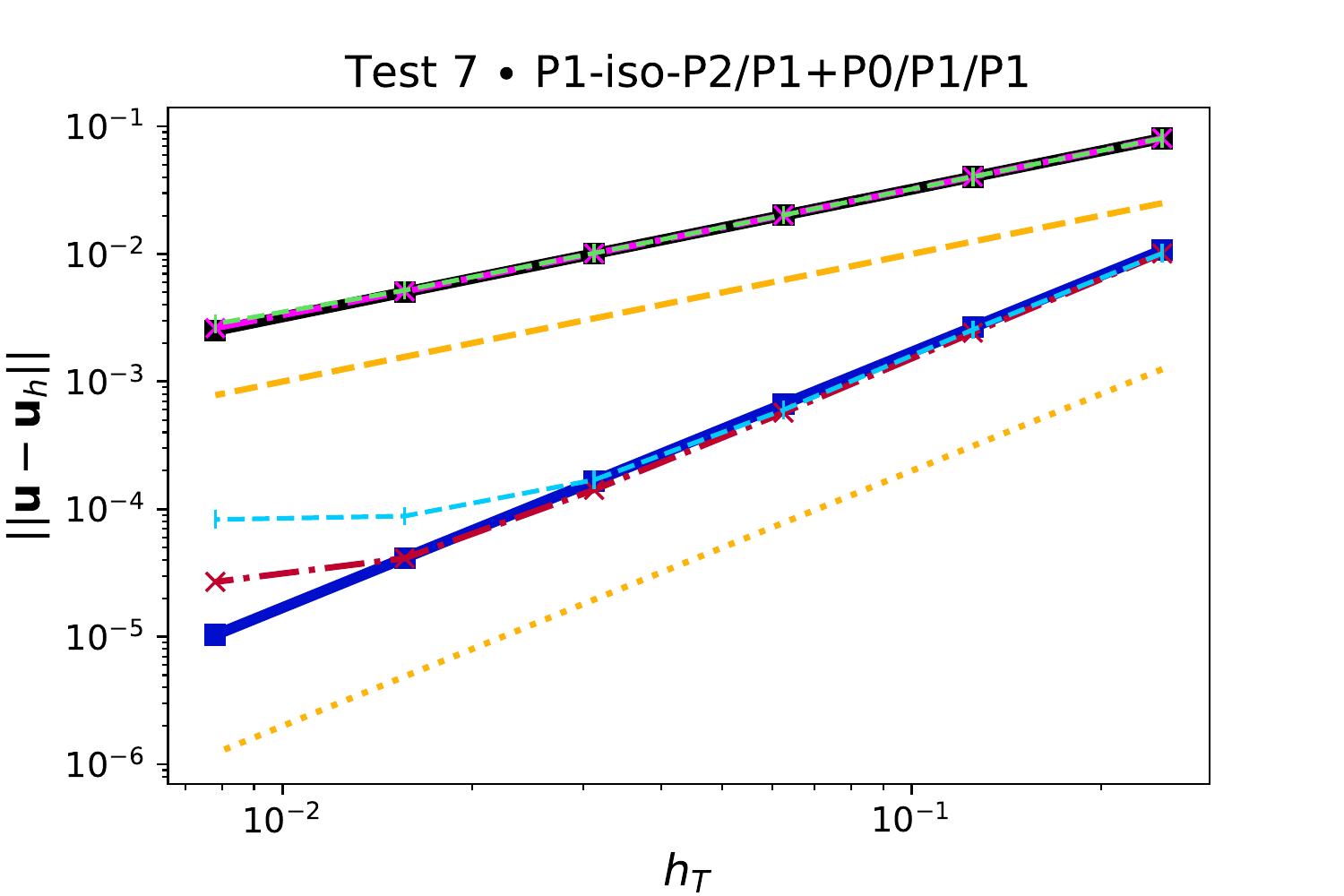}
		\includegraphics[width=0.24\linewidth]{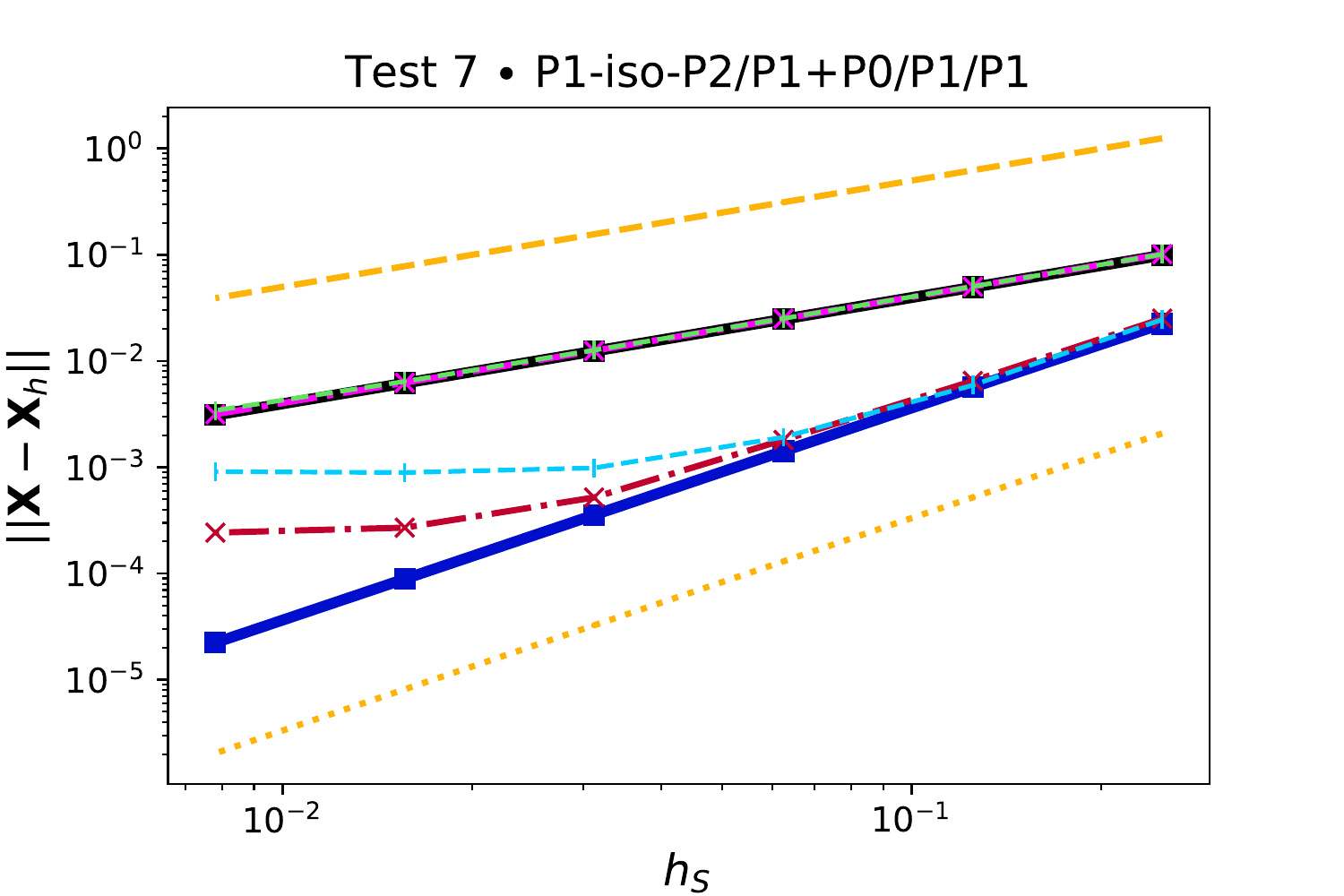}
		\includegraphics[width=0.24\linewidth]{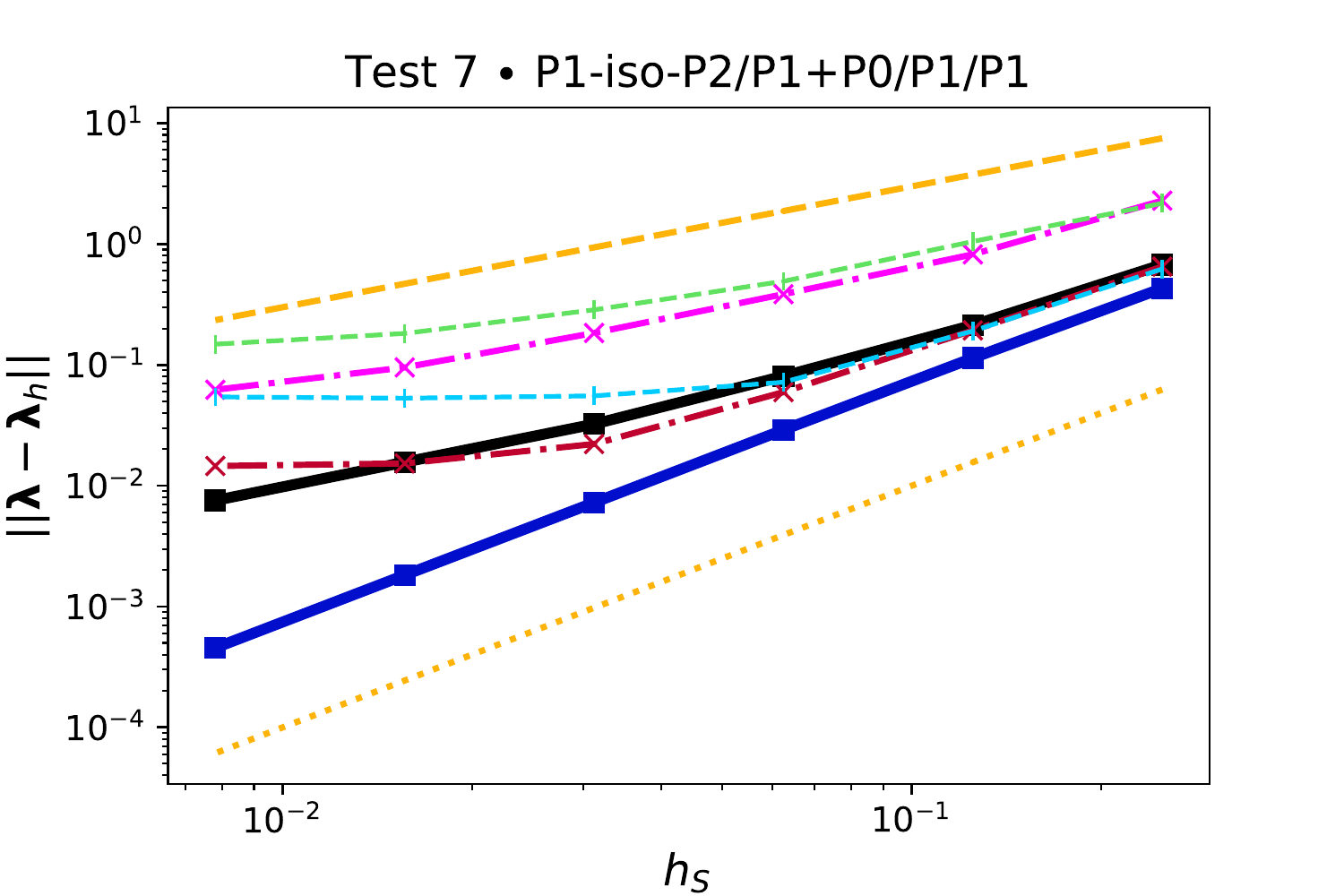}\\
		\includegraphics[trim=150 16 50 50,width=1.4\linewidth]{figures/label_p0-eps-converted-to}
		\caption{Convergence plots of Test 7 with $\Pcal_1-iso-\Pcal_2/\Pcal_1+\Pcal_0/\Pcal_1/\Pcal_1$}
		\label{fig:test7_p1p0}
	\end{figure}

	The substantial difference between this test and the previous ones is that the evaluation of the solid basis functions is done after mapping back the nodes on the reference domain $\B$. Also in this case, the behavior observed previously is confirmed with the intersection approach prevailing on the non-intersection one (see Figures \ref{fig:test7_p1} and \ref{fig:test7_p1p0}).
	
	\subsection{Test 8, square}\label{sub:test8}
	This test is a variation of Test 2 and 3: we consider the same geometry but in the case of a non trivial mapping. Indeed, the solid reference domain is chosen to be the unit square $\B=[0,1]^2$, while the actual body is an another square, precisely $\Os=[-0.62,1.38]^2$. We use again uniform meshes (right and left oriented respectively) for both $\Omega=[-2,2]^2$ and $\B$. \lg The initial mapping $\Xbar$ is defined by
	\begin{equation*}
		\Xbar(x,y) = (-0.62+2x,\,-0.62+2y),
	\end{equation*}
	and the solutions are those presented in \eqref{eq:standard_sol}. \gl
	
\begin{figure}[!h]
		\centering
		\includegraphics[width=0.24\linewidth]{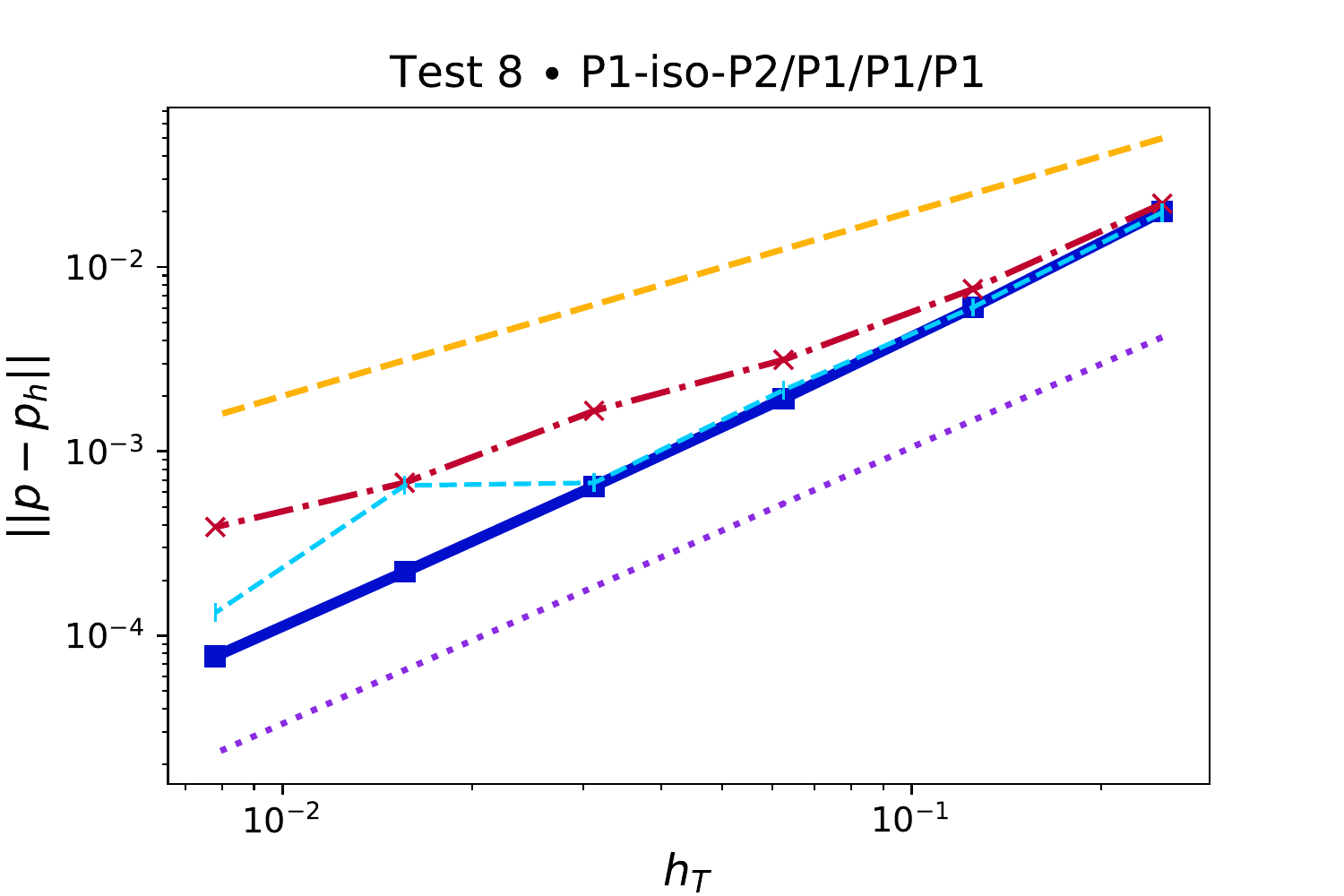}
		\includegraphics[width=0.24\linewidth]{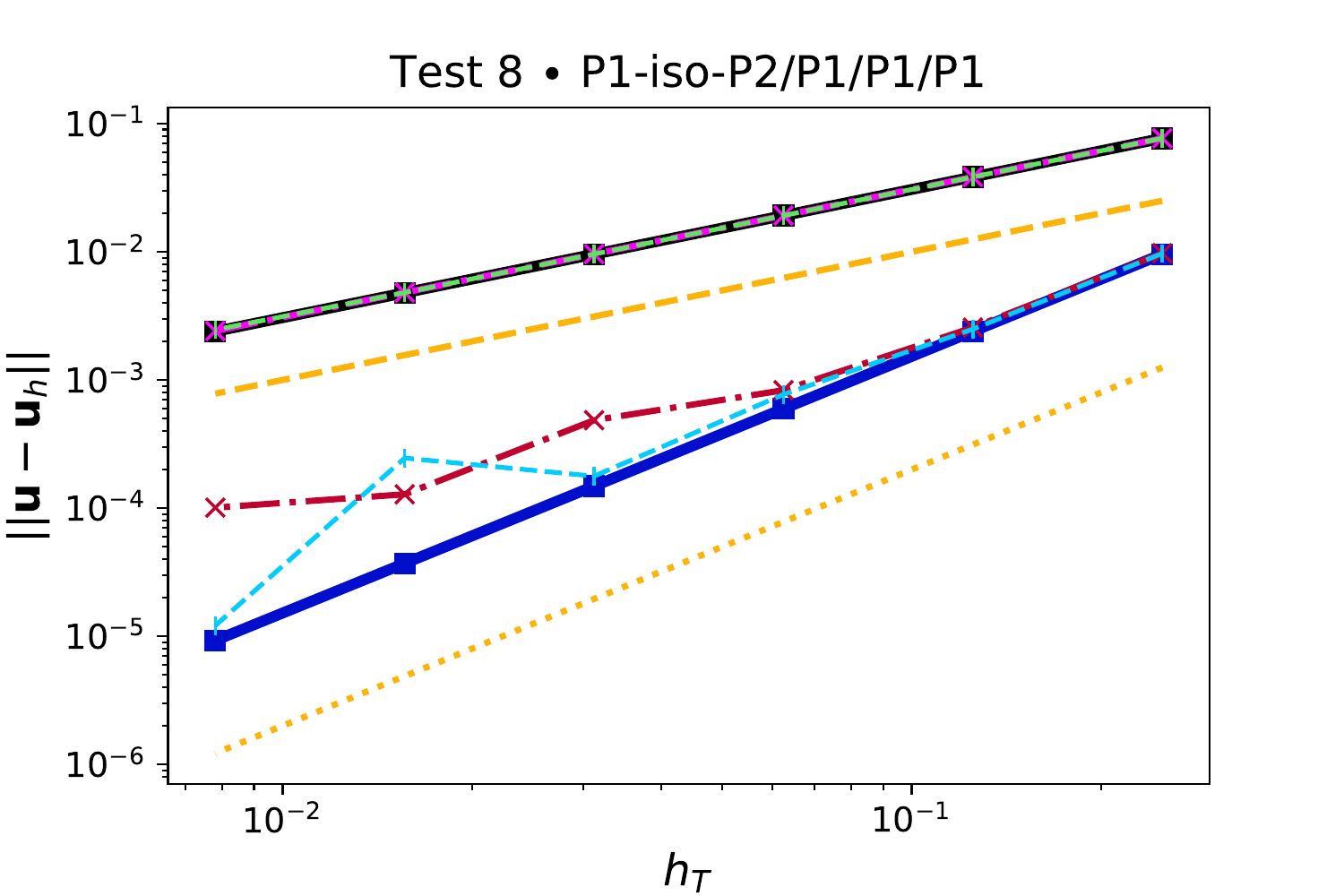}
		\includegraphics[width=0.24\linewidth]{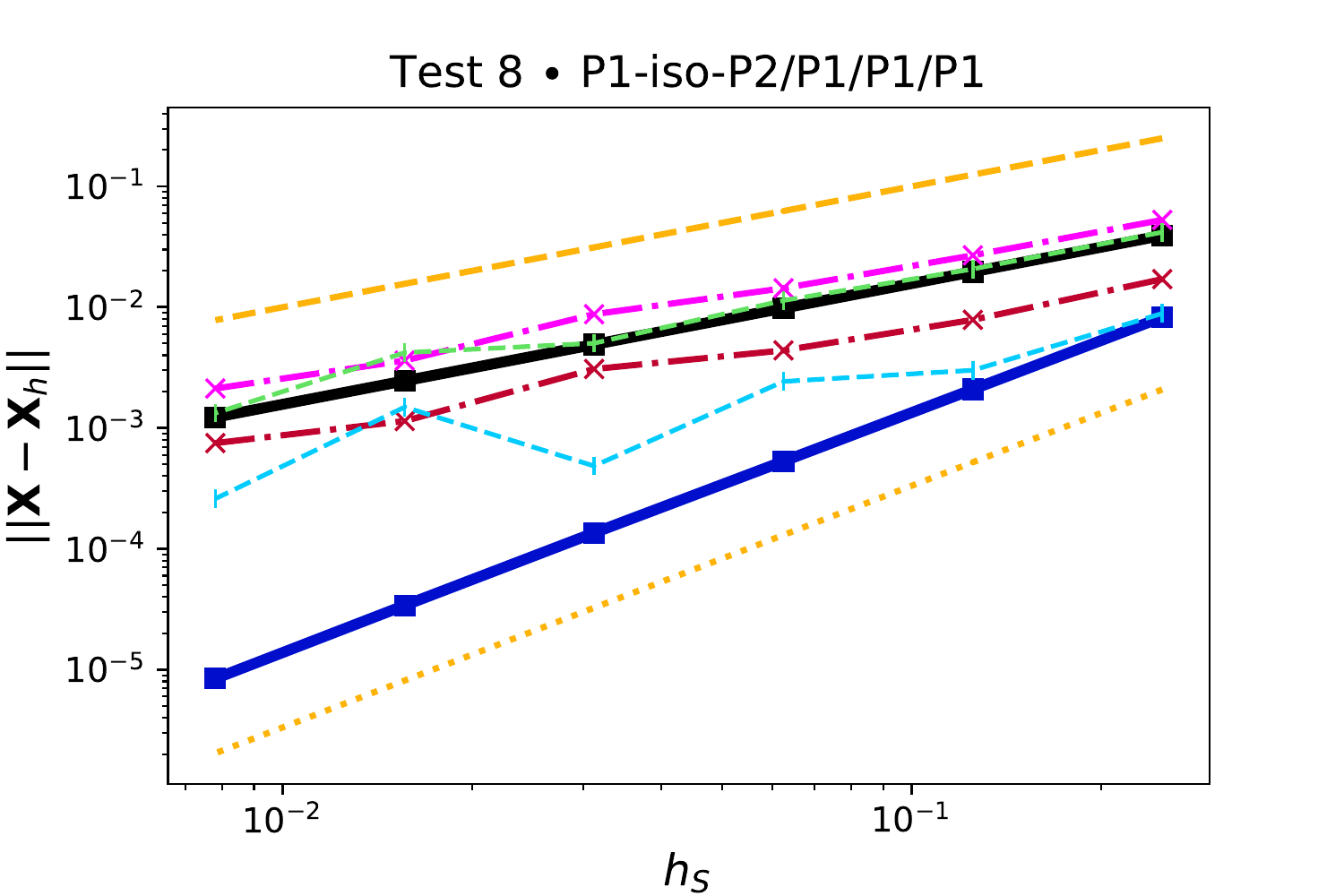}
		\includegraphics[width=0.24\linewidth]{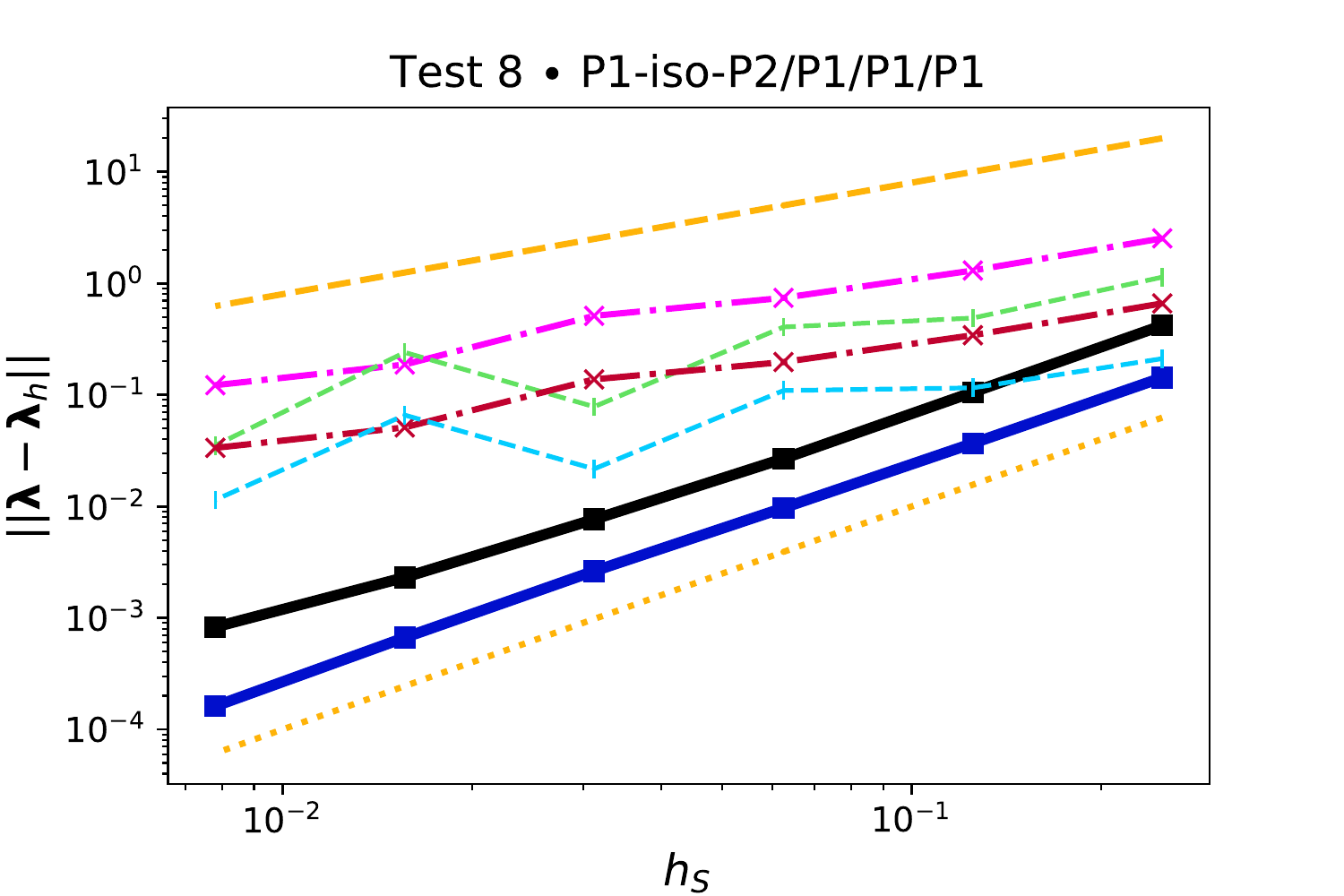}\\
		\includegraphics[trim=180 16 50 50,width=1.35\linewidth]{figures/label_p1-eps-converted-to}
		\caption{Convergence plots of Test 8 with $\Pcal_1-iso-\Pcal_2/\Pcal_1/\Pcal_1/\Pcal_1$}
		\label{fig:test8_p1}
	\end{figure}

\begin{figure}[!h]
		\centering
		\includegraphics[width=0.24\linewidth]{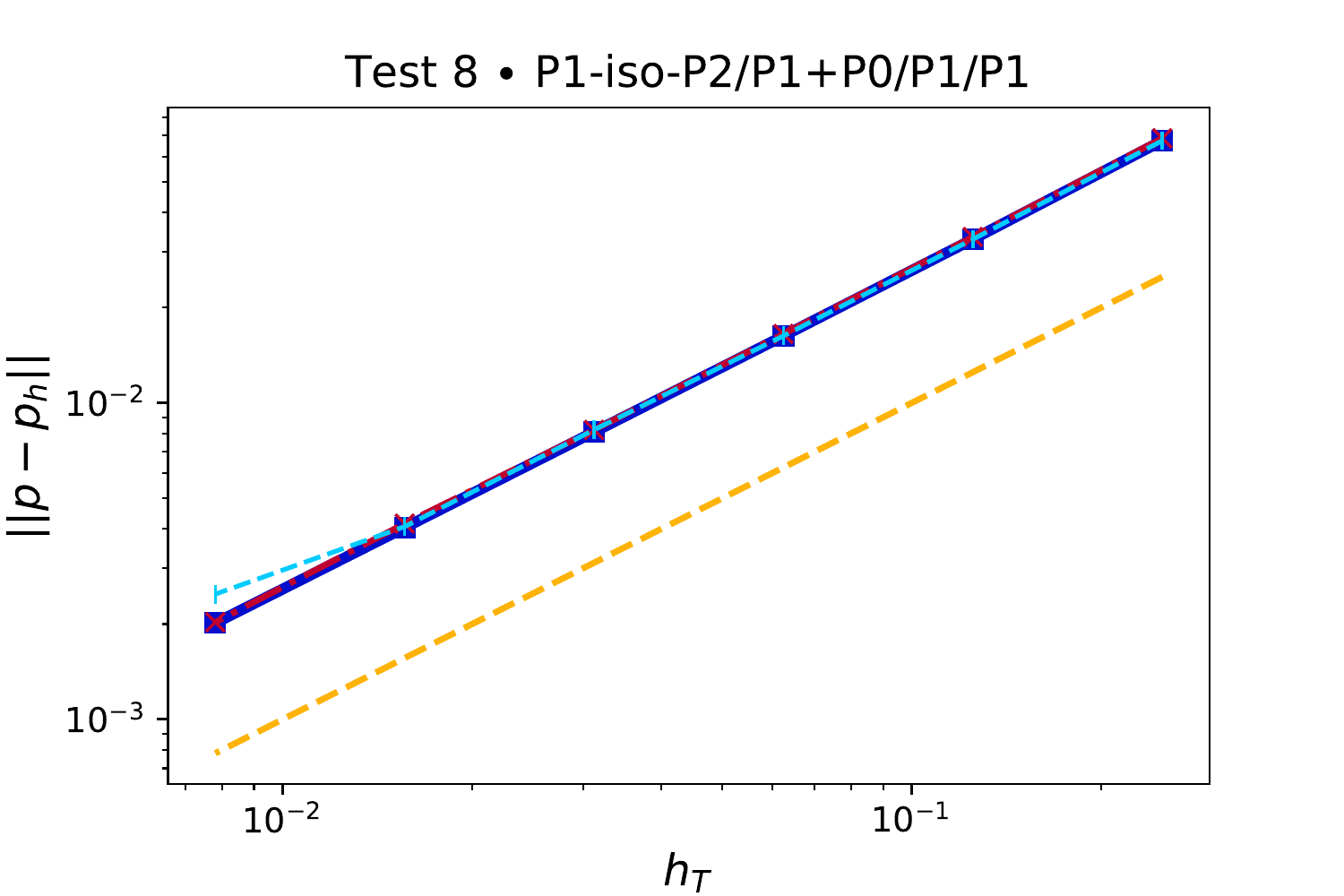}
		\includegraphics[width=0.24\linewidth]{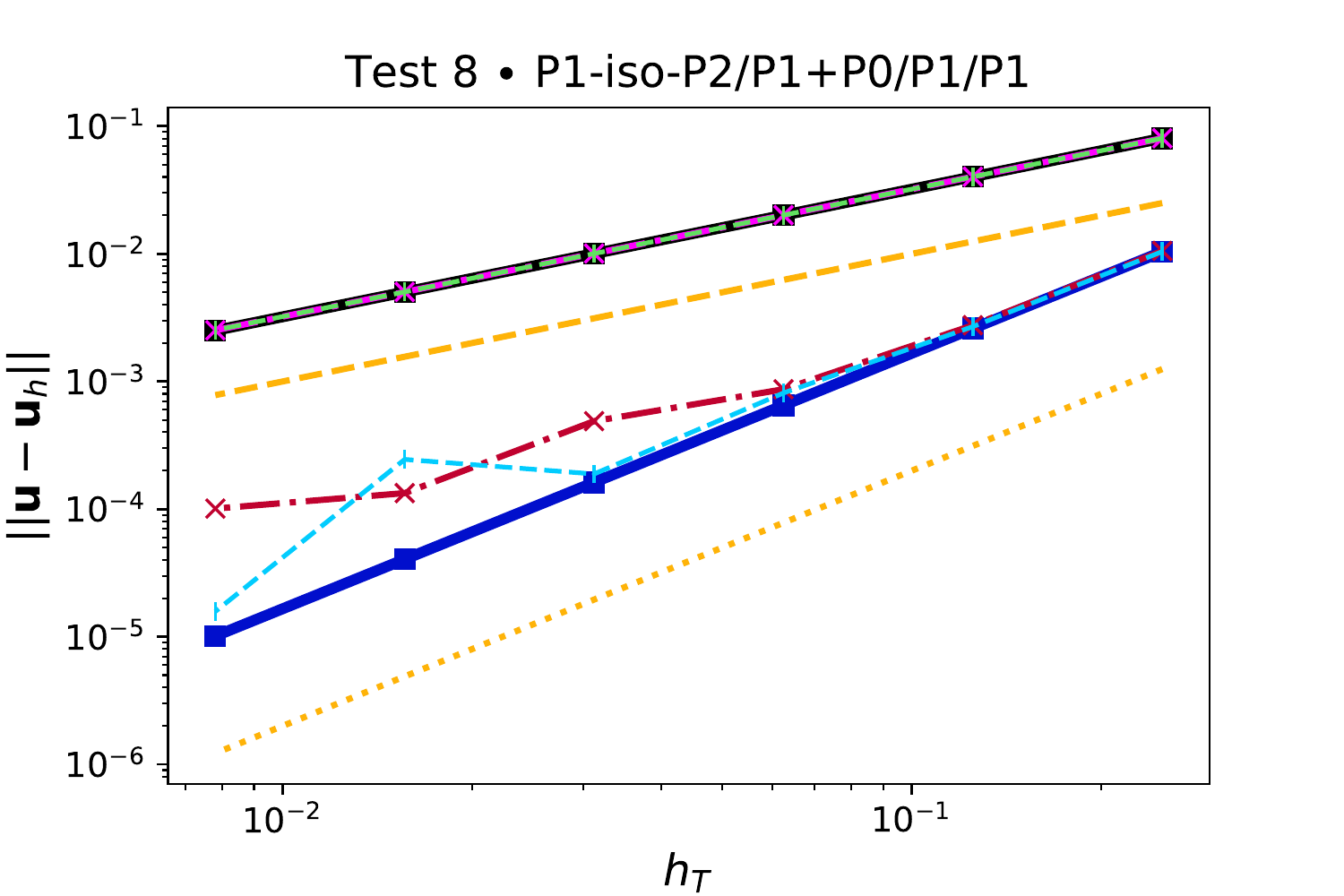}
		\includegraphics[width=0.24\linewidth]{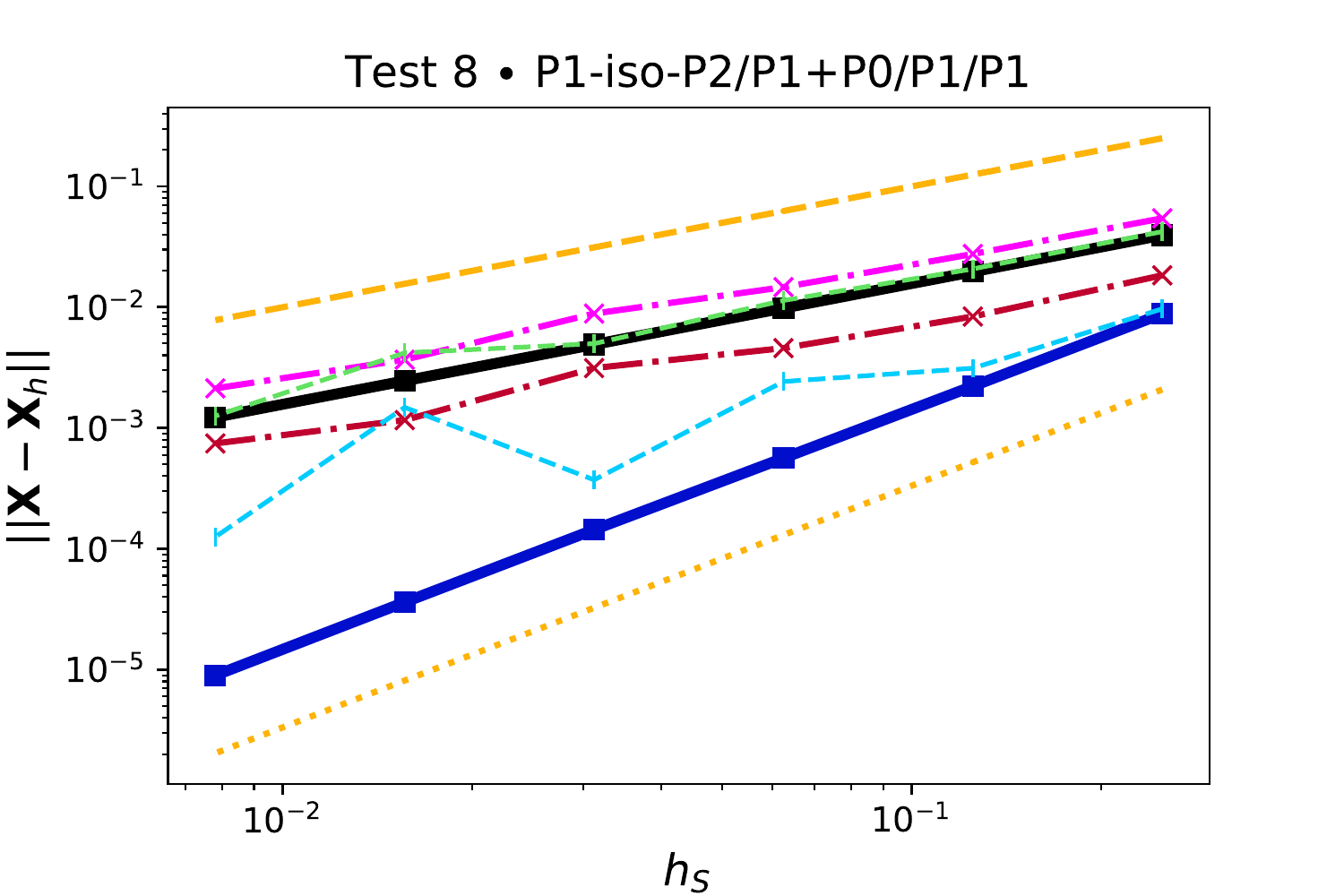}
		\includegraphics[width=0.24\linewidth]{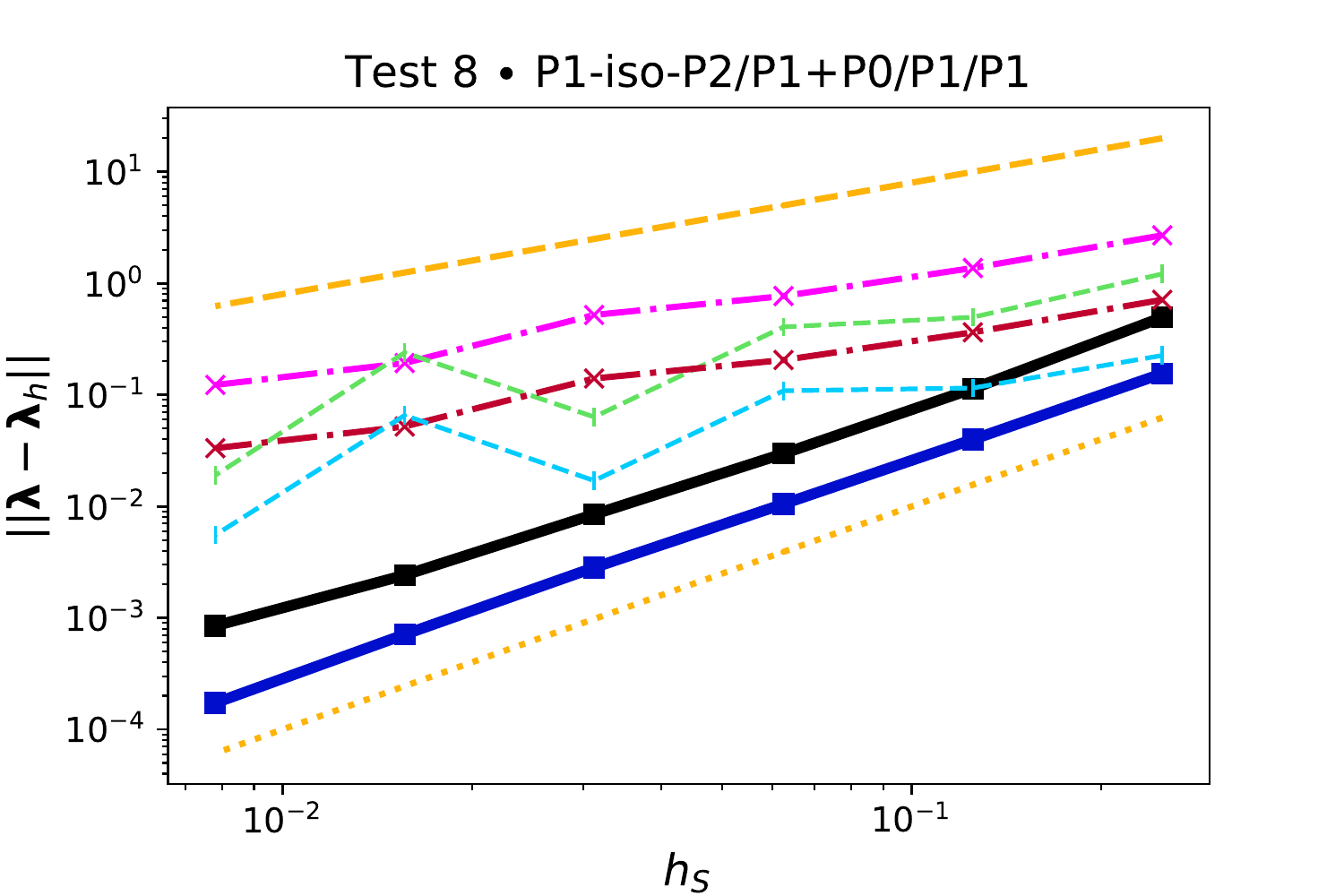}\\
		\includegraphics[trim=150 16 50 50,width=1.4\linewidth]{figures/label_p0-eps-converted-to}
		\caption{Convergence plots of Test 8 with $\Pcal_1-iso-\Pcal_2/\Pcal_1+\Pcal_0/\Pcal_1/\Pcal_1$}
		\label{fig:test8_p1p0}
	\end{figure}

	The results, shown in Figures \ref{fig:test8_p1} and \ref{fig:test8_p1p0}, are in agreement with Test 7.
	
	\section{Conclusions}
	In the opening sections of this paper, we recalled the formulation for fluid-structure interaction problem introduced in \cite{2015}: this approach, born in the spirit of the fictitious domain approach, is based on the use of a distributed Lagrange multiplier with the aim of enforcing the motion condition.
	
	In particular, the main contribution of this paper is given by the comparison of two approaches to be used for the assembling of the interface matrix between fluid and structure: the first approach is based on the computation of the intersection between the solid and the fluid meshes, while the second one consists in computing the entries of the matrix with a direct integration on each solid element. Several numerical tests showed that the optimal convergence rate is reached only when the intersection is computed. On the other hand, when the intersection of the two meshes is not computed, the method is not fully performing and the increase of precision of the quadrature rule does not produce any significant improvement.

	\bibliographystyle{plain}
	\bibliography{biblio}
	\nocite{2011}
	\nocite{2015}
	\nocite{mixedFEM}
	
\end{document}